\documentclass[12pt,reqno]{amsart}
\usepackage[margin=1.25in]{geometry}
\usepackage[utf8]{inputenc}
\usepackage{amsmath, amssymb, amsthm}
\usepackage{mathtools}
\usepackage{anyfontsize}
\usepackage{lmodern}
\usepackage{microtype}
\usepackage{enumitem}
\usepackage{ifthen}
\usepackage{environ}
\usepackage{xfrac}
\usepackage{pdflscape}
\usepackage{esvect}
\usepackage{appendix}

\usepackage{bbm}
\usepackage{bm}
\usepackage{makecell}
\usepackage{tikz}
    \usetikzlibrary{calc}
    \usetikzlibrary{knots}
    \usetikzlibrary{math}
    \usetikzlibrary{shapes}
    \usetikzlibrary{arrows}
    \usetikzlibrary{cd}
    \usetikzlibrary{intersections}
\usepackage{xcolor} % look into colorblind palettes
    \colorlet{DarkGreen}{green!50!black}
    \colorlet{DarkRed}{red!90!black}
    \colorlet{DarkBlue}{blue!90!black}

\usepackage[pdfencoding=unicode,pdfusetitle]{hyperref}
    \hypersetup{colorlinks=true,
        linkcolor=blue,
        filecolor=purple,      
        urlcolor=[rgb]{0 0 .6},
        psdextra}

\renewcommand{\arraystretch}{1.5}

\newcommand{\id}{\textsf{id}}
\newcommand{\1}{\mathbbm{1}}

\newcommand{\s}{\mathcal}
\newcommand{\bb}{\mathbb}

\DeclareMathOperator{\Hom}{Hom}
\DeclareMathOperator{\End}{End}
\DeclareMathOperator{\Aut}{Aut}

\DeclareMathOperator{\Fun}{Fun}

\DeclareMathOperator{\Rep}{Rep}

\DeclareMathOperator{\op}{op}
\DeclareMathOperator{\Vect}{Vect}
\DeclareMathOperator{\fd}{fd}
\DeclareMathOperator{\Gal}{Gal}
\DeclareMathOperator{\sgn}{sgn}

\DeclareMathOperator{\QF}{QF}
\DeclareMathOperator{\GO}{GO}

\makeatletter
\newtheorem*{rep@theorem}{\rep@title}
\newcommand{\newreptheorem}[2]{%
\newenvironment{rep#1}[1]{%
 \def\rep@title{#2 \ref{##1}}%
 \begin{rep@theorem}}%
 {\end{rep@theorem}}}
\makeatother

\theoremstyle{definition} % changed for now since italics are hard to read
\newtheorem{theorem}{Theorem}[section]
\newreptheorem{theorem}{Theorem}
\newtheorem{proposition}[theorem]{Proposition}
\newtheorem{corollary}[theorem]{Corollary}
\newtheorem{lemma}[theorem]{Lemma}
\theoremstyle{definition}
\newtheorem{definition}[theorem]{Definition}
\newtheorem{example}[theorem]{Example}

\newtheorem{remark}[theorem]{Remark}

\newtheorem{notation}[theorem]{Notation}

\newtheorem{fact}[theorem]{Fact}

\NewEnviron{tikzineqn}[1][]{\vcenter{\hbox{\begin{tikzpicture}[#1]
            \BODY \end{tikzpicture}}}}
            
\newcommand{\arxiv}[1]{\href{http://arxiv.org/abs/#1}{\tt arXiv:\nolinkurl{#1}}}
\newcommand{\arXiv}[1]{\href{http://arxiv.org/abs/#1}{\tt arXiv:\nolinkurl{#1}}}
\newcommand{\doi}[1]{\href{http://dx.doi.org/#1}{{\tt DOI:#1}}}

\newcommand{\mathscinet}[1]{\href{http://www.ams.org/mathscinet-getitem?mr=#1}{\tt #1}}
\newcommand{\googlebooks}[1]{(preview at \href{http://books.google.com/books?id=#1}{google books})}

\newcommand{\eqnscale}{0.4}
\newcommand{\tscale}{0.8}

\title{Braidings for Non-Split Tambara-Yamagami Categories over the Reals}
\author{David Green${}^{1, \dagger}$}\thanks{${}^{1}$\texttt{Department of Mathematics, The Ohio State University, Columbus OH, USA}}
\author{Yoyo Jiang${}^{2}$}\thanks{${}^{2}$\texttt{Department of Mathematics, Johns Hopkins University, Baltimore MD, USA}}
\author{Sean Sanford${}^{3}$}\thanks{${}^{3}$\texttt{School of Mathematics, The University of Edinburgh, Edinburgh Scotland, UK}}\thanks{\texttt{${}^\dagger$(Corresponding Author) Email:\href{mailto:green.2116@buckeyemail.osu.edu}{green.2116@buckeyemail.osu.edu}}}
\begin{document}

\begin{abstract}
Non-split Real Tambara-Yamagami categories are a family of fusion categories over the real numbers that were recently introduced and classified by Plavnik, Sanford, and Sconce.  We consider which of these categories admit braidings, and classify the resulting braided equivalence classes. We also prove some new results about the split real and split complex Tambara-Yamagami Categories. 
\end{abstract}
\maketitle

\tikzmath{
    \x=1;
    \topratio=2/3;
    \beadsizenum=\x/2;
}
\def\beadsize{\beadsizenum cm}

\tikzstyle{strand x} = [thick,DarkRed]
\tikzstyle{strand y} = [thick,DarkGreen]
\tikzstyle{strand z} = [thick,orange]
\tikzstyle{strand w} = [thick,violet!80]
\tikzstyle{strand a} = [thick,DarkRed]
\tikzstyle{strand 1} = [thick,DarkRed,dotted]
\tikzstyle{strand b} = [thick,DarkGreen]
\tikzstyle{strand c} = [thick,orange]
\tikzstyle{strand ab} = [thick,orange]
\tikzstyle{strand bc} = [thick,orange]
\tikzstyle{strand abc} = [thick,DarkBrown]
\tikzstyle{strand m} = [thick,black]
\tikzstyle{node a}   = [DarkRed]
\tikzstyle{node b}   = [DarkGreen]
\tikzstyle{node c}   = [orange]
\tikzstyle{node x}   = [DarkRed]
\tikzstyle{node y}   = [DarkGreen]
\tikzstyle{node z}   = [orange]
\tikzstyle{node w}   = [violet!80]
\tikzstyle{node ab} = [orange]
\tikzstyle{node bc} = [orange]
\tikzstyle{node abc} = [DarkBrown]
\tikzstyle{node m}   = [black]
\tikzstyle{beadLabel}   = [font=\tiny]
\tikzstyle{smallbead} = [circle,
                         fill=blue!20,
                         draw=black,
                         inner sep=0,
                         minimum size=\beadsize*0.7,
                         font=\tiny]
\tikzstyle{bead}     = [circle,
                        fill=blue!20,
                        draw=black,
                        inner sep=0,
                        minimum size=\beadsize,
                        font=\tiny]
\tikzstyle{emptybead} = [circle,
                         thick,
                         fill=blue!40,
                         draw=black,
                         inner sep=0,
                         minimum size=\beadsize*0.5,
                         font=\tiny]
\tikzstyle{longbead} = [rectangle,
                        fill=blue!20,
                        rounded corners=2mm,
                        draw=black,
                        inner sep=1mm,
                        minimum size=\beadsize,
                        font=\tiny]
% \tikzstyle{grid} = [black!20] % draw grid
\tikzstyle{grid} = [draw=none] % do not draw grid

\newcommand{\TrivalentVertex}[3]{
    \coordinate (mid)            at (0,0);
    \coordinate (top)            at (0,1);
    \coordinate (bottom left)   at (-1,-1);
    \coordinate (bottom right)    at (1,-1);
    \draw[strand #1] (mid) to (bottom left) node[below left] {$#1$};
    \draw[strand #2] (mid) to (bottom right) node[below right] {$#2$};
    \draw[strand #3] (mid) to (top) node[above] {$#3$};
    }

\newcommand{\DagTrivalentVertex}[3]{
    \coordinate (mid)            at (0,0);
    \coordinate (bot)            at (0,-1);
    \coordinate (top left)   at (-1,1);
    \coordinate (top right)    at (1,1);
    \draw[strand #1] (mid) to (top left) node[above left] {$#1$};
    \draw[strand #2] (mid) to (top right) node[above right] {$#2$};
    \draw[strand #3] (mid) to (bot) node[below] {$#3$};
    }

\newcommand{\TetraTransformBeads}[7]{
    % grid lines for testing
    % \draw[draw=black!20] (-\x,-\x) grid (\x,\x);
    % coordinates - do not change
    \coordinate (mid)            at (0,0);
    \coordinate (top)            at (0,\topratio*\x);
    \coordinate (bottom left)    at (-\x,-\x);
    \coordinate (bottom right)   at (\x,-\x);
    \coordinate (bottom mid)     at (0,-\x);
    \coordinate (right  vertex)  at ($1/2*(bottom right)$);
    \coordinate (left   vertex)  at ($1/2*(bottom left)$);
    \draw[strand #2] (mid)           to (top);
    \draw[strand #3] (mid)           to (left  vertex);
    \draw[strand #4] (mid)           to (right vertex);
    \draw[strand #5] (left  vertex)  to (bottom left);
    \draw[strand #7] (right vertex)  to (bottom right);
    \ifthenelse{%if
        \equal{#1}{left}}
        {%then
        \draw[strand #6] (left  vertex) to (bottom mid);
        }{%else
        \draw[strand #6] (right vertex) to (bottom mid);
        }
    \node[node #2][above]  at (top)           {$#2$};
    \node[node #5][below]  at (bottom left)   {$#5$};
    \node[node #6][below]  at (bottom mid)    {$#6$};
    \node[node #7][below]  at (bottom right)  {$#7$};
    \ifthenelse{%if
        \equal{#1}{left}}
        {%then
        \node[node #3][above left]
            at ($(0,0)!1/2!(left vertex)$)  {$#3$};
        }{%else
        \node[node #4][above right]
            at ($(0,0)!1/2!(right vertex)$) {$#4$};
        }
}

\newcommand{\TetraTransform}[7]{
\begin{tikzineqn}
    % grid lines for testing
    % \draw[draw=black!20] (-\x,-\x) grid (\x,\x);
    % coordinates - do not change
    \coordinate (mid)            at (0,0);
    \coordinate (top)            at (0,\topratio*\x);
    \coordinate (bottom left)    at (-\x,-\x);
    \coordinate (bottom right)   at (\x,-\x);
    \coordinate (bottom mid)     at (0,-\x);
    \coordinate (right  vertex)  at ($1/2*(bottom right)$);
    \coordinate (left   vertex)  at ($1/2*(bottom left)$);
    \draw[strand #2] (mid)           to (top);
    \draw[strand #3] (mid)           to (left  vertex);
    \draw[strand #4] (mid)           to (right vertex);
    \draw[strand #5] (left  vertex)  to (bottom left);
    \draw[strand #7] (right vertex)  to (bottom right);
    \ifthenelse{%if
        \equal{#1}{left}}
        {%then
        \draw[strand #6] (left  vertex) to (bottom mid);
        }{%else
        \draw[strand #6] (right vertex) to (bottom mid);
        }
    \node[node #2][above]  at (top)           {$#2$};
    \node[node #5][below]  at (bottom left)   {$#5$};
    \node[node #6][below]  at (bottom mid)    {$#6$};
    \node[node #7][below]  at (bottom right)  {$#7$};
    \ifthenelse{%if
        \equal{#1}{left}}
        {%then
        \node[node #3][above left]
            at ($(0,0)!1/2!(left vertex)$)  {$#3$};
        }{%else
        \node[node #4][above right]
            at ($(0,0)!1/2!(right vertex)$) {$#4$};
        }
\end{tikzineqn}
}

\newcommand{\DrawBead}[4][]{
\node[bead,#1] at ($(#2)!1/2!(#3)$) {$#4$};
}

\newcommand{\DrawSmallBead}[4][]{
\node[smallbead,#1] at ($(#2)!1/2!(#3)$) {$#4$};
}

\newcommand{\DrawLongBead}[4][]{
\node[longbead,#1] at ($(#2)!1/2!(#3)$) {$#4$};
}
%to be used inside Tikzinequation environments
\newcommand{\AMBraidCrossing}{\begin{knot}[clip width=10]
		\strand[strand a] (-1,-1) node[below] {$a$} to (1,1);
		\strand[strand m] (1,-1) node[below] {$m$} to (-1,1);
\end{knot}}
\newcommand{\MABraidCrossing}{\begin{knot}[clip width=10]
		\strand[strand m] (-1,-1) node[below] {$m$} to (1,1);
		\strand[strand a] (1,-1) node[below] {$a$} to (-1,1);
\end{knot}}

\newcommand\blfootnote[1]{
	\begingroup
	\renewcommand\thefootnote{}\footnote{#1}
	\addtocounter{footnote}{-1}
	\endgroup
}
\section{Introduction}
In \cite{pss23}, Plavnik, Sconce and our third author introduced and classified three infinite families of fusion categories over the real numbers.
These categories are analogues of the classical Tambara-Yamagami fusion categories introduced and classified in \cite{ty98}.
This new version of Tambara-Yamagami (TY) categories allowed for non-split simple objects: simples whose endomorphism algebras are division algebras, and not just $\mathbb R$.
These non-split TY categories generalize classical examples such as $\Rep_{\mathbb R}(Q_8)$ and $\Rep_{\mathbb R}(\mathbb Z/4\mathbb Z)$, but also include many new fusion categories that fail to admit a fiber functor, i.e. they are not even $\Rep(H)$ for a semisimple Hopf-algebra.
This paper provides a classification of all possible braidings that exist on these new non-split TY categories.

Since their introduction, TY categories have been studied and generalized extensively (including the closely related notion of \textit{near-group} categories) \cite{Tambara2000,MR2677836,Izumi_2021,GALINDO_2022,SchopierayNonDegenExtension,galindo2024modular}.  Their complexity lies just above the pointed fusion categories, and well below that of general fusion categories.
This intermediate complexity allows for deep analysis of their structure, while simultaneously providing examples of interesting properties that cannot be observed in the more simplistic pointed categories.
For example, in \cite{Nikshych2007NongrouptheoreticalSH} Nikshych showed that some TY categories provide examples of non-group-theoretical (i.e. not even Morita equivalent to pointed) fusion categories that admit fiber functors.

The physical motivation for extending this theory of TY categories to the real numbers comes from time reversal symmetry.
Though there is some debate (see \cite{MR3675715}) about the exact nature of time reversal in quantum mechanics, it is generally agreed upon that time reversal acts by an anti-unitary operator (see for example \cite[Section 1.1]{10.21468/SciPostPhys.5.1.006}), so in particular, it performs complex conjugation on scalars.
Thus, without assuming any properties beyond anti-unitarity, we say that a time reversal symmetry on a fusion category $\mathcal C$ over $\mathbb C$ is a categorical action of $\mathbb Z/2\mathbb Z$ by $\mathbb R$-linear monoidal functors on $\mathcal C$, that behaves as complex conjugation on $\End(\1)$.
By Etingof and Gelaki's theory of categorical descent \cite{MR2946231}, real fusion categories must arise as the equivariantization $\mathcal C^{\mathbb Z/2\mathbb Z}$ of $\mathcal C$ with respect to such a time reversal action.
In condensed matter terminology, fusion categories describe the topological quantum field theory (TQFT) that arises in the low-energy limit of a gapped quantum field theory in (1+1)D (see \cite[Section I]{MR4268163} in TQFT language, and \cite[Section X B]{kong2014braidedfusioncategoriesgravitational} in the language of topological orders). 
Thus real fusion categories describe time reversal symmetric TQFTs in (1+1)D.
In the (2+1)D setting, time reversal symmetric TQFTs should be described by \emph{braided} fusion categories over the reals.

With an eye toward time reversal symmetry in (2+1)D, in this paper we classify all possible braidings admitted by non-split TY categories over $\mathbb R$.
We proceed in the style of Siehler \cite{sie00}, by distilling invariants of a braiding that follow from the hexagon equations.
Next, we leverage the description of monoidal equivalences given in \cite{pss23} in order to determine which braiding invariants produce braided equivalent categories, thus establishing a classification.
Along the way we describe all braided classifications for split real and split complex TY categories as well.

In Section \ref{sec:CrossedBraided}, we explain why complex/complex (see section for terminology) TY categories can never admit a braiding, due to the presence of Galois-nontrivial objects.
\subsection{Results} 

Tambara-Yamagami categories over the reals are roughly classified by four pieces of data: The division algebras $\End(1)$ and $\End(m)$ of endomorphisms of the unit object $1$ and the non-invertible object respectively, a nondegenerate bilinear form $\chi$ on the group of invertible objects $A$, and a square root $\tau$ of $|A|$. When  $\End(m)$ is the field of complex numbers, the additional data of a Galois automorphism of $\mathbb C / \mathbb R$ is required to fully specify the monoidal equivalence class.   

For all the split and non-split real Tambara-Yamagami categories over $\mathbb R$, there turns out to be a unique family of bicharacters $\chi$ such that the associated Tambara-Yamagami category can possibly admit a braiding.
As has appeared previously in the literature, the classification is in terms of $\Aut(A, \chi)$ orbits of \textit{$\chi$-admissible forms}, these are quadratic forms with coboundary $\chi$. The results are summarized below in Table \ref{table:Table1}, under the assumption that the group of invertible objects is not trivial (see the theorem statements for precise results in these cases).
We label the cases by
``$\End(1)$/$\End(m)$"
(for instance, ``real/complex"
or ``$\mathbb R$/$\mathbb C$"
when $\End(1)\cong\mathbb R$
and $\End(m)\cong\mathbb C$)
as well as the chosen Galois
automorphism in the
$\mathbb R$/$\mathbb C$ case.

The other piece of data appearing in Table \ref{table:Table1} is the number $\sigma_3(1)$, defined as the element of $\End(m)$ satisfying 

\[
        \begin{tikzineqn}[scale=\eqnscale]
            \draw[strand a, dotted] (0,0) to ++(0,1) node[above] {$1$};
            \begin{knot}[clip width=10]
                \strand[strand m] (0,0)
                to ++(1,-1)
                to ++(-2,-2) node[below left] {$m$};
                \strand[strand m] (0,0)
                to ++(-1,-1)
                to ++(2,-2) node[below right] {$m$};
            \end{knot}
        \end{tikzineqn}
        := \
        \begin{tikzineqn}[scale=\tscale] 
            \TrivalentVertex{m}{m}{1} 
            \DrawLongBead{mid}{bottom right}{\sigma_3(1)}
        \end{tikzineqn}. 
\]

This element is frequently, but not always, an invariant of the braided structure. 
 
\begin{center}
\begin{table}[!h]
	\caption{Number of braidings on Tambara-Yamagami Categories}
	\label{table:Table1}
	\begin{tabular}{|c|c|c|c|c|c|}
		\hline
		Case & Split Real & $\mathbb{R} / \mathbb{C}, \id$ & $\mathbb{R} / \mathbb{C}, \bar \cdot $ & $\mathbb{R} / \mathbb{H}$  \\ \hline
        Article section & \ref{sec:SplitReal} & \ref{sec:Real/Complex} & \ref{sec:Real/Complex} & \ref{sec:RealQuaternionic} \\
		\hline 
		$\chi$-admissible orbits & 2 & 2 & 2 & 2  \\ \hline 
		Orbits extending to braidings & 1 & 2 & 2 & 1  \\ \hline 
		Braidings per orbit & 2 & Varies & 2 & 2   \\ \hline 
		Total braidings & 2  & 3 & 4 & 2   \\ \hline
		Is $\tau$ an invariant? & Yes & No & Yes & Yes  \\  \hline 
		Is $\sigma_3(1)$ an invariant? & Yes & No & Yes & Yes  \\ \hline 
	\end{tabular}
\end{table}
\end{center}
 In contrast to the real case, there are three families of bicharacters (not all of which are defined on a given 2-group) on the split complex Tambara-Yamagami categories. If we let $\ell$ be the quadratic form on $\mathbb{Z}/2\mathbb{Z}$ with $\ell(g,g) = -1$, then the three families of bicharacters are distinguished by the multiplicity (modulo 3) of $\ell$ in their direct sum decomposition. We write $|\ell|$ for this number. In this case all orbits of quadratic forms extend to braidings. The results (derived in Section \ref{sec:SplitComplex}) are summarized below, under the assumption that the group of invertibles is not too small (see the theorem statements for precise results in these cases).
\begin{center}
	\begin{table}[!h]
		\caption{Braidings on split complex Tambara-Yamagami categories}
		\label{table:Table2}
	\begin{tabular}{|c|c|c|c|}
	\hline
	$|\ell|$ & 0 & 1 & 2 \\
	\hline 
	$\chi$-admissible orbits & 2 & 4 & 4 \\ \hline 
	Braidings per orbit & 2 & 2 & 2 \\ \hline 
	Total braidings & 4  & 8 & 8 \\ \hline
\end{tabular} 
\end{table}
\end{center}
Here $\tau$ and $\sigma_3(1)$ are always invariants, and the classification is up to \textit{complex}-linear functors.

Next, we collect a table describing when the various braidings we define are symmetric or non-degenerate. If we let $\sigma$ be a quadratic form inducing $\chi$ (see Section \ref{sec:QFAnalysis}), its \textit{sign} $\sgn(\sigma)$ is defined in terms of the Gauss sum (Notation \ref{not:QF}). In some cases, this sign must match the sign of $\tau$ (i.e, if $\tau$ is a positive square root of $|A|$ then $\sigma$ must have positive sign) in order for the braiding to exist at all. In the other cases, the difference between the signs of $\sigma$ and $\tau$ controls whether or not the categories under consideration are symmetric or nondegenerate. Nondegeneracy occurs only when the group of invertible objects is very small.
	\begin{center}
		\begin{table}[!h]
			\caption{Symmetry and nondegeneracy of the braidings}
			\label{table:Table3}
    \begin{tabular}{|c|c|c|} \hline
    Case & Symmetric? & Nondegenerate? \\ \hline
    Split Real & Always & Never \\ \hline 
    Real/Quaternionic & Always & Never \\  \hline 
    \makecell{Real/Complex, $g = \id_\mathbb{C},$ \\ $\sgn(\sigma) = \sgn(\tau)$ }& Never & Never \\ \hline 
        \makecell{Real/Complex, $g = \id_\mathbb{C},$ \\ $\sgn(\sigma) = -\sgn(\tau)$ }& Never & Only when $A = *$ \\ \hline 
    Real/Complex, $g = \bar \cdot$ & Always & Never \\ \hline
    Split Complex, $|\ell| = 0$ & Only when $\sgn(\sigma) = \sgn(\tau)$ & \makecell{Only when $A = *$ and  \\$\sgn(\sigma) = -\sgn(\tau)$} \\ \hline
    Split Complex, $|\ell| = 1$ & Never & Never \\ \hline
    Split Complex, $|\ell| = 2$ & Never & Never \\ \hline
    \end{tabular}
\end{table}
\end{center}
Some cases include multiple equivalence classes of braidings, but in all cases, the results in the table above are immediate from the classifications of braidings we give.
The nondegenerate split complex categories are the well-known semion and reverse semion categories respectively. Finally, in Section \ref{sec:CrossedBraided}, we show that in the complex/complex case, when $m$ is \textit{Galois nontrivial}, there are no braidings.  

\subsection{Acknowledgments} This project began during Summer 2023 as part of the Research Opportunities in Mathematics for Underrepresented Students, supported by NSF
grants DMS CAREER 1654159 and DMS 2154389. DG would like to thank the Isaac Newton Institute for Mathematical Sciences, Cambridge, for support and hospitality during the \textit{Topology, Representation theory and Higher Structures} programme where work on this paper was undertaken. This work was supported by EPSRC grant no EP/R014604/1. YJ was supported by the Woodrow Wilson Research Fellowship at Johns Hopkins University. Rui Wen provided critical corrections to a previous edition of this manuscript. 
DG, SS, and YJ would all also like to thank David Penneys for his guidance and support.  

\section{Background}

We refer the reader to \cite{EGNO15} for
the basic theory of fusion categories and to
\cite{pss23} and \cite{MR4806973} for the basics of
(non-split) fusion categories over non-algebraically closed fields.

\begin{definition}\label{defn:BraidedMonodialCategory}
A braiding on a monoidal category $\mathcal{C}$ is a set of isomorphisms
\[
\{\beta_{x,y}\colon x\otimes y \xrightarrow{} y\otimes x\}_{x,y\in \mathcal{C}}
\]
such that the following diagrams commute (omitting $\otimes$)

\begin{equation}\begin{tikzcd}\label{defn:ForwardsHexagon}
	& {x(yz)} & {(yz)x} \\
	{(xy)z} &&& {y(zx)} \\
	& {(yx)z} & {y(xz)}
	\arrow["\alpha_{x,y,z}", from=2-1, to=1-2]
	\arrow["{\beta_{x,yz}}", from=1-2, to=1-3]
	\arrow["\alpha_{y,z,x}", from=1-3, to=2-4]
	\arrow["{\beta_{x,y}\otimes \id_z}"', from=2-1, to=3-2]
	\arrow["\alpha_{y,x,z}"', from=3-2, to=3-3]
	\arrow["{\id_y \otimes \beta_{x,z}}"', from=3-3, to=2-4]
\end{tikzcd}\end{equation}

\begin{equation}\begin{tikzcd}\label{defn:BackwardsHexagon}
	& {(xy)z} & {z(xy)} \\
	{x(yz)} &&& {(zx)y} \\
	& {x(zy)} & {(xz)y}
	\arrow["\alpha^{-1}_{x,y,z}", from=2-1, to=1-2]
	\arrow["{\beta_{xy,z}}", from=1-2, to=1-3]
	\arrow["\alpha^{-1}_{z,x,y}", from=1-3, to=2-4]
    \arrow["{\id_x \otimes \beta_{y,z}}"', from=2-1, to=3-2]
	\arrow["\alpha^{-1}_{x,z,y}"', from=3-2, to=3-3]
    \arrow["{\beta_{x,z}\otimes \id_y}"', from=3-3, to=2-4]
\end{tikzcd}\end{equation}

for all objects $x,y,z\in \mathcal{C}$, where $\alpha_{x,y,z}$ is the associator.
We will refer to the commutativity of the top diagram as the hexagon axiom
and of the bottom diagram as the inverse hexagon axiom.
Note that these encode commutative diagrams of natural transformations.
\end{definition}
Our goal is to classify braiding structures on a fusion category $\mathcal{C}$ with a fixed monoidal structure.
To do this, we will use the Yoneda lemma to show that
the data defining abstract braiding isomorphisms
is given by a finite set of linear maps between Hom-spaces,
which we can then specify by their values on basis vectors.

Specifically, a braiding on $\mathcal{C}$ is given by a natural transformation
\[
    \beta\colon (-)\otimes (=) \Rightarrow
    (=)\otimes (-),
\]
a morphism in the category
of linear functors from $\mathcal{C}\times \mathcal{C}\to \mathcal{C}$.
By semisimplicity, it suffices to consider the components of $\beta$ on simple objects, and by the Yoneda lemma,
this data is given by a natural transformation in $\Fun(\mathcal{S}_{\mathcal{C}}^{\op}\times \mathcal{S}_{\mathcal{C}}^{op}\times \mathcal{S}_{\mathcal{C}}, \Vect_k^{\fd})$, i.e. a finite set of linear maps
\[
\Hom_{\mathcal{C}}(s\otimes t,u)\xrightarrow[]{\beta_{t,s}^{*}}
\Hom_{\mathcal{C}}(t\otimes s,u)
\]
natural in simple objects $s,t,u\in \mathcal{C}$.
Furthermore, by Schur's lemma, it suffices to
check naturality on endomorphisms of $s$, $t$ 
and $u$,
which is in particular vacuous if the category is split.
After fixing a set of basis vectors for
the Hom spaces, this reduces to a set of matrix coefficients,
which we will refer to as the braiding coefficients.

Similarly, to check that $\beta$ satisfies
the hexagon axioms, it suffices to check that
for any $s,t,u,v\in \mathcal{C}$ simple,
the two linear maps
\[
\Hom_\mathcal{C}(t(us),v)\xrightarrow[]{}
\Hom_\mathcal{C}((st)u,v)
\]
obtained by precomposing the top and bottom
paths of \eqref{defn:ForwardsHexagon} are equal,
and similarly for the inverse hexagon axiom.
With the choice of a basis for Hom spaces,
this condition is given by the set of polynomial equations in terms in the braiding coefficients,
which we will refer to as the braiding equations.

\section{Quadratic forms on elementary abelian 2-groups}\label{sec:QFAnalysis}

In this section we gather several facts that we will need regarding quadratic forms.
We follow the convention of the fusion category literature (see for example \cite{MR2609644}, \cite{MR3039775}, and \cite{MR3022755}) in writing the these forms multiplicatively.
The full classification of quadratic forms on finite groups was done in \cite{wall63}.
We will need this level of generality when discussing forms with values in $\mathbb C^\times$.

When the forms take values in $\mathbb R^\times$, the theory reduces to that of quadratic forms on vector spaces over $\mathbb F_2$.
The theory of quadratic forms on vector spaces is much older, and very well understood (see for example \cite{MR770063}, \cite{MR310083}, \cite{MR506372} for theorems, or \cite{MR1803370} for a historical account). The account here is intended to be self contained and convenient for our purposes.
We have supplied proofs of our own, when necessary, but make no claims of originality of the results in this section.

Given a field $\mathbb K$, a quadratic form on a finite abelian group $A$ is a function $\sigma:A\to\mathbb K^\times$ such that $\sigma(x^{-1})=\sigma(x)$, and 
\[(\delta\sigma)(a,b)\,:=\frac{\sigma(ab)}{\sigma(a)\sigma(b)}\]
is a bicharacter.
When equipped with a quadratic form $\sigma$, the pair $(A,\sigma)$ is called a pre-metric group, and is called a metric group in the case where $\delta\sigma$ is nondegenerate.

Pointed braided fusion categories 
$(\mathcal C,\{\beta_{X,Y}\}_{X,Y})$ over $\mathbb K$
(recall that a fusion category is pointed if all its simple objects are invertible)
are determined up to equivalence by their group of invertible objects $\mathrm{Inv}(\mathcal C)$ and the quadratic form $\sigma:\mathrm{Inv}(\mathcal C)\to\mathbb K^\times$ given by the formula
\[\beta_{g,g}=\sigma(g)\cdot\id_{g^2}\,.\]
In fact, this classification arises from an equivalence of categories, and is due to Joyal and Street in \cite[Section 3]{MR1250465} (their terminology differs from ours).
This equivalence of categories implies that two pointed braided fusion categories are equivalent if and only if their corresponding pre-metric groups are isometric.

Any braided TY category contains a pointed braided subcategory, and thus gives rise to a pre-metric group.
Our analysis in the non-split TY cases will mirror that of the split cases, and it is interesting to note that the quadratic form that gives rise to a braiding on a TY category is a square root of the quadratic form on its own pointed subcategory.

\begin{definition}\label{defn:ChiAdmissibleFunction}
    Given a bicharacter $\chi:A\times A\to\mathbb K^\times$, a quadratic form $\sigma:A\to\mathbb K^\times$ is said to be $\chi$-admissible if $\delta\sigma\,=\,\chi$.
    The collection of all $\chi$-admissible quadratic forms will be denoted $\QF_{\mathbb K}(\chi)$.  For the majority of the paper, we are concerned with $\QF_{\mathbb R}(\chi)$, and so we simply write $\QF(\chi)$ when $\mathbb K=\mathbb R$.
\end{definition}

\begin{remark}
    In the literature the coboundary $\delta\sigma$ is often referred to as the associated bicharacter of the quadratic form $\sigma$ (see e.g. \cite[Section 2.11.1]{MR2609644}).  Thus ``$\sigma$ is $\chi$-admissible'' is synonymous with ``the associated bicharacter of $\sigma$ is $\chi$''.

    We caution that our coboundary is inverted in order to align with the hexagon equations that appear later, though this is immaterial from a formal standpoint.
    Furthermore, in some conventions the phrase ``associated bicharacter'' or ``associated bilinear form'' refers to the square root of $\delta\sigma$ (see e.g. \cite[Section 7]{wall63}).
    Our general feeling is that while this square root is irrelevant for odd groups, it complicates the analysis unnecessarily for 2-groups, which are the main application in this paper.
\end{remark}

The group $\Aut(A, \chi)$ of automorphisms preserving the bicharacter acts on $\QF(\chi)$ by the formula $(f.\sigma)(g):=\sigma\big(f^{-1}(a)\big)$.
We will be particularly concerned with the Klein four-group $K_4:=(\mathbb Z/2\mathbb Z)^2$ and powers $(\mathbb Z/2\mathbb Z)^n$ generally.  We will occasionally think of $(\mathbb Z/2\mathbb Z)^n$ as an $\mathbb F_2$ vector space in order to refer to a basis, but we will still write the group multiplicatively.

A key tool in the analysis of quadratic forms is the Gauss sum.

\begin{definition}
    Given a quadratic form $\sigma:A\to\mathbb K^\times$, the Gauss sum $\Sigma(\sigma)\in\mathbb K$ of $\sigma$ is the sum $\Sigma_{a\in A}\sigma(a)$.
    Occasionally we will write this as $\Sigma(A)$, when the quadratic form can be inferred.
\end{definition}

Recall that a subgroup $H\leq A$ is said to be \emph{isotropic} if $\sigma|_H=1$.
Isotropic subgroups automatically satisfy $H\leq H^\perp$, where $H^\perp$ is the orthogonal compliment of $H$ with respect to $\delta\sigma$.
A metric group $(A,\sigma)$ is said to be \emph{anisotropic} if $\sigma(x)=1$ implies $x=1$.
An isotropic subgroup is said to be \emph{Lagrangian} if $H=H^\perp$, and a pre-metric group is said to be \emph{hyperbolic} if it contains a Lagrangian subgroup. The following lemma records some important properties of Gauss sums with respect to isotropic subgroups.

\begin{lemma}[{\cite[cf. Section 6.1]{MR2609644}}]\label{lem:GaussSumProperties}
    Let $(A,\sigma)$ be a pre-metric group.
    \begin{enumerate}[label=(\roman*)]
        \item For any isotropic subgroup $H\leq A$, $\Sigma(A)=|H|\cdot\Sigma(H^\perp/H)$.
        \item If $A$ is hyperbolic, then $\Sigma(A)$ is a positive integer.
        \item If $\Sigma(A)$ is a positive integer, and $|A|$ is a prime power, then $A$ is hyperbolic.
        \item The Gauss sum is multiplicative with respect to orthogonal direct sums, i.e. $\Sigma\left(\bigoplus_iA_i\right)=\prod_i\Sigma(A_i)\,.$
    \end{enumerate}
\end{lemma}

The following pre-metric groups will appear throughout this article, and so we give them some notation

\begin{definition}\label{def:StandardHyperbolic}
    The \emph{standard hyperbolic} pairing on $K_4=\langle a,b\rangle$ is the nondegenerate bicharacter $h(a^ib^j,a^kb^\ell)=(-1)^{i\ell + jk}$.
    There are two isometry classes of $h$-admissible quadratic forms over $\mathbb R$, and they are distinguished by the rules:
    \begin{itemize}
        \item $q_+(x)=-1$ for exactly 1 element $x\in K_4$, or
        \item $q_-(x)=-1$ for all $x\in K_4\setminus\{1\}$.
    \end{itemize}
    We will call the corresponding metric groups $K_{4,\pm}=(K_4,q_\pm)$ respectively.  Note that $K_{4,+}$ is hyperbolic, whereas $K_{4,-}$ is anisotropic.
\end{definition}

\begin{remark}
    The terms hyperbolic, (an)isotropic, and Lagrangian all have analogues for bilinear forms, but the connection between the bilinear form terminology and the quadratic form terminology can be subtle.
    For example, an element $a\in A$ is called isotropic with respect to $\chi$ if $\chi(a,-)$ is trivial, and this does not imply that $\sigma(a)=1$ in the case that $\chi=\delta\sigma$.
    The use of the word \emph{hyperbolic} in Definition \ref{def:StandardHyperbolic} refers to the fact that $h$ has a Lagrangian subgroup \emph{as a bilinear form} (bicharacter).
    Note in particular that non-hyperbolic quadratic forms can give rise to hyperbolic bicharacters.  
\end{remark}

Observe that for any pre-metric group $(A,\sigma)$, its `norm-square' $(A,\sigma)\oplus(A,\sigma^{-1})$ is hyperbolic via the diagonal embedding, so in particular $(K_{4,-})^2$ is hyperbolic.
In fact, more can be said.
The isomorphism that sends the ordered basis $(a_1,b_1,a_2,b_2)$ for $(K_{4,+})^2$ to $(a_1a_2,a_1b_2,b_1a_2b_2,a_1b_1a_2b_2)$ in $(K_{4,-})^2$ is an isometry $(K_{4,+})^2\cong(K_{4,-})^2$.  This observation leads to the following result.

\begin{proposition} \label{prop:OrbitEquivalenceCharacterization}
    Suppose $\mathbb K=\mathbb R$, and that there is some basis for $K_4^n$ with respect to which $\delta\sigma=h^n$.
    The metric group $(K_{4}^n,\sigma)$ is hyperbolic if and only if $\Sigma(\sigma)=2^n$, and in this case, $(K_{4}^n,\sigma)\cong(K_{4,+})^n$.
    If not, then $\Sigma(\sigma)=-2^n$ and $(K_{4}^n,\sigma)\cong K_{4,-}\oplus (K_{4,+})^{n-1}$.
\end{proposition}

\begin{proof}
    By hypothesis, we can choose some basis for which $\delta\sigma=h^n$, and in this way, establish an isometry $(K_4^n,\sigma)\cong(K_{4,-})^k\oplus(K_{4,+})^{n-k}$.
    By our previous observation, $(K_{4,-})^2\cong(K_{4,+})^2$, and so copies of $(K_{4,-})$ can be canceled out in pairs until there is at most one copy left.
    The Gauss sum condition then follows from Lemma \ref{lem:GaussSumProperties} parts (ii) and (iii) and (iv).
\end{proof}

Because the sign of the Gauss sum of the pre-metric group $(K_4^n,\sigma)$ determines its isometry class (assuming $\delta\sigma=h^n$), it will be convenient to establish some notation.

\begin{notation}\label{not:QF}
    For any $\sigma\in\QF(h^n)$, the sign $\sgn(\sigma)$ of the quadratic form $\sigma\colon K_4^n\to\mathbb R^\times$ is
    \[\sgn(\sigma):=\frac{\Sigma(\sigma)}{|\Sigma(\sigma)|}\, .\]
    We write $\QF_+^n$ and $\QF_-^n$ for the sets of $h^{n}$-admissibles with positive and negative sign, respectively. 
\end{notation}

\begin{proposition} \label{prop:StabilizerCombinatorics}
Let $n \geq 0$, and $\GO^n_\pm$ be the stabilizers in $\Aut(K_4^n, h^{n})$ of elements in $\QF^n_\pm$. Then 
\begin{align*}
	|\GO^n_+| &= 2^{n^2 -n + 1}(2^n - 1)\prod_{i=1}^{n - 1}(2^{2i} - 1) \\
	|\GO^n_-| &= 2^{n^2 -n + 1}(2^n + 1)\prod_{i=1}^{n - 1}(2^{2i} - 1) .
\end{align*}
Consequently, 
\begin{align*}
|\QF_+^n| &= 2^{n - 1}(2^n + 1) \\
|\QF^n_-| &= 2^{n - 1}(2^n - 1) = 2^{2n} - |\QF^n_+|. 
\end{align*}
\end{proposition}
\begin{proof}
This is a straightforward application of the orbit stabilizer theorem, once the orders of the groups preserving $h^n$ and elements in $\QF_\pm$, respectively are known. The first can be obtained from Theorem 6.18 of \cite{jacobson2009basic} as
\[
|\Aut(A, \chi)| = 2^{n^2}\prod_{i = 1}^{n} (2^{2i} - 1).
\]
 Locating a reference for $|\GO^n_\pm|$ is slightly more difficult, as many texts omit even characteristic from their treatments. Regardless, the results are quite classic - the stated orders are derived on page 206 of \cite{MR104735} (whose first edition was published in 1901!) 
\end{proof}
\begin{remark}
    The groups $\GO_\pm^n$ are the \textit{general orthogonal groups} of vector spaces over the finite field $\mathbb{F}_2$. That there are two families of such groups is related to the fact that the correspondence between quadratic forms and symmetric bilinear forms is no longer bijective in even characteristic. 
\end{remark}
Let $\ell$ be the unique nondegenerate bicharacter on $\mathbb{Z}/2\mathbb{Z}$. Observe that $\QF_{\mathbb R}(\ell^2)=\emptyset$, whereas  $|\QF_{\mathbb C}(\ell^2)|=4$.

Two of these forms over $\mathbb C$ are isometric to one another, so
we find that there are exactly three isometry classes of quadratic forms on $K_4$ inducing $\ell^{2}$.

\begin{proposition}\label{prop:StabilizerCombinatorics2ElectricBoogaloo}
Let $n > 0$. Then there are exactly four equivalence classes of complex-valued quadratic forms on $K_4^n \times K_4$ inducing $h^{n} \oplus \ell^{2}$.  When $n = 0$, there are three. 
\end{proposition} 
\begin{proof}
By the remark preceding the proof, we may assume $n > 0$. A quadratic form on $K_4^n \times K_4$ with coboundary $h^{n} \oplus \ell^{2}$, determines and is uniquely determined by a pair of quadratic forms on $K_4^n$ and $K_4$ with coboundaries $h^{n}$ and $\ell^2$ respectively. So there are at most six equivalence classes of quadratic forms with coboundary $h^{n} \oplus \ell^{2}$. We claim there are exactly four. Let us fix some notation.

We label the elements of the first factor $K_4^n$ by $a_k$ and $b_k$ respectively, and we let $g_1, g_2$ be the two elements of the second factor with self-pairing $-1$. Given a triple of signs $(\kappa, \epsilon_1, \epsilon_2)$ we denote by $\sigma(\kappa,\epsilon_1, \epsilon_2)$ the quadratic form with 
$$\sgn(\sigma(\kappa,\epsilon_1, \epsilon_2)|_{K_4^n}) = \kappa, \quad [\sigma(\kappa,\epsilon_1, \epsilon_2)](g_k) = i\epsilon_k.$$

Using the multiplicativity of the Gauss sum from Lemma \ref{lem:GaussSumProperties}, the Gauss sums of these forms are given by the formula
\[\Sigma\big(\sigma(\kappa,\epsilon_1,\epsilon_2)\big)\;=\;(\kappa\cdot2^n)\cdot(1+i\epsilon_1)\cdot(1+i\epsilon_2)\,.\]
We collect the various values $\Sigma\big(\sigma(\kappa,\epsilon_1,\epsilon_2)\big)$ into a table:
\begin{center}
    \begin{tabular}{|c|c|c|c|c|c|c|} \hline 
        $(\kappa, \epsilon_1, \epsilon_2)$ &   $(+--)$ & $(+++)$ & $(+-+)$ & $(---)$ 
        &$(-++)$ & $(--+)$ \\   \hline
        $\Sigma\big(\sigma(\kappa, \epsilon_1, \epsilon_2)\big)$ & $-2^{n + 1}i$ & $2^{n + 1}i$ & $2^{n + 1}$ & $2^{n + 1}i$ & $-2^{n + 1}i$ & $-2^{n + 1}$ \\\hline
    \end{tabular}
\end{center}

Now let $f$ be the automorphism with $$f(a_1) = a_1g_1g_2, f(b_1) = b_1g_1g_2, f(g_1) = a_1b_1g_1, f(g_2) = a_1b_1g_2$$ and which fixes $a_j, b_j$ for $j > 1$.  
Direct computations show that $f$ interchanges the forms $(---)$ and $(+++)$, as well as $(+--)$ and $(-++)$, fixes the remaining two equivalence classes, and preserves $h^{n} \oplus \ell ^{2}$.  The calculations of the Gauss sums in the above table show the resulting equivalence classes are indeed distinct.
\end{proof}

We conclude with a recognition theorem for the powers of the standard hyperbolic pairing $h^n$. 
\begin{theorem}[{\cite[Theorem 6.3]{jacobson2009basic}}] \label{thm:WallClassification} Let $\chi$ be a symmetric nondegenerate bilinear form on $(\mathbb Z /2\mathbb Z)^n$. Suppose moreover that $\chi(a, a) = 1$ for all $a \in (\mathbb Z /2\mathbb Z)^n$.  Then $((\mathbb Z /2\mathbb Z)^n, \chi)$ is isomorphic to a power of the standard hyperbolic pairing. In particular, $n$ must be even.
\end{theorem}

\section{Braidings on Split Real Tambara-Yamagami Categories} \label{sec:SplitReal}
In this section we examine the split real case with the primary purpose of setting a foundation for the non-split cases and illustrating the method. We obtain some new results, but much of the analysis in this section is originally due to Siehler \cite{sie00}, with a more contemporary perspective on the results due to Galindo \cite{GALINDO_2022}.
We begin by recalling the classification of monoidal
structures on split Tambara-Yamagami categories in
\cite{ty98}:

\begin{theorem}[{\cite[Theorem 3.2]{ty98}}]
Let $A$ be a finite group, let
$\tau=\frac{\pm 1}{\sqrt{|A|}}$,
and let $\chi\colon A\times A\to k^{\times }$ be a symmetric nondegenerate
bicharacter.
We define a split fusion category $\mathcal{C}_{\mathbb{R}}(A,\chi,\tau)$ by taking
the underlying fusion ring to be $\mathsf{TY}(A)$, the unitor isomorphisms
to be identity, and the associators to be
\begin{align*}
    \alpha_{a,b,c} &= 1_{abc}, \\
    \alpha_{a,b,m} = \alpha_{m,a,b} &= 1_{m}, \\
    \alpha_{a,m,b} &= \chi(a,b)\cdot 1_{m}, \\
    \alpha_{a,m,m} = \alpha_{m,m,a} &= \bigoplus_{b\in A} 1_{b}, \\
    \alpha_{m,a,m} &= \bigoplus_{b\in A} \chi(a,b)\cdot 1_b, \\
    \alpha_{m,m,m} &= (\tau\chi(a,b)^{-1}\cdot 1_m)_{a,b}.
\end{align*}
All split fusion categories over $k$ with fusion ring $\mathsf{TY}(A)$
arise this way, and two fusion categories
$\mathcal{C}_{\mathbb{R}}(A,\chi,\tau)$ and $\mathcal{C}_{\mathbb{R}}(A',\chi',\tau')$ are equivalent
if and only if $\tau=\tau'$ and there is a group isomorphism
$\phi\colon A\to A'$ such that $\chi(\phi(a),\phi(b))=\chi'(a,b)$
for all $a,b\in A$.
\end{theorem}

In the split case,
\mbox{$\End(X)\cong \mathbb{R}$}
for all simple objects $X\in \mathcal{C}$,
and each Hom space is spanned by a single non-zero vector.
The associators are computed in \cite{ty98} using a set of fixed normal bases, denoted in string diagrams by trivalent vertices:

\newcommand{\TSize}{0.45}
\newcommand{\abNode}{
\begin{tikzineqn}[scale=\TSize]
\coordinate (top)            at (0,1);
\coordinate (bottom left)   at (-1,-1);
\coordinate (bottom right)    at (1,-1);
\draw[strand a]  (0,0)  to (bottom left)
    node[below left] {$a$};
\draw[strand b]  (0,0)  to (bottom right)
    node[below right, yshift=0.1cm] {$b$}; 
\draw[strand ab] (0,0)  to (top)
    node[above] {$ab$};
\end{tikzineqn}}

\[
\begin{matrix}
[a,b] & = & \abNode
\quad&\quad
[a,m] & = &
\begin{tikzineqn}[scale=\TSize]
\TrivalentVertex{a}{m}{m} \end{tikzineqn}
\\
[m,a] & = &
\begin{tikzineqn}[scale=\TSize] \TrivalentVertex{m}{a}{m} \end{tikzineqn}
\quad&\quad
[a] & = &
\begin{tikzineqn}[scale=\TSize] \TrivalentVertex{m}{m}{a} \end{tikzineqn}
\end{matrix}
\]

Using the basis vectors, our set of non-trivial linear isomorphisms
$(\beta_{x,y}^{*})_{z}\in \mathrm{GL}_1(\mathbb{R})$
can be written as a set of coefficients in $\mathbb{R}^{\times }$
\begin{align*}
(\beta_{a,b}^{*})_{ab}([b,a]) &:= \sigma_{0}(a,b) [a,b] \\
(\beta_{a,m}^{*})_{m}([m,a]) &:= \sigma_{1}(a) [a,m] \\
(\beta_{m,a}^{*})_{m}([a,m]) &:= \sigma_{2}(a) [m,a] \\
(\beta_{m,m}^{*})_{a}([a]) &:= \sigma_{3}(a) [a]
\end{align*}
thus defining coefficient functions $\sigma_i$ that take inputs in $A$ and produce outputs in $\mathbb{R}^{\times}$.

\begin{remark}
Since $\chi\colon A\times A\to \mathbb{R}^{\times}$ is a bicharacter
and $A$ is a finite group, the image of $\chi$ is a finite subgroup of
$\mathbb{R}^{\times}$, so it is a subset of $\{\pm 1\}$.
This implies that for all $a\in A$, we have
\[
  \chi(a^2,-) = \chi(a,-)^2 = 1,
\]
and by nondegeneracy we have $a^2=1_{A}$.
Thus, $A$ is an elementary abelian 2-group with
$A\cong (\mathbb{Z}/2\mathbb{Z})^{m}$ for some $m\in \mathbb{Z}_{\ge 0}$.
In particular, we have $a^{-1}=a$ for all $a\in A$,
so we may freely drop inverse signs on group elements
and on $\chi$.
\end{remark}

\subsection{The hexagon equations}
After fixing bases for the Hom spaces,
we obtain a set of real valued equations
by performing precomposition on our chosen basis vectors
using graphical calculus.
The resulting unsimplified hexagon equations are as follows:
(hexagon equations)
\begin{align}
   \sigma_0(c,ab)
       &= \sigma_0(c,a)\sigma_0(c,b), \label{eqn:hexR1} \\
   \sigma_2(ab)
       &= \sigma_2(a)\chi(a,b)\sigma_2(b), \label{eqn:hexR2} \\
   \sigma_0(b,a)\sigma_1(b)
       &= \sigma_1(b)\chi(a,b), \label{eqn:hexR3} \\
   \sigma_1(b)\sigma_0(b,a)
       &= \chi(b,a)\sigma_1(b), \label{eqn:hexR4} \\
   \chi(a,b)\sigma_3(b)
       &= \sigma_2(a)\sigma_3(a^{-1}b), \label{eqn:hexR5} \\
   \sigma_3(b)\chi(a,b)
       &= \sigma_3(ba^{-1})\sigma_2(a), \label{eqn:hexR6} \\
   \sigma_0(a,ba^{-1})
       &= \sigma_1(a)\chi(a,b)\sigma_1(a), \label{eqn:hexR7} \\
   \sigma_3(a)\tau\chi(a,b)^{-1}\sigma_3(b)
       &= \sum_{c\in A}\tau\chi(a,c)^{-1}\sigma_2(c)\tau\chi(c,b)^{-1},
          \label{eqn:hexR8}
\end{align}
(inverse hexagon equations)
\begin{align}
    \sigma_0(c,a)\sigma_0(b,a)
        &= \sigma_0(bc,a), \label{eqn:hexR9} \\
    \chi(b,a)^{-1}\sigma_2(a)
        &= \sigma_2(a)\sigma_0(b,a), \label{eqn:hexR10} \\
    \sigma_0(b,a)\sigma_2(a)
        &= \sigma_2(a)\chi(a,b)^{-1}, \label{eqn:hexR11} \\
    \sigma_1(b)\chi(a,b)^{-1}\sigma_1(a)
        &= \sigma_1(ab), \label{eqn:hexR12} \\
    \sigma_0(a^{-1}b,a)
        &= \sigma_2(a)\chi(a,b)^{-1}\sigma_2(a), \label{eqn:hexR13} \\
    \sigma_3(a^{-1}b)\sigma_1(a)
        &= \sigma_3(b)\chi(a,b)^{-1}, \label{eqn:hexR14} \\
    \sigma_1(a)\sigma_3(ba^{-1})
        &= \chi(a,b)^{-1}\sigma_3(b), \label{eqn:hexR15} \\
    \sigma_3(a)\tau \chi(a,b)\sigma_3(b)
        &= \sum_{c\in A} \tau \chi(a,c)\sigma_1(c)\tau \chi(c,b).
           \label{eqn:hexR16}
\end{align}

\subsection{Reduced hexagon equations}
The following six equations are algebraically equivalent
to the sixteen unsimplified hexagon equations:
\begin{align}
    &\sigma_0(a,b) = \chi(a,b), \label{eqn:reducedR1} \\
    &\sigma_1(a)^2 = \chi(a,a), \label{eqn:reducedR2} \\
    &\sigma_1(ab)  = \sigma_1(a)\sigma_1(b)\chi(a,b), \label{eqn:reducedR3} \\
    &\sigma_2(a)   = \sigma_1(a), \label{eqn:reducedR4} \\
    &\sigma_3(1)^2 = \tau \sum_{c\in A}\sigma_1(c), \label{eqn:reducedR5} \\
    &\sigma_3(a)   = \sigma_3(1)\sigma_1(a)\chi(a,a). \label{eqn:reducedR6}
\end{align}
The process of eliminating redundancies is as follows.
First, we may eliminate any term that appears on both sides
of any equation, as all functions are valued in the $\{\pm1\}$.
Then, we have the following implications:
\begin{center}
\begin{tabular}{|c|c|c|}
\hline
(\ref{eqn:hexR3})$\implies$ (\ref{eqn:reducedR1})
&
(\ref{eqn:hexR12})$\implies$ (\ref{eqn:reducedR3})
&
(\ref{eqn:hexR16}),
$a=b=1$
$\implies$ (\ref{eqn:reducedR5})
\\ \hline
(\ref{eqn:hexR7})$\implies$ (\ref{eqn:reducedR2})
&
(\ref{eqn:hexR6}),
(\ref{eqn:hexR15})
$\implies$ (\ref{eqn:reducedR4})
&
(\ref{eqn:hexR14}),
$a=b$
$\implies$ (\ref{eqn:reducedR6})
\\ \hline
\end{tabular}
\end{center}
To check that the reduced equations are indeed equivalent
to the original sixteen, first note that
the equality $\sigma_2=\sigma_1$ from equation
(\ref{eqn:reducedR4}) identifies each of
(\ref{eqn:hexR9})-(\ref{eqn:hexR16}) with one of
(\ref{eqn:hexR1})-(\ref{eqn:hexR8}), so it suffices to prove
the first eight hexagons from the reduced equations.
Equations
(\ref{eqn:hexR1}),
(\ref{eqn:hexR3}) and
(\ref{eqn:hexR4})
follows from equation
(\ref{eqn:reducedR1})
which identifies $\sigma_0=\chi$ to be a bicharacter.
Equation (\ref{eqn:hexR2}) follows from
(\ref{eqn:reducedR3}) and (\ref{eqn:reducedR4}).
Equation 
(\ref{eqn:hexR7}) follows from
(\ref{eqn:reducedR2}).
Equations
(\ref{eqn:hexR5}) and
(\ref{eqn:hexR6}) can be derived by expanding both sides in terms of $\sigma_1$ and $\chi$ using equations \eqref{eqn:reducedR4} and \eqref{eqn:reducedR6}. 

It remains to derive equation 
(\ref{eqn:hexR8}). First, equation \eqref{eqn:reducedR3} implies 
\begin{equation} \label{eqn:Sigma1Expansion}
    \sigma_1(a)\sigma_1(b)\sigma_1(d) = \frac{\sigma_1(abd)}{\chi(a, bd)\chi(b,d)}
\end{equation}

Finally we derive an equivalent form of \eqref{eqn:hexR8} from the reduced equations, along with the fact that $\chi$ is a $\{\pm 1\}$-valued symmetric bicharacter. 
\begin{align*}
    \sigma_3(a)\chi(a,b)^{-1}\sigma_3(b) &\overset{\eqref{eqn:reducedR6}}{=} \sigma_3(1)^2\sigma_1(a)\sigma_1(b)\chi(a,a)\chi(b,b)\chi(a,b)^{-1} \\ 
    &\overset{\eqref{eqn:reducedR5}}{=} 
    \tau \sum_{d\in A}\sigma_1(d)\sigma_1(a)\sigma_1(b)\chi(a,a)\chi(b,b)\chi(a,b)^{-1}\\
    &\overset{\eqref{eqn:Sigma1Expansion}}{=} 
    \tau \sum_{d\in A}\sigma_1(abd)\frac{\chi(a,a)\chi(b,b)}{\chi(a,b)\chi(a, bd)\chi(b,d)} \\
    &\overset{c := abd}{=} 
    \tau \sum_{c\in A}\sigma_1(c)\frac{\chi(a,a)\chi(b,b)}{\chi(a,b)\chi(a, a^{-1}c)\chi(b,b^{-1}a^{-1}c)}\\
    &\overset{\eqref{eqn:reducedR4}}{=}  \tau\sum_{c\in A}\chi(a,c)^{-1}\sigma_2(c)\chi(c,b)^{-1}
\end{align*}

\subsection{Classification of Braidings}

By equation (\ref{eqn:reducedR2}) and the fact
that all coefficients are real, we have the restriction that
$\chi(a,a)>0$ for all $a\in A$. We conclude using Theorem \ref{thm:WallClassification}:
\begin{proposition}\label{thm:SplitClassification}
If $\mathcal{C}_{\mathbb{R}}(A,\chi,\tau)$ admits a braiding,
then $A\cong K_4^{n}$ for some $n\in \mathbb{Z}_{\ge 0}$
and $\chi$ is the hyperbolic pairing $h^{n}$.
\end{proposition}

From the simplified hexagon equations, we have the
following classification of braidings
on a split TY category over $\mathbb{R}$.

\begin{theorem}\label{thm:split-class-sols}
A braiding on
$\mathcal{C}_{\mathbb{R}}(K_4^n,h^{n},\tau)$
is given by a $\chi$-admissible function $\sigma$
with $\sgn\sigma=\sgn\tau$ and a coefficient
$\epsilon\in \{\pm 1\}$.
In other words, the set of braidings on
$\mathcal{C}_{\mathbb{R}}(K_4^n,h^{n},\tau)$
is in bijection 
with $\QF_{\sgn\tau}^n \times \{\pm 1\}$.
\end{theorem}

\begin{proof}
Given a braiding on
$\mathcal{C}_{\mathbb{R}}(K_4^n,h^{n},\tau)$, we deduce from the reduced hexagon equations (namely \ref{eqn:reducedR3}) that $\sigma_1 \in \QF(h^{n})$
Equation (\ref{eqn:reducedR5}) gives the constraint
\[
    \tau \sum_{c\in A}\sigma_1(c)
    = 2^{n}\tau\sgn{\sigma_1}>0,
\]which tells us that $\sigma_1 \in \QF^n_{\sgn(\tau)}$. We may also extract a sign $\epsilon$ which is defined by the equation 
\begin{equation} \label{eqn:RealSigma31Definition}
    \sigma_3(1) = \epsilon \sqrt{2^{n}\tau\sgn{\sigma_1}} .
\end{equation} We thus obtain an element $(\sigma_1, \epsilon)  \in \QF^n_{\text{sgn}(\tau)} \times \{\pm 1\}$.

Conversely, given an element $(\sigma, \epsilon) \in \QF^n_{\text{sgn}(\tau)} \times \{\pm 1\}$, we let $\sigma_1 = \sigma_2 = \sigma$, $\sigma_0 = h^{n}$ and $\sigma_3(1)$ by Equation \eqref{eqn:RealSigma31Definition}. We can then extend $\sigma_3(1)$ to a function $\sigma_3(a)$ by equation \eqref{eqn:reducedR6}. Equations \eqref{eqn:reducedR1}-\eqref{eqn:reducedR4} and \eqref{eqn:reducedR6} hold by our definitions along with that fact that $\sigma \in \QF(h^{n})$. The remaining constraint \eqref{eqn:reducedR5} holds by Proposition \ref{prop:OrbitEquivalenceCharacterization}, our choice of $\sigma_3(1)$ and the definition of $\QF^n_{\text{sgn}(\tau)}$. Finally, we observe that these procedures are, by construction, mutually inverse. 
\end{proof}

Note that when $n=0$, $\sgn(\sigma)$ is automatically equal to 1. If we additionally impose $\tau < 0$, then the proof above shows that $\sigma_3(1)$ must be purely imaginary, and thus such categories can only exist over fields containing a square root of $-1$.
Over $\mathbb C$, $\sigma_3(1)=i$ gives the semion category, and $\sigma_3(1)=-i$ gives the reverse semion.

Over $\mathbb R$, \eqref{eqn:RealSigma31Definition} cannot be satisfied when $n=0$ and $\tau<0$, and so this category admits no braidings (i.e. $\QF^0_{-}=\emptyset$).

As a consequence of Theorem \ref{thm:split-class-sols}, the following braidings are coherent.
\begin{definition}\label{defn:ExplicitSplitRealBraidings}
Given an element $(\sigma, \epsilon)$ of $\QF_{\sgn\tau}^n\times \{\pm 1\}$, we define a braided structure $\mathcal{C}_\mathbb{R}(K_4^n,h^{n},\tau,\sigma,\epsilon)$ on $\mathcal{C}_\mathbb{R}(K_4^n,h^{n},\tau)$ by:
\begin{align*}
	\beta_{a,b} &= \chi(a,b)\cdot \id_{ab}, \\
	\beta_{a,m} &= \beta_{m,a} = \sigma(a)\cdot \id_{m}, \\
	\beta_{m,m} &= \sum_{a\in K_4^{n}}
	\epsilon\,\sigma(a) [a]^{\dag}[a].
\end{align*}
Since the group $K_4^n$, bicharacter $h^{n}$, and coefficient $\tau$ are determined from context, we will abbreviate $\mathcal{C}_\mathbb{R}(K_4^n,h^{n},\tau,\sigma,\epsilon) := \mathcal{C}_\mathbb{R}(\sigma,\epsilon)$. 
\end{definition}

We next analyze when $\mathcal{C}_\mathbb{R}(\sigma,\epsilon)$ is braided equivalent to $\mathcal{C}_\mathbb{R}(\sigma', \epsilon')$, by analyzing the properties of certain categorical groups attached to these categories.

\begin{notation}\label{not:CatGrp}
The autoequivalences of any ($\star=$ plain, monoidal, braided, etc.) category $\mathcal C$ form a categorical group $\Aut_{\star}(\mathcal C)$.
The objects of $\Aut_{\star}(\mathcal C)$ are $\star$-autoequivalences of $\mathcal C$, and the morphisms are $\star$-natural isomorphisms.
For any categorical group $\mathcal G$, the group of isomorphism classes of objects is denoted by $\pi_0\mathcal G$, and the automorphisms of the identity are denoted by $\pi_1\mathcal G$.
\end{notation}

\begin{lemma}\label{lem:SplitRealFunctorClassification}
	$$\pi_0\Aut_\otimes\big(\mathcal{C}_\mathbb{R}(K_4^n,h^{n},\tau)\big) \cong \Aut(K_4^n,h^{n})$$
\end{lemma}
\begin{proof}
This fact appears in several places in the literature (for instance \cite[Proposition 1]{Tambara2000}, \cite[Proposition 2.10]{Nikshych2007NongrouptheoreticalSH}, and \cite[Lemma 2.16]{EDIEMICHELL2022108364}) and is proved with arguments that do not depend on the algebraic closure of the field in question. They do, however, assume that the underlying semisimple category is split. We will see in future sections that this does affect the validity of the conclusion. 
\end{proof}
\begin{proposition}\label{prop:RealFunctorBraided}
The monoidal functor $F(f)$ determined by an automorphism $f\in\Aut(K_4^n,h^{n})$ forms a braided monoidal equivalence $\mathcal{C}_\mathbb{R}(\sigma,\epsilon) \to \mathcal{C}_\mathbb{R}(\sigma',\epsilon')$ if and only if $f \cdot \sigma = \sigma'$ and $\epsilon = \epsilon'$.
\end{proposition}
\begin{proof}
	 Using Definition \ref{defn:ExplicitSplitRealBraidings}, the required constraints for $F(f)$ to be braided are
	\begin{align*}
		h^{n}(f(a), f(b)) &= h^{n}(a, b)  \\
		\sigma'(f(a)) &= \sigma(a) \\ 
		\epsilon' &= \epsilon.
	\end{align*}
    These equations are indeed equivalent to $f \cdot \sigma = \sigma'$ and $\epsilon = \epsilon'$.
\end{proof}
The following theorem strengthens \cite{GALINDO_2022} in the split real case. 
\begin{theorem}\label{thm:SplitCaseEquivalence}

There is a braided equivalence $\mathcal{C}_\mathbb{R}(\sigma,\epsilon) \sim \mathcal{C}_\mathbb{R}(\sigma',\epsilon')$ if and only if $\epsilon = \epsilon'$. In particular, there are exactly two equivalence classes of braidings on $\mathcal{C}_\mathbb{R}(K_4^n,h^{n},\tau)$ when $n > 0$, or when $n = 0$ and $\tau > 0$, and zero otherwise.
\end{theorem}
\begin{proof}
By Lemma \ref{lem:SplitRealFunctorClassification}, the functors $F(f)$ form a complete set of representatives for $\pi_0(\Aut(\mathcal{C}_\mathbb{R}(K_4^n,h^{n},\tau)))$. Therefore it suffices to check when some $F(f)$ is a braided equivalence  $\mathcal{C}_\mathbb{R}(\sigma,\epsilon) \to \mathcal{C}_\mathbb{R}(\sigma',\epsilon')$. By Proposition \ref{prop:RealFunctorBraided}, this occurs exactly when $\epsilon = \epsilon'$ and $\sigma$ is orbit equivalent to $\sigma'$. This last condition always holds by Proposition \ref{prop:OrbitEquivalenceCharacterization} since the sign of $\sigma$ is determined by $\tau$ (part of the underlying monoidal structure). 
\end{proof}
Taking $\epsilon = \epsilon'$ and $\sigma = \sigma'$ in Proposition \ref{prop:RealFunctorBraided}, and using the notation of Proposition \ref{prop:StabilizerCombinatorics}, we obtain:
\begin{proposition}\label{prop:SplitRealBraidedFunctorClassification}
	$$\pi_0(\Aut_{\text{br}}(\mathcal{C}_\mathbb{R}(\sigma, \epsilon))) \cong \GO^n_{\sgn \sigma}.$$
\end{proposition}
 Note that by Propositions \ref{prop:SplitRealBraidedFunctorClassification} and  \ref{prop:StabilizerCombinatorics}, $|\pi_0\Aut_{\text{br}}(\mathcal{C}_\mathbb{R}(\sigma, \epsilon)|$ depends on $\tau$, while Lemma \ref{lem:SplitRealFunctorClassification} shows that $|\pi_0\Aut_\otimes(\mathcal{C}_\mathbb{R}(K_4^n,h^{n},\tau))|$ does not. 
 \begin{remark}
	When $n = 1$ (but $\tau$ is not fixed), braidings on the split complex Tambara-Yamagami categories were classified in \cite[Example 2.5.2, Figures 3-5]{SchopierayNonDegenExtension}. We can see that the four symmetrically braided categories appearing ibid., Figure 3, are defined over the reals; our results here show that these are in fact the only possibilities.
\end{remark}
We conclude with a lemma on twist morphisms for these braidings. 
\begin{lemma}
There are exactly two families of twist morphisms for any $\mathcal{C}_\mathbb{R}(\sigma,\epsilon)$, corresponding to a sign $\rho \in \{\pm 1\}$. These twists are indeed ribbon structures (in the sense of \cite[Definition 8.10.1]{EGNO15}).
\end{lemma}
\begin{proof}
The first part of the remark is due to \cite{sie00}, who gives the components $\theta_x$ of the twist as $\theta_a = 1, \theta_m = \rho \sigma_3(1)^{-1}$. Since every simple object is self dual, the required axiom is simply $\theta_m = \theta_m^*$. But this holds as a result of the linearity of composition.
\end{proof}

\section{Braidings on Real/Quaternionic Tambara-Yamagami Categories}\label{sec:RealQuaternionic}

We will now examine the case where $\End(\mathbbm{1})\cong \mathbb{R}$
and $\End(m)\cong \mathbb{H}$.
We first note that the four dimensional $\mathbb{R}$ vector spaces
$\Hom(a\otimes m,m)$, $\Hom(m\otimes a,m)$ and $\Hom(m\otimes m,a)$
can be endowed with the structure of $(\mathbb{H},\mathbb{H})$-bimodules
under pre- and postcomposition with quaternions.
By naturality, the effect of precomposing with braiding isomorphisms
for each of these hom-spaces is determined on an ($\mathbb{H},\mathbb{H}$)-basis.
A preferred system of basis vectors
(over $\mathbb{R}$ for $\Hom(a\otimes b,ab)$ and over $\mathbb{H}$ 
for the others) is chosen in \cite[Section 5.1]{pss23},
depicted again as trivalent vertices:
\[
    \begin{matrix}
        [a,b] & = & \abNode
        \quad&\quad
        [a,m] & = &
        \begin{tikzineqn}[scale=\TSize]
        \TrivalentVertex{a}{m}{m} \end{tikzineqn}
        \\
        [m,a] & = &
        \begin{tikzineqn}[scale=\TSize] \TrivalentVertex{m}{a}{m} \end{tikzineqn}
        \quad&\quad
        [a] & = &
        \begin{tikzineqn}[scale=\TSize] \TrivalentVertex{m}{m}{a} \end{tikzineqn}
    \end{matrix}
\]
Splittings to each $[a]$ is chosen in
\cite[Proposition 4.4]{pss23} and will be denoted
by
\[        
[a]^\dagger =
\begin{tikzineqn}[scale=\TSize,yscale=-1]
    \coordinate (mid)            at (0,0);
    \coordinate (top)            at (0,1);
    \coordinate (bottom left)   at (-1,-1);
    \coordinate (bottom right)    at (1,-1);
    \draw[strand m] (mid) to (bottom left) node[above left] {$m$};
    \draw[strand m] (mid) to (bottom right) node[above right] {$m$};
    \draw[strand a] (mid) to (top) node[below] {$a$};
\end{tikzineqn}
\]
such that
\[
[a](\id_m\otimes h)[a]^\dag
\quad=\quad
\begin{tikzineqn}[scale=0.5]
    \draw[strand m] (-1,0) -- (0,1) -- (1,0) -- (0,-1) -- (-1,0);
    \draw[strand a] (0,-2) -- (0,-1);
    \draw[strand a] (0,2) -- (0,1);
    \node[smallbead] at (1,0) {$h$};
\end{tikzineqn}
\quad=\quad
\mathfrak{Re}(h)\cdot
\begin{tikzineqn}[scale=0.5]
    \draw[strand a] (0,-2) -- (0,2);
\end{tikzineqn}
\quad=\quad
\mathfrak{Re}(h)\cdot\id_a\;,
\]
where $\mathfrak{Re}(w+xi+yj+zk)=w$ is the real part of the quaternion, and also such that
\[
\id_{m\otimes m}
\quad=\quad
\begin{tikzineqn}
    \draw[strand m] (0,0) -- (0,2);
    \draw[strand m] (1,0) -- (1,2);
\end{tikzineqn}
\quad=\quad
\sum_{\substack{a\in A\\ s\in S}}
\begin{tikzineqn}[scale=0.5]
    \draw[strand a] (0,0) -- (0,2);
    \draw[strand m] (0,2) -- ++(1,1);
    \draw[strand m] (0,2) -- ++(-1,1);
    \draw[strand m] (0,0) -- ++(1,-1);
    \draw[strand m] (0,0) -- ++(-1,-1);
    \node[smallbead] at (0.5,2.5) {$s$};
    \node[smallbead] at (0.5,-0.5) {$\overline{s}$};
\end{tikzineqn}
\quad=\quad
\sum_{\substack{a\in A\\ s\in S}}
(\id_m\otimes s)[a]^{\dag}[a](\id_m\otimes \overline{s})
\]
where $S:=\{1,i,j,k\}$.
By \cite[Proposition 5.1]{pss23}, the basis vectors
satisfy the convenient property that they commute
\newcommand{\beadedTSize}{0.7}
\[
    \begin{tikzineqn}[scale=\beadedTSize]
        \TrivalentVertex{a}{m}{m} 
        \DrawSmallBead{mid}{top}{v}
    \end{tikzineqn}
    \ = \
    \begin{tikzineqn}[scale=\beadedTSize]
        \TrivalentVertex{a}{m}{m} 
        \DrawSmallBead{mid}{bottom right}{v}
    \end{tikzineqn}
    \quad\quad
    \begin{tikzineqn}[scale=\beadedTSize]
        \TrivalentVertex{m}{a}{m} 
        \DrawSmallBead{mid}{top}{v}
    \end{tikzineqn}
    \ = \
    \begin{tikzineqn}[scale=\beadedTSize]
        \TrivalentVertex{m}{a}{m} 
        \DrawSmallBead{mid}{bottom left}{v}
    \end{tikzineqn}\;\,,
\]
or conjugate-commute
\[
    \begin{tikzineqn}[scale=\beadedTSize]
        \TrivalentVertex{m}{m}{a} 
        \DrawSmallBead{mid}{bottom left}{v}
    \end{tikzineqn}
    \ = \
    \begin{tikzineqn}[scale=\beadedTSize]
        \TrivalentVertex{m}{m}{a} 
        \DrawSmallBead{mid}{bottom right}{\overline{v}}
    \end{tikzineqn}
\]
with all quaternions $v\in \mathbb{H}$.
We can now recall the classification of associators
on these categories using the chosen bases.

\begin{theorem}[{\cite[Theorem 5.4]{pss23}}]
   Let $A$ be a finite group, let $\tau=\frac{\pm1}{\sqrt{4|A|}}$, and let $\chi:A\times A\to \mathbb R^\times$ be a nondegenerate symmetric bicharacter on $A$.
A triple of such data gives rise to a non-split Tambara-Yamagami category \mbox{$\mathcal{C}_{\bb H}(A,\chi,\tau)$}, with $\End(\1)\cong\bb R$ and $\End(m)\cong\bb H$, whose associators for $a, b, c\in A$ are given as follows:
	\begin{gather*}
	    \alpha_{a,b,c}=\id_{abc}\,,\\
	    \alpha_{a,b,m}=\alpha_{m,b,c}=\id_{m}\,,\\
	    \alpha_{a,m,c}=\chi(a,c)\cdot\id_{m},\\
	    \alpha_{a,m,m}=\alpha_{m,m,c}=\id_{m\otimes m}\,,\\
	    \alpha_{m,b,m}=\bigoplus_{a\in A}\chi(a,b)\cdot\id_{a^{\oplus4}}\,,\\
	    \alpha_{m,m,m}=\tau\cdot\sum_{\substack{a,b\in A\\s,t\in S}}\chi(a,b)^{-1}\cdot(s\otimes(\id_m\otimes\overline{t}))(\id_m\otimes[a]^\dagger)([b]\otimes\id_m)((\id_m\otimes s)\otimes t),
	\end{gather*}
    where $S:=\{1,i,j,k\}\subseteq \mathbb{H}$.
    The unit is $\1=\mathbf 1_A$, the identity in $A$, and the unit constraints are identity morphisms.
    Furthermore, all equivalence classes of such categories arise in this way.  Two categories $\mathcal{C}_{\bb H}(A,\chi,\tau)$ and $\mathcal{C}_{\bb H}(A',\chi',\tau')$ are equivalent if and only if $\tau=\tau'$ and there exists an isomorphism $f:A\to A'$ such that for all $a,b\in A$,
    \[\chi'\big(f(a),f(b)\big)\;=\;\chi(a,b)\,.\]
\end{theorem}

We can now write down our braiding coefficients, some of which are a priori
quaternions:

    
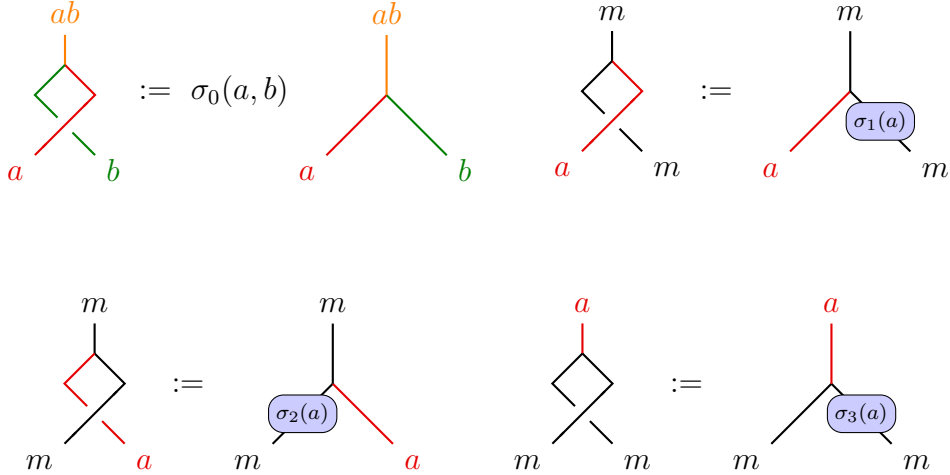
\[
        \begin{tikzineqn}[scale=\eqnscale]
            \draw[strand ab] (0,0) to ++(0,1) node[above] {$ab$};
            \begin{knot}[clip width=10]
                \strand[strand a] (0,0)
                to ++(1,-1)
                to ++(-2,-2) node[below left] {$a$};
                \strand[strand b] (0,0)
                to ++(-1,-1)
                to ++(2,-2) node[below right,yshift=0.1cm] {$b$};
            \end{knot}
        \end{tikzineqn}
        := \
        \sigma_0(a,b)
        \begin{tikzineqn}[scale=\tscale]
            \coordinate (top)            at (0,1);
            \coordinate (bottom left)   at (-1,-1);
            \coordinate (bottom right)    at (1,-1);
            \draw[strand a]  (0,0)  to (bottom left)
                node[below left] {$a$};
            \draw[strand b]  (0,0)  to (bottom right)
                node[below right, yshift=0.1cm] {$b$};
            \draw[strand ab] (0,0)  to (top)
                node[above] {$ab$};
        \end{tikzineqn}
        \quad\quad
        \begin{tikzineqn}[scale=\eqnscale]
            \draw[strand m] (0,0) to ++(0,1) node[above] {$m$};
            \begin{knot}[clip width=10]
                \strand[strand a] (0,0)
                to ++(1,-1)
                to ++(-2,-2) node[below left] {$a$};
                \strand[strand m] (0,0)
                to ++(-1,-1)
                to ++(2,-2) node[below right] {$m$};
            \end{knot}
        \end{tikzineqn}
        := \
        \begin{tikzineqn}[scale=\tscale] 
            \TrivalentVertex{a}{m}{m} 
            \DrawLongBead{mid}{bottom right}{\sigma_1(a)}
        \end{tikzineqn} 
    \]
    \vspace{-0.2cm}
    \[
        \begin{tikzineqn}[scale=\eqnscale]
            \draw[strand m] (0,0) to ++(0,1) node[above] {$m$};
            \begin{knot}[clip width=10]
                \strand[strand m] (0,0)
                to ++(1,-1)
                to ++(-2,-2) node[below left] {$m$};
                \strand[strand a] (0,0)
                to ++(-1,-1)
                to ++(2,-2) node[below right] {$a$};
            \end{knot}
        \end{tikzineqn}
        := \
        \begin{tikzineqn}[scale=\tscale] 
            \TrivalentVertex{m}{a}{m} 
            \DrawLongBead{mid}{bottom left}{\sigma_2(a)}
        \end{tikzineqn} \quad\quad \
        \begin{tikzineqn}[scale=\eqnscale]
            \draw[strand a] (0,0) to ++(0,1) node[above] {$a$};
            \begin{knot}[clip width=10]
                \strand[strand m] (0,0)
                to ++(1,-1)
                to ++(-2,-2) node[below left] {$m$};
                \strand[strand m] (0,0)
                to ++(-1,-1)
                to ++(2,-2) node[below right] {$m$};
            \end{knot}
        \end{tikzineqn}
        := \
        \begin{tikzineqn}[scale=\tscale] 
            \TrivalentVertex{m}{m}{a} 
            \DrawLongBead{mid}{bottom right}{\sigma_3(a)}
        \end{tikzineqn} 
    \]
It is clear that if the braiding coefficients
are natural if they are real-valued.
It turns out the the converse is true, in that naturality forces
all braiding coefficients to be real.

\begin{lemma} \label{lem:RQSigma12Real}
   The functions $\sigma_1$ and $\sigma_2$ are real-valued.
\end{lemma}

\begin{proof}
For any $v\in \mathbb{H}$ and any $a\in A$, consider the following diagram:
% https://q.uiver.app/#q=WzAsOCxbMCwwLCJtIl0sWzEsMSwiYVxcb3RpbWVzIG0iXSxbMiwxLCJtXFxvdGltZXMgYSJdLFsxLDIsImFcXG90aW1lcyBtIl0sWzIsMiwibVxcb3RpbWVzIGEiXSxbMywwLCJtIl0sWzAsMywibSJdLFszLDMsIm0iXSxbMSwyLCJcXGFscGhhIl0sWzIsNCwiaFxcb3RpbWVzIFxcaWRfYSJdLFsxLDMsIlxcaWRfYVxcb3RpbWVzIGgiLDJdLFszLDQsIlxcYWxwaGEiLDJdLFsxLDAsImZfMSIsMl0sWzIsNSwiZl8yIl0sWzMsNiwiZl8xIl0sWzQsNywiZl8yIiwyXSxbMCw1LCJcXHNpZ21hXzEoYSkiXSxbNSw3LCJoIl0sWzAsNiwiaCIsMl0sWzYsNywiXFxzaWdtYV8xKGEpIiwyXV0=
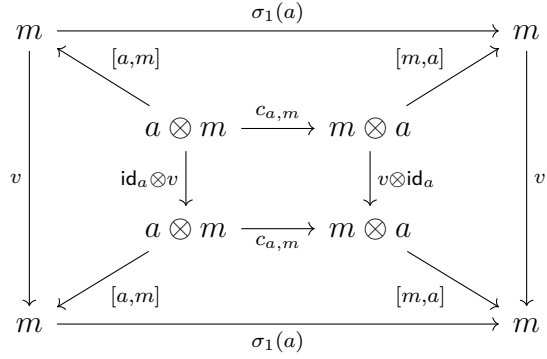
\[\begin{tikzcd}
	m &&& m \\
	& {a\otimes m} & {m\otimes a} \\
	& {a\otimes m} & {m\otimes a} \\
	m &&& m
	\arrow["c_{a,m}", from=2-2, to=2-3]
	\arrow["{v\otimes \id_a}", from=2-3, to=3-3]
	\arrow["{\id_a\otimes v}"', from=2-2, to=3-2]
	\arrow["c_{a,m}"', from=3-2, to=3-3]
	% \arrow["{f_1}"', from=2-2, to=1-1]
	% \arrow["{f_2}", from=2-3, to=1-4]
	% \arrow["{f_1}", from=3-2, to=4-1]
	% \arrow["{f_2}"', from=3-3, to=4-4]
	\arrow["{[a,m]}"', from=2-2, to=1-1]
	\arrow["{[m,a]}", from=2-3, to=1-4]
	\arrow["{[a,m]}", from=3-2, to=4-1]
	\arrow["{[m,a]}"', from=3-3, to=4-4]
	\arrow["{\sigma_1(a)}", from=1-1, to=1-4]
	\arrow["v", from=1-4, to=4-4]
	\arrow["v"', from=1-1, to=4-1]
	\arrow["{\sigma_1(a)}"', from=4-1, to=4-4]
\end{tikzcd}\]

The middle diagram commutes by the naturality of the braiding,
while the top and bottom quadrangles commute by the definition of $\sigma_1$.
As our chosen basis vector $[a,m]$ commutes with quaternions,
we have
\[
    v\circ f_1 = v \triangleright [a,m]
    = [a,m] \triangleleft v = f_1 \otimes (\id_a\otimes v)
,\]
so the left quadrangle commutes,
and the same argument can be made for the right quadrangle
using the vector $[m,a]$.
Since both $[a,m]$ and $[m,a]$ are isomorphisms,
we have the commutativity of the outer rectangle,
and thus we have that
\[
(\forall v\in \mathbb{H}) \quad \sigma_1(a)\circ v = v \circ \sigma_1(a)
\]
or that $\sigma_1(a)$ lies in the center of $\mathbb{H}$.
Alternatively, we can present the proof using graphical calculus.
We first introduce a ``bubble" by precomposing with our basis vector
and its inverse, and commute the quaternion through the trivalent vertex:
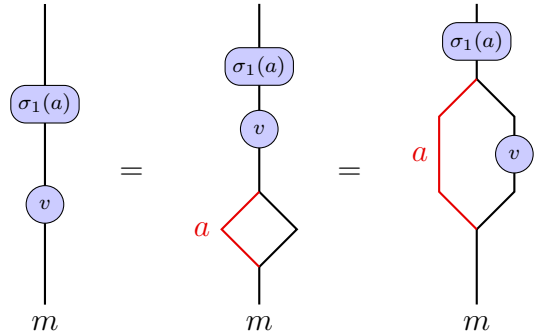
\newcommand{\lemmascale}{1}
\[
    \begin{tikzineqn}[scale=\lemmascale]
        \coordinate (bot) at (0,-2); 
        \coordinate (mid) at (0,0); 
        \coordinate (top) at (0,2); 
        \coordinate (bead1) at ($(bot)!1/3!(top)$);
        \coordinate (bead2) at ($(bot)!2/3!(top)$);
        \draw[strand m] (top) to (bot) node[below] {$m$};
        \node[bead] at (bead1) {$v$};
        \node[longbead] at (bead2) {$\sigma_1(a)$};
    \end{tikzineqn} 
    \quad=\quad
    \begin{tikzineqn}[scale=\lemmascale]
        \draw[strand m]
        node[below] {$m$}
        (0,0) to ++(0,1/2) coordinate (vert)
              to ++(1/2,1/2)
              to ++(-1/2,1/2) coordinate (triv)
              to (0,4);
        \draw[strand a]
        (vert) to ++(-1/2,1/2) node[left] {$a$}
               to ++(1/2,1/2);
        \node[bead] at ($(triv)!1/3!(0,4)$) {$v$};
        \node[longbead] at ($(triv)!2/3!(0,4)$) {$\sigma_1(a)$};
    \end{tikzineqn} 
    \quad = \quad
    \begin{tikzineqn}[scale=\lemmascale]
    \begin{knot}[clip width=10]
        \strand[strand m]
        node[below] {$m$}
        (0,0) to ++(0,1)
              to ++(1/2,1/2)
              to ++(0,1)
              to ++(-1/2,1/2)
              to ++(0,1);
        \strand[strand a]
        (0,1) to ++(-1/2,1/2)
              to ++(0,1)
              to ++(1/2,1/2);
    \end{knot}
    \node[node a,left] at (-1/2,2) {$a$};
    \node[longbead] at (0,3.5) {$\sigma_1(a)$};
    \node[bead] at (1/2,2) {$v$};
    \end{tikzineqn}
    \]
Then, by the definition of $\sigma_1$ and naturality, we have
    
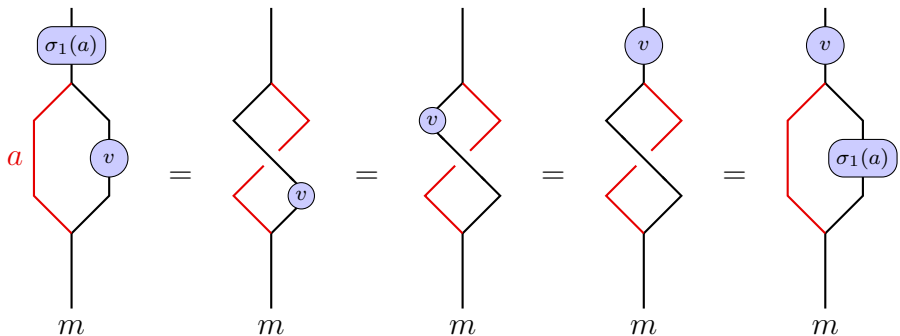
\[
    \begin{tikzineqn}[scale=\lemmascale]
    \begin{knot}[clip width=10]
        \strand[strand m]
        node[below] {$m$}
        (0,0) to ++(0,1)
              to ++(1/2,1/2)
              to ++(0,1)
              to ++(-1/2,1/2)
              to ++(0,1);
        \strand[strand a]
        (0,1) to ++(-1/2,1/2)
              to ++(0,1)
              to ++(1/2,1/2);
    \end{knot}
    \node[node a,left] at (-1/2,2) {$a$};
    \node[longbead] at (0,3.5) {$\sigma_1(a)$};
    \node[bead] at (1/2,2) {$v$};
    \end{tikzineqn}
    \quad =\quad
    \begin{tikzineqn}[scale=\lemmascale]
    \begin{knot}[clip width=10]
        \strand[strand m]
        node[below] {$m$}
        (0,0) to ++(0,1)
              to ++(1/2,1/2)
              to ++(-1,1)
              to ++(1/2,1/2)
              to ++(0,1);
        \strand[strand a]
        (0,1) to ++(-1/2,1/2)
              to ++(1,1)
              to ++(-1/2,1/2);
    \end{knot}
    \node[smallbead,xshift=-0.1cm] at (1/2,3/2) {$v$};
    \end{tikzineqn} 
    \quad = \quad
    \begin{tikzineqn}[scale=\lemmascale]
    \begin{knot}[clip width=10]
        \strand[strand m]
        node[below] {$m$}
        (0,0) to ++(0,1)
              to ++(1/2,1/2)
              to ++(-1,1)
              to ++(1/2,1/2)
              to ++(0,1);
        \strand[strand a]
        (0,1) to ++(-1/2,1/2)
              to ++(1,1)
              to ++(-1/2,1/2);
    \end{knot}
    \node[smallbead,xshift=0.1cm] at (-1/2,5/2) {$v$};
    \end{tikzineqn}
    \quad=\quad
    \begin{tikzineqn}[scale=\lemmascale]
    \begin{knot}[clip width=10]
        \strand[strand m]
        node[below] {$m$}
        (0,0) to ++(0,1)
              to ++(1/2,1/2)
              to ++(-1,1)
              to ++(1/2,1/2)
              to ++(0,1);
        \strand[strand a]
        (0,1) to ++(-1/2,1/2)
              to ++(1,1)
              to ++(-1/2,1/2);
    \end{knot}
    \node[bead] at (0,3.5) {$v$};
    \end{tikzineqn}
    \quad=\quad
    \begin{tikzineqn}[scale=\lemmascale]
        \draw[strand m]
        node[below] {$m$}
        (0,0) to ++(0,1)
              to ++(1/2,1/2)
              to ++(0,1)
              to ++(-1/2,1/2)
              to ++(0,1);
        \draw[strand a]
        (0,1) to ++(-1/2,1/2)
              to ++(0,1)
              to ++(1/2,1/2);
        \node[bead] at (0,3.5) {$v$};
        \node[longbead] at (1/2,2) {$\sigma_1(a)$};
    \end{tikzineqn} \]
and we can pass $\sigma_1(a)$ through the trivalent vertex to get
    
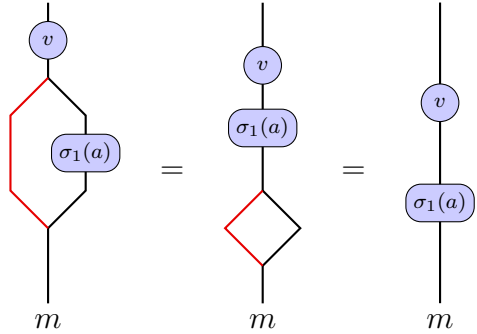
\[
    \begin{tikzineqn}[scale=\lemmascale]
        \draw[strand m]
        node[below] {$m$}
        (0,0) to ++(0,1)
              to ++(1/2,1/2)
              to ++(0,1)
              to ++(-1/2,1/2)
              to ++(0,1);
        \draw[strand a]
        (0,1) to ++(-1/2,1/2)
              to ++(0,1)
              to ++(1/2,1/2);
        \node[bead] at (0,3.5) {$v$};
        \node[longbead] at (1/2,2) {$\sigma_1(a)$};
    \end{tikzineqn}
    \quad=\quad
    \begin{tikzineqn}[scale=\lemmascale]
        \draw[strand m]
        node[below] {$m$}
        (0,0) to ++(0,1/2) coordinate (vert)
              to ++(1/2,1/2)
              to ++(-1/2,1/2) coordinate (triv)
              to (0,4);
        \draw[strand a]
        (vert) to ++(-1/2,1/2)
               to ++(1/2,1/2);
        \node[bead] at ($(triv)!2/3!(0,4)$) {$v$};
        \node[longbead] at ($(triv)!1/3!(0,4)$) {$\sigma_1(a)$};
    \end{tikzineqn}
    \quad=\quad
    \begin{tikzineqn}[scale=\lemmascale]
        \coordinate (bot) at (0,-2); 
        \coordinate (mid) at (0,0); 
        \coordinate (top) at (0,2); 
        \coordinate (bead1) at ($(bot)!1/3!(top)$);
        \coordinate (bead2) at ($(bot)!2/3!(top)$);
        \draw[strand m] (top) to (bot) node[below] {$m$};
        \node[bead] at (bead2) {$v$};
        \node[longbead] at (bead1) {$\sigma_1(a)$};
    \end{tikzineqn}
    \]
as desired.
A similar argument using either method can be applied to show that $\sigma_2$ is also real-valued.
\end{proof}

\begin{lemma}\label{lem:RQSigma3Real}
    The function $\sigma_3$ is real-valued.
\end{lemma}

\begin{proof}
Let $a\in A$.
We want to show that $\sigma_3(a)$ is in the center of $\mathbb{H}$.
First, we will use the naturality of the braiding to show that
\[
(\forall v\in \mathbb{H}) \quad
[a]\triangleleft \big(\sigma_3(a)\cdot v\big)
= [a]\triangleleft \big(v\cdot \sigma_3(a)\big)
.\]
First, we use naturality and the property of the trivalent vertex
to get
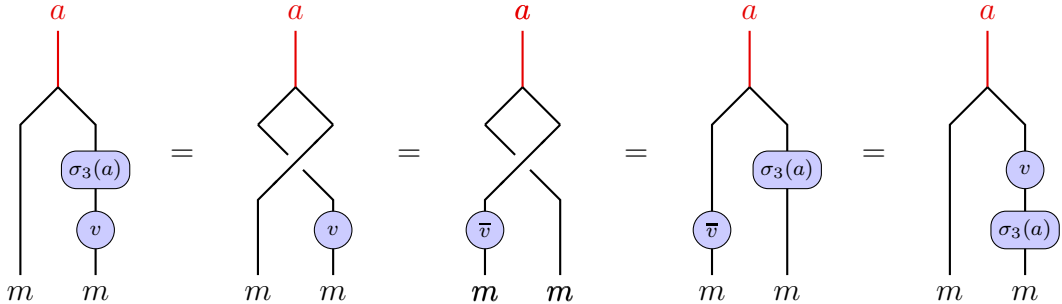
\[
\begin{tikzineqn}[scale=0.5]
    \draw[strand a] (0,0) -- (0,1.5);
    \draw[strand m] (0,0) -- (1,-1) -- ++(0,-4);
    \draw[strand m] (0,0) -- (-1,-1) -- ++(0,-4);
    \node[longbead] at (1,-2.2) {$\sigma_3(a)$};
    \node[bead] at (1,-3.8) {$v$};
    \node[below] at (-1,-5) {$m$};
    \node[below] at (1,-5) {$m$};
    \node[strand a,above] at (0,1.5) {$a$};
\end{tikzineqn}
\quad=\quad
\begin{tikzineqn}[scale=0.5]
    \draw[strand a] (0,0) -- (0,1.5);
    \draw[strand m] (0,0) -- (1,-1);
    \draw[strand m] (0,0) -- (-1,-1);
    \begin{knot}[clip width = 10]
        \strand[strand m] (1,-1) -- ++(-2,-2) -- ++(0,-2);
        \strand[strand m] (-1,-1) -- ++(2,-2) -- ++(0,-2);
    \end{knot}
    \node[bead] at (1,-3.8) {$v$};
    \node[below] at (-1,-5) {$m$};
    \node[below] at (1,-5) {$m$};
    \node[strand a,above] at (0,1.5) {$a$};
\end{tikzineqn}
\quad=\quad
\begin{tikzineqn}[scale=0.5]
    \draw[strand a] (0,0) -- (0,1.5);
    \draw[strand m] (0,0) -- (1,-1);
    \draw[strand m] (0,0) -- (-1,-1);
    \begin{knot}[clip width = 10]
        \strand[strand m] (1,-1) -- ++(-2,-2) -- ++(0,-2);
        \strand[strand m] (-1,-1) -- ++(2,-2) -- ++(0,-2);
    \end{knot}
    \node[bead] at (-1,-3.8) {$\overline{v}$};
    \node[below] at (-1,-5) {$m$};
    \node[below] at (1,-5) {$m$};
    \node[strand a,above] at (0,1.5) {$a$};
    \node at (-1,-5.5) {$m$};
    \node at (1,-5.5) {$m$};
    \node[strand a] at (0,2) {$a$};
\end{tikzineqn}
\quad=\quad
\begin{tikzineqn}[scale=0.5]
    \draw[strand a] (0,0) -- (0,1.5);
    \draw[strand m] (0,0) -- (1,-1) -- ++(0,-4);
    \draw[strand m] (0,0) -- (-1,-1) -- ++(0,-4);
    \node[longbead] at (1,-2.2) {$\sigma_3(a)$};
    \node[bead] at (-1,-3.8) {$\overline{v}$};
    \node[below] at (-1,-5) {$m$};
    \node[below] at (1,-5) {$m$};
    \node[strand a,above] at (0,1.5) {$a$};
\end{tikzineqn}
\quad=\quad
\begin{tikzineqn}[scale=0.5]
    \draw[strand a] (0,0) -- (0,1.5);
    \draw[strand m] (0,0) -- (1,-1) -- ++(0,-4);
    \draw[strand m] (0,0) -- (-1,-1) -- ++(0,-4);
    \node[bead] at (1,-2.2) {$v$};
    \node[longbead] at (1,-3.8) {$\sigma_3(a)$};
    \node[below] at (-1,-5) {$m$};
    \node[below] at (1,-5) {$m$};
    \node[strand a,above] at (0,1.5) {$a$};
\end{tikzineqn}
\]
By self duality of $m$, we may ``rotate" the diagram up to a non-zero
quaternionic constant by composing with the coevaluation map on the left strand, yielding
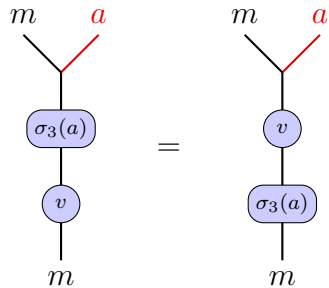
\[
\begin{tikzineqn}[scale=0.5]
    \draw[strand a] (0,0) -- (1,1) node[above] {$a$};
    \draw[strand m] (0,0) -- (-1,1) node[above] {$m$};
    \draw[strand m] (0,0) -- (0,-5) node[below] {$m$};
    \node[longbead] at (0,-1.5) {$\sigma_3(a)$};
    \node[bead] at (0,-3.5) {$v$};
\end{tikzineqn}
\quad=\quad
\begin{tikzineqn}[scale=0.5]
    \draw[strand a] (0,0) -- (1,1) node[above] {$a$};
    \draw[strand m] (0,0) -- (-1,1) node[above] {$m$};
    \draw[strand m] (0,0) -- (0,-5) node[below] {$m$};
    \node[longbead] at (0,-3.5) {$\sigma_3(a)$};
    \node[bead] at (0,-1.5) {$v$};
\end{tikzineqn}
\]
which we may compose with the inverse to the trivalent vertex
to conclude the desired result.
\end{proof}

\subsection{The Hexagon Equations}
Since all the braiding coefficients are real, the only difference in
the braiding equations arises from the fact that
$m\otimes m\cong 4\bigoplus_{a\in A} a$ 
rather than $\bigoplus_{a\in A} a$.
The graphical computations remain mostly the same except for
the hexagon diagrams involving $\alpha_{m,m,m}$.
The resulting braiding equations are equations
(\ref{eqn:hexR1}) through (\ref{eqn:hexR7}),
(\ref{eqn:hexR9}) through (\ref{eqn:hexR15}),
and the following two, which differ from 
(\ref{eqn:hexR8}) and (\ref{eqn:hexR16}) by a coefficient of $-2$:
\begin{equation}
   \sigma_3(a)\tau\chi(a,b)^{-1}\sigma_3(b)
       = -2\sum_{c\in A}\tau\chi(a,c)^{-1}\sigma_2(c)\tau\chi(c,b)^{-1},
       \tag{10'}\label{eqn:hexH8}
\end{equation}
\begin{equation}
    \sigma_3(a)\tau \chi(a,b)\sigma_3(b)
        = -2\sum_{c\in A} \tau \chi(a,c)\sigma_1(c)\tau \chi(c,b).
        \tag{18'}\label{eqn:hexH16}
\end{equation}

The presence of the $-2$ does not affect the algebraic reduction process, and the reduced hexagon equations are thus
\begin{align}
	&\sigma_0(a,b) = \chi(a,b), \label{eqn:RQreducedR1} \\
	&\sigma_1(a)^2 = \chi(a,a), \label{eqn:RQreducedR2} \\
	&\sigma_1(ab)  = \sigma_1(a)\sigma_1(b)\chi(a,b), \label{eqn:RQreducedR3} \\
	&\sigma_2(a)   = \sigma_1(a), \label{eqn:RQreducedR4} \\
	&\sigma_3(1)^2 = -2\tau \sum_{c\in A}\sigma_1(c), \label{eqn:RQreducedR5} \\
	&\sigma_3(a)   = \sigma_3(1)\sigma_1(a)\chi(a,a), \label{eqn:RQreducedR6}
\end{align}
which coincide with (\ref{eqn:reducedR1}) through (\ref{eqn:reducedR6})
except for the added $-2$ in (\ref{eqn:RQreducedR5}).

\subsection{Classification}
With the notation of Proposition \ref{prop:OrbitEquivalenceCharacterization}, we have: 
\begin{theorem} \label{thm:RQ-class-sols}
	Braidings on $\mathcal{C}_{\mathbb{H}}(K_4^n, h^{n}, \tau)$ are in bijection with $\QF^n_{-\text{sgn}(\tau)}\times \{\pm 1\}$.
\end{theorem}
\begin{proof} The argument is exactly parallel to the proof of Theorem \ref{thm:split-class-sols}, except that the extra factor of $-2$ in \eqref{eqn:RQreducedR5} gives $\sgn(\sigma_1) = -\sgn(\tau)$. 
\end{proof}
\begin{theorem}
A real/quaternionic Tambara-Yamagami category $\mathcal{C}_{\mathbb{H}}(A, \chi, \tau)$ admits a braiding if and only if either $(A, \chi) \cong (K_4^n, h^{n})$ for $n > 0$ or $(A, \chi)$ is trivial and $\tau < 0$.
\end{theorem}
\begin{proof}
By Theorem \ref{thm:WallClassification}, we know $(A, \chi) \cong (K_4^n, h^{n})$. The conclusion then follows from the previous theorem, observing that $\QF^n_{-\text{sgn}(\tau)}$ is always nonempty except when $n = 0$ and $\tau > 0$. 
\end{proof}
Since the group $K_4^n$, bicharacter $h^{\oplus n}$ and scaling coefficient $\tau$ are determined by context, we denote the braiding on $\mathcal{C}_{\mathbb{H}}(K_4^n, h^{n}, \tau)$ corresponding to $(\sigma, \epsilon) \in \QF^n_{-\text{sgn}(\tau)} \times \{\pm 1\}$ by $\mathcal{C}_{\mathbb{H}}(\sigma_{1}, \epsilon)$.
\begin{definition}\label{defn:ExplicitRealQuaternionicBraidings}
	Given an element $(\sigma, \epsilon)$ of $\QF_{-\sgn\tau}\times \{\pm 1\}$, we define a braided structure $\mathcal{C}_\mathbb{H}(\sigma,\epsilon)$ on $\mathcal{C}_\mathbb{H}(K_4^n,h^{n},\tau)$ by: 
	\begin{align*}
		\beta_{a,b} &= \chi(a,b)\cdot \id_{ab}, \\
		\beta_{a,m} &= \beta_{m,a} = \sigma(a)\cdot \id_{m}, \\
		\beta_{m,m} = 
        \sum_{\substack{s\in S\\a\in K_4^n}}
		\epsilon\,&\sigma(a) 
        (\id_m \otimes \bar{s})[a]^{\dag}[a]
        (s \otimes \id_m).
	\end{align*}
\end{definition}
As before, we now turn to the question of when $\mathcal{C}_\mathbb{H}(\sigma,\epsilon)$ and $\mathcal{C}_\mathbb{H}(\sigma',\epsilon')$ are braided equivalent. 
\begin{definition}
	Let $f \in \Aut(A, \chi)$ and $\kappa \in \{\pm1\}$. We let $F(f,\kappa)$ be the monoidal endofunctor of $\mathcal{C}_\mathbb{H}(K_4^n,h^{n},\tau)$ whose underlying action on grouplike simples is $f$ and fixes $m$ and $\End(m)$. The tensorator coefficients are:
	
	$$J_{a,b} = \id_{f(a)f(b)}, \quad J_{a,m} = \id_{f(a)} \otimes \id_m, \quad J_{m,a} = \id_m \otimes \id_{f(a)}, \quad J_{m,m} = \kappa\cdot\id_m \otimes \id_m.$$
\end{definition}
\begin{lemma}\label{lem:RealQuaternionicFunctorClassification} For any $A,\chi, \tau$,
	$$\pi_0\Aut_\otimes\big(\mathcal{C}_\mathbb{H}(A,\chi,\tau)\big) \cong \Aut(A, \chi) \times \mathbb{Z}/2\mathbb{Z},$$ with representatives given by $F(f,\kappa)$. 
\end{lemma}
\begin{proof} 
We first remark that every functor in $\Aut(\mathcal{C}_\mathbb{H}(A, \chi,\tau))$ is naturally equivalent to one which fixes $\End(m)$; the action of $F$ on $\End(m)$ must be conjugation by some quaternion, and this same quaternion forms the requisite natural transformation together with the identity on the invertible objects. 

Let $\psi$ and $\omega$ be functions $A \to \mathbb{R}^\times$ with $\phi(a)\omega(a)$ constant.  We define $F(f, \psi, \omega)$ to be the monoidal functor whose underlying homomorphism is $f$ and has
\begin{align*}
J_{a,b} = \delta \psi(a,b) \cdot \id_{f(a)f(b)}, &\quad J_{a,m} = \psi(a)\cdot \id_{f(a)} \otimes \id_m,\\ \quad J_{m,a} = \psi(a)\cdot \id_m \otimes \id_{f(a)}, &\quad J_{m,m} =  \id_m \otimes \omega(a)\id_m.
\end{align*}

The proof of Theorem 5.4 of \cite{pss23} shows us that $F(f, \psi, \omega)$ is a monoidal functor and every monoidal functor with underlying homomorphism $f$ is monoidally isomorphic to $F(f, \psi, \omega)$ for some $\psi, \omega$. 

The consistency equations for a monoidal natural isomorphism $\mu \colon F(f, \psi, \omega) \to F(f, \psi', \omega')$ are: 
\begin{align*}
\phi'(a) &= \phi(a)\mu_a \\
\omega'(a) &= \frac{\overline{\mu_m}\mu_m}{\mu_a}\omega(a)
\end{align*} 
By setting $\mu_a = \phi(a)^{-1}$, and using that $\phi(a)\omega(a)$ is constant, we see that $\mu$ defines a natural isomorphism to $F(f, \sgn(\omega(1)))$. 

Moreover, these same consistency conditions rule out any natural isomorphisms $F(f, 1) \to F(f,-1)$; we must have $\mu_1 = 1$ and so would obtain $-1 = |\mu_m|^2$, a contradiction. 
\end{proof}
The proofs of the following proposition and theorem are identical to those of Proposition \ref{prop:RealFunctorBraided} and Theorem \ref{thm:SplitCaseEquivalence} upon replacing Lemma \ref{lem:SplitRealFunctorClassification} with Lemma \ref{lem:RealQuaternionicFunctorClassification}.
\begin{proposition}\label{prop:QuaternionincFunctorBraided}
	The monoidal functor $F(f, \kappa)$ forms a braided monoidal equivalence $\mathcal{C}_\mathbb{H}(\sigma,\epsilon) \to \mathcal{C}_\mathbb{H}(\sigma',\epsilon')$ if and only if $f \cdot \sigma = \sigma'$ and $\epsilon = \epsilon'$.
\end{proposition}
\begin{theorem}\label{thm:RealQuaternionicBraidedEquivalence}	There is a braided monoidal equivalence $\mathcal{C}_\mathbb{H}(\sigma,\epsilon) \sim \mathcal{C}_\mathbb{H}(\sigma',\epsilon')$ if and only if $\epsilon = \epsilon'$. In particular, there is no braiding on $\mathcal{C}_\mathbb{H}(K_4^n,h^{\oplus n},\tau)$ when $n = 0$ and $\tau > 0$, and in all other cases there are exactly two equivalence classes of braidings. 
\end{theorem}
\begin{remark}
In the split real case, the $\Aut(A, \chi)$ orbit which extends to a braiding has the same sign as $\tau$. Here, the sign is reversed. In both cases the scalar $\sigma_3(1)$ is a braided invariant, and indeed determines the equivalence class. 
\end{remark}

\begin{example}\label{eg:Q+HasNoBraiding}
    Let $\mathcal Q_{\pm}:=\mathcal C_{\mathbb H}(K_4^0,h^{\oplus0},\pm\tfrac12)$.  It can be shown by direct computation\footnote{The direct computation referenced here is analogous to our analysis of hexagons, but where only forward hexagons are analyzed for the sake of finding half-braidings instead of full braidings.} that as a fusion category, $\mathcal Z(\mathcal Q_+)\simeq\mathcal C_{\mathbb C}(\mathbb Z/2\mathbb Z,\id_{\mathbb C},\textit{triv}\,,\tfrac12)$.  
    In particular, $\mathcal Z(\mathcal Q_+)$ contains no quaternionic object, and therefore cannot contain $\mathcal Q_+$ as a fusion subcategory.  This is equivalent to the observation that $\mathcal Q_+$ cannot have a braiding, as indicated by Theorem \ref{thm:RealQuaternionicBraidedEquivalence}.
    This is directly analogous to the fact that $\mathcal{C}_{\mathbb{R}}(K_4^0,h^{\oplus 0},-1)$ also admits no braiding.

    Here is yet another way to see why there cannot be a braiding in this case.
    The category $\mathcal Q_+$ can be realized as the time reversal equivariantization of $\Vect_{\mathbb C}^\omega(\mathbb Z/2\mathbb Z)$, where $0\neq[\omega]\in H^3(\mathbb Z/2\mathbb Z;\mathbb C^\times)$ (see \cite{MR2946231} for further details on categorical Galois descent).
    The time reversal symmetry that produces $\mathcal Q_+$ is anomalous in the sense that it uses a nontrivial tensorator $T_1\circ T_1\cong T_0=\id$.
    This anomaly is what causes the presence of a quaternionic object, because without it, equivariantization would just produce $\Vect_{\mathbb R}^\omega(\mathbb Z/2\mathbb Z)$.  If $\mathcal Q_+$ were to admit a braiding, then by base extension it would produce one of the two braidings on the category $\Vect_{\mathbb C}^\omega(\mathbb Z/2\mathbb Z)$ \textemdash~ either the semion or reverse semion.
    However, the time reversal functor $T_1$ is not braided (it swaps these two braidings), and so neither of these braidings could have come from $\mathcal Q_+$.
\end{example}

Taking $\sigma = \sigma'$ and $\epsilon = \epsilon'$ in Proposition \ref{prop:QuaternionincFunctorBraided}, we obtain:
\begin{corollary}
$$\pi_0\Aut_{br}\big(\mathcal{C}_{\mathbb{H}}(K_4^n , h^{\oplus n}, \tau, \sigma, \epsilon)\big) \cong \GO_{\sgn(\sigma)}^n \times \mathbb{Z}/2\mathbb{Z}$$
\end{corollary}

\begin{lemma}
	There are exactly two families of twist morphisms for any $\mathcal{C}_{\mathbb{H}}(\sigma, \epsilon)$, corresponding to a sign $\rho \in \{\pm 1\}$. These twists are ribbon structures.
\end{lemma}
\begin{proof}
	 Denoting the components of the twist by $\theta_x$, the required equations can be derived identically to \cite[Section 3.7]{sie00}, and algebraically reduced in an identical way using that $\mathbb{H}$ is a division algebra and $\sigma$ is real valued and so the values $\sigma(a)$ commute with $\theta_m$. The results are (still): 
	 
	 \begin{align*}
	 	\theta_{ab}& = \theta_a\theta_b\\
	 	\theta_a &= \sigma(a)^2 = 1\\
	 	\theta_a &= \theta_m^2\sigma_3(a)^2
	 	\end{align*}
	 
	Thus, the square root required to define $\theta_m$ is always of a positive real number and therefore still determined by a sign. Since every simple object is self dual, the required axiom is simply $\theta_m = \theta_m^*$. But this holds as a result of the (real) linearity of composition.
\end{proof}

\section{Braidings on Real/Complex Tambara-Yamagami Categories}\label{sec:Real/Complex}

In the case where the invertibles are real and $m$ is complex, the analysis in \cite{pss23} was much more involved than in the other cases. Part of this complexity arises due to the fact that $m$ can be either directly or conjugately self dual, and this property is a monoidal invariant, necessitating some degree of casework.

The classification of split real Tambara-Yamagami categories relies on some background on (abelian) groups and nondegenerate bilinear forms; we begin this section with the analogous background for the real/complex Tambara-Yamagami categories. If $A$ is an abelian group, we let $D_A$ denote the generalized dihedral group on $A$, the semidirect product of $A$ by $\mathbb{Z}/2\mathbb{Z}$ (where the nontrivial element, which we will denote by $w$, acts by inversion).   

Identifying $\text{Gal}(\mathbb C/\mathbb R)$ with $\mathbb{Z}/2\mathbb{Z}$, the group $D_A$ acts on $\mathbb{C}$ by Galois automorphisms by way of the quotient map $\pi_A \colon D_A \to \mathbb{Z}/2\mathbb{Z}$. We will use superscripts to denote this action, i.e., for $\lambda\in \mathbb{C}$ we define $\lambda^x = [\pi_A(x)](\lambda)$. 

For each $g\in \mathbb{Z}/2\mathbb{Z} = \text{Gal}(\mathbb C/\mathbb R)$, its \textit{degree}, denoted $|g|$ is defined as:
\[
    |g|:= \begin{cases}
        0, \ g \text{ is trivial}, \\
        1, \ g \text{ is nontrivial}.
    \end{cases}\;
\] We now extend $|\cdot|$ to $x \in D_A$ by the formula $$|x| := |\pi_A(x)|.$$  We will also occasionally use  products of elements $x,y \in D_A$ and elements $g \in \text{Gal}(\mathbb C/\mathbb R)$, where we declare $gxy := g\pi_A(xy)$ and extend the degree notation accordingly.  

Following \cite{pss23}, a \textit{bicocycle} is a function $\chi \colon D_A \times D_A \to \mathbb{C}^\times$ such that 
$$\chi(x, yz) = \chi(x,y)\chi(x,z)^y, \quad \text{and} \quad \chi(xy,z) = \chi(x,z)^y\chi(y,z).$$ A bicocycle said to be \textit{symmetric} with respect to $g \in \text{Gal}(\mathbb C/\mathbb R)$ if $$\chi(x,y) = \chi(y,x)^{gxy}.$$
We can now recall the classification for associators: 

\begin{theorem}[{\cite[Theorem 6.10]{pss23}}]\label{thm:RealComplexFromPSS}
Let $\tau=\sfrac{\pm 1}{\sqrt{2|A|}}$, let $(-)^g\in\text{Gal}(\mathbb C/\mathbb R)$, and let $\chi:D_A\times D_A\to \mathbb C^\times$ be a symmetric bicocycle on $D_A$ with respect to $(-)^g$, whose restriction $\chi\mid_{A \times A}$ is a nondegenerate bicharacter.
A quadruple of such data gives rise to a non-split Tambara-Yamagami category $\mathcal{C}_{\bb C}(D_A,g,\chi,\tau)$, with $\End(\mathbbm{1})\cong\mathbb{R}$ and $\End(m)\cong\mathbb{C}$.
Furthermore, all equivalence classes of such categories arise in this way.
More explicitly, two categories $\mathcal{C}_{\bb C}(D_A,g,\chi,\tau)$ and $\mathcal{C}_{\mathbb{C}}(D_{A'},g',\chi',\tau')$ are equivalent if and only if $g=g'$, and there exists
the following data:
    \begin{enumerate}[label = \roman*)]
        \item an isomorphism $f:D_A\to D_{A'}$,
        \item a map $(-)^h:\mathbb{C}\to\mathbb{C}$, either the identity or complex conjugation,
        \item a scalar $\lambda\in S^1\subset \mathbb C$,
    \end{enumerate}
    satisfying the following conditions for all $a,b\in D_A$
    \begin{gather}
        \chi'\big(f(a),f(b)\big)=\frac{\lambda\cdot\lambda^{ab}}{\lambda^a\cdot\lambda^b}\cdot\chi(a,b)^h\;,\label{EquivCond1}\\
        \frac{\tau'}{\tau}=\frac{\lambda}{\lambda^g}\label{EquivCond2}\,.
    \end{gather}

The components of the associator for $\s C$ are defined by the following equations
    \begin{gather*}
	    \alpha_{a,b,c}=\id_{abc},\\
	    \alpha_{a,b,m}=\alpha_{m,a,b}\;=\;\id_m,\\
	    \alpha_{a,m,b}=\chi(a,b)^{ab}\cdot\id_m,\\
	   % \alpha_{a,b,m}=1\\
	    \alpha_{m,m,a}=\alpha_{a,m,m}\;=\;\id_{m\otimes m},\\
	    \alpha_{m,a,m}=\sum_{\substack{b\in D_A\\t\in\{1,i\}}}\left(\id_m\otimes \chi(a,b)^bt^b\right)\big(\iota_{b,1}\pi_{b,t}\big),\\
	   % \alpha_{m,m,a}=\id_{m\otimes m},\\
		\alpha_{m,m,m}=\sum_{\substack{a,b\in D_A\\s,t\in\{1,i\}}}(\id_m\otimes\iota_{a,t})\circ\left(\frac{\overline{s}^{gab}t^{b}\tau}{\chi(a,b)^g}\right)\circ(\pi_{b,s}\otimes\id_m)\,.
	\end{gather*}
where
$\pi_{b,s}=[b](\id_m\otimes \overline{s})$ denotes
the projection maps
$m\otimes m\to a$,
and $\iota_{a,t}=(\id_m\otimes t)[a]^\dagger$
the inclusion maps.
  The left and right unitors $\ell_X$ and $r_X$ are identities for all simple objects $X$.
\end{theorem}

We continue to recall conventions and facts about
$\mathcal{C}_{\mathbb{C}}(D_A,g,\tau,\chi)$ 
from \cite[Section 6]{pss23}. The group of invertible objects in
$\mathcal{C}_{\mathbb{C}}(D_A,g,\tau,\chi)$ 
is $D_A$.
The spaces
\[
\Hom(x\otimes m,m),\quad
\Hom(m\otimes x,m) \ \text{and}\
\Hom(m\otimes m,x)
\]
are 1-dimensional complex bimodules,
which are isomorphic to the trivial bimodule
$\mathbb{C}$ if $|x|=0$
or the conjugating bimodule $\overline{\mathbb{C}}$ 
if $|x|=1$ in which the left and right actions
differ by conjugation.
In practice, when a scalar $\lambda\in \mathbb{C}$
passes through a trivalent vertex in graphical calculus (representing a chosen basis vector for the Hom space), 
it ``conjugates by $x$",
i.e. we have the relations
\begin{equation}\label{eqn:passThroughVertexConjugate}
    \begin{tikzineqn}[scale=\beadedTSize]
        \TrivalentVertex{x}{m}{m} 
        \DrawSmallBead{mid}{top}{\lambda}
    \end{tikzineqn}
    \ = \
    \begin{tikzineqn}[scale=\beadedTSize]
        \TrivalentVertex{x}{m}{m} 
        \DrawSmallBead{mid}{bottom right}{\lambda^{x}}
    \end{tikzineqn}
    \quad\quad
    \begin{tikzineqn}[scale=\beadedTSize]
        \TrivalentVertex{m}{x}{m} 
        \DrawSmallBead{mid}{top}{\lambda}
    \end{tikzineqn}
    \ = \
    \begin{tikzineqn}[scale=\beadedTSize]
        \TrivalentVertex{m}{x}{m} 
        \DrawSmallBead{mid}{bottom left}{\lambda^{x}}
    \end{tikzineqn}\;\,
\end{equation}
\begin{equation}\label{eqn:passThroughVertexConjugateMM}
    \begin{tikzineqn}[scale=\beadedTSize]
        \TrivalentVertex{m}{m}{x} 
        \DrawSmallBead{mid}{bottom left}{\lambda}
    \end{tikzineqn}
    \ = \
    \begin{tikzineqn}[scale=\beadedTSize]
        \TrivalentVertex{m}{m}{x} 
        \DrawSmallBead{mid}{bottom right}{\lambda^{gx}}
    \end{tikzineqn}
    \;\;.
\end{equation}

We begin our analysis of braidings on
$\mathcal{C}_{\mathbb{C}}(D_A,g,\chi,\tau)$
by imposing some conditions on $A$ and $\chi$:
\begin{lemma}\label{lem:RCChiProperties}
Suppose $\mathcal{C}_{\mathbb{C}}(D_A,g,\tau,\chi)$
admits a braiding,
with $D_A\cong A \rtimes (\mathbb{Z}/2\mathbb{Z})\langle w \rangle$.
Then, $A \cong (\mathbb{Z}/2\mathbb{Z})^{n}$
is an elementary abelian 2-group with
$n\in \mathbb{Z}_{\ge 0}$,
and the symmetric bicocycle $\chi$ satisfies
the following:
\begin{enumerate}[label=(\roman*)]
\item For all $a\in A$ and all $x\in D_A$,
    $\chi(a,x)$ is real-valued;
\item $\chi$ is symmetric;
\item $\chi(x,y)=\chi(x,y)^{gxy}=\chi(x,y)^{g}$ for all $x,y\in D_A$.
\end{enumerate}
\end{lemma}
\begin{proof}
If $\mathcal{C}_{\mathbb{C}}(D_A,g,\tau,\chi)$ admits a braiding,
then $D_A$ is an abelian generalized dihedral group,
so for any $x\in D_A$ we have
\[
    x=ww^{-1}x=wxw^{-1}=x^{-1} \implies x^2=1.
\]
Now we use the cocycle condition to see that
for all $x\in D_A$,
\[
    \chi(1,x)=\chi(1,x)^2 \implies \chi(1,x)=1,
\]
and by the same argument in the other coordinate
we have $\chi(x,1)=1$.
Then, since $a^2=1$, we have
\[
    1=\chi(a^2,x)=\chi(a,x)^{a}\chi(a,x)=\chi(a,x)^2,
\]
which tells us that $\chi(a,x)\in \{\pm 1\}$
(and similarly $\chi(x,a)\in \{\pm 1\}$).
This implies that $\chi$ is real valued and therefore symmetric on $(D_A\times A)\cup (A\times D_A)$ as a consequence of the symmetric cocyle condition.

For condition (ii), we check that for any $a,b\in A$,
\begin{align*}
\chi(aw,bw)&=\chi(a,bw)^{w}\chi(w,bw) \\
&=\chi(a,b)\chi(a,w)^{b}\chi(w,b)\chi(w,w)^{b}\\
&=\chi(a,b)\chi(a,w)\chi(w,b)\chi(w,w),
\end{align*}
which gives us symmetry of $\chi$.
Note that in particular $\chi(aw,aw)=\chi(a,a)\chi(w,w)$.

It suffices to check condition (iii)
on $Aw\times Aw$, since $\chi$ is real-valued on the remainder of $D_A \times D_A$.
We use the symmetric cocycle and symmetric
conditions to get that $\chi(x,y)=\chi(x,y)^{gxy}$,
and since $|xy|=0$ we have the desired result.
\end{proof}

At this point, we have been using a choice of isomorphism $D_A\cong A\rtimes (\mathbb{Z}/2\mathbb{Z})\langle w \rangle$, which amounts to choosing an element $w\in D_A\setminus A$.  It turns out that there is a canonical way to choose this element.

\begin{lemma}\label{lem:CanonicalW}
    There is a unique $w\in D_A\setminus A$ with the property that $\chi(w,-)$ is trivial when restricted to $A$.	Moreover restriction to $A$ gives an isomorphism $\Aut(D_A, \chi)$ to $\Aut(A, \chi|_{A \times A})$.
\end{lemma}

\begin{proof}
    At first, let $w'\in D_A\setminus A$ be any element.  Since $\chi_{A\times A}$ is nondegenerate, there exists a unique $c\in A$ such that $\chi(w',a)=\chi(c,a)$ for every $a\in A$.  It follows that $w=cw'\in D_A\setminus A$ is an element that satisfies
    \[\chi(w,a)=\chi(c,a)\chi(w',a)=\chi(w',a)^2=1\,,\]
    where the last equality follows from Lemma \ref{lem:RCChiProperties} parts (i) and (ii).

    Any other choice is of the form $bw$ for $b\in A$, and would satisfy $\chi(bw,a)=\chi(b,a)\chi(w,a)=\chi(b,a)$ for every $a\in A$.  Again by nondegeneracy, $\chi(bw,-)$ can only be trivial when $b=1$, so $w$ is unique. For the second part of the lemma, the defining property of $w$ implies $w$ is fixed by every $f \in \Aut(D_A,\chi)$, so that $f$ is completely determined by the homomorphism property together with its restriction to $A$. 
\end{proof}

\begin{lemma} \label{lem:RCChiWWPositive}
Up to monoidal equivalence, $\chi(w,w)$ can be taken to be 1 when $|g|=0$. 
\end{lemma}
\begin{proof}
By Theorem \ref{thm:RealComplexFromPSS}, for any $\lambda\in S^1\subset\mathbb C^\times$ there exists an equivalence $(\id_{\mathcal C},\id_{\mathbb C},\lambda):\mathcal C_{\mathbb C}(D_A,\id,\chi,\tau)\to\mathcal C_{\mathbb C}(D_A,\id,\chi',\tau)$, where $\chi'$ is the bicocycle defined by the equation
\[\chi'(a,b)=\frac{\lambda\cdot\lambda^{ab}}{\lambda^a\cdot\lambda^b}\cdot\chi(a,b)\,.\]
Whenever $|a|=0$ or $|b|=0$, it follows that $\chi'=\chi$.  When both arguments conjugate, the bicocycles are related by $\chi'=\lambda^4\chi$.  In particular, by setting $\lambda^4=\chi(w,w)^{-1}$, we can force $\chi'(w,w)=1$.
\end{proof}

\subsection{Hexagon Equations}
As before, we will use the Yoneda embedding and a set of normal bases to produce equations. Our conventions for the functions $\sigma_i$ are similar to the previous sections: 
 \[
        \begin{tikzineqn}[scale=\eqnscale]
            \draw[strand ab] (0,0) to ++(0,1) node[above] {$ab$};
            \begin{knot}[clip width=10]
                \strand[strand a] (0,0)
                to ++(1,-1)
                to ++(-2,-2) node[below left] {$a$};
                \strand[strand b] (0,0)
                to ++(-1,-1)
                to ++(2,-2) node[below right,yshift=0.1cm] {$b$};
            \end{knot}
        \end{tikzineqn}
        := \
        \sigma_0(a,b)
        \begin{tikzineqn}[scale=\tscale]
            \coordinate (top)            at (0,1);
            \coordinate (bottom left)   at (-1,-1);
            \coordinate (bottom right)    at (1,-1);
            \draw[strand a]  (0,0)  to (bottom left)
                node[below left] {$a$};
            \draw[strand b]  (0,0)  to (bottom right)
                node[below right, yshift=0.1cm] {$b$};
            \draw[strand ab] (0,0)  to (top)
                node[above] {$ab$};
        \end{tikzineqn}
        \quad\quad
        \begin{tikzineqn}[scale=\eqnscale]
            \draw[strand m] (0,0) to ++(0,1) node[above] {$m$};
            \begin{knot}[clip width=10]
                \strand[strand a] (0,0)
                to ++(1,-1)
                to ++(-2,-2) node[below left] {$a$};
                \strand[strand m] (0,0)
                to ++(-1,-1)
                to ++(2,-2) node[below right] {$m$};
            \end{knot}
        \end{tikzineqn}
        := \
        \begin{tikzineqn}[scale=\tscale] 
            \TrivalentVertex{a}{m}{m} 
            \DrawLongBead{mid}{bottom right}{\sigma_1(a)}
        \end{tikzineqn} 
    \]
    \vspace{-0.2cm}
    \[
        \begin{tikzineqn}[scale=\eqnscale]
            \draw[strand m] (0,0) to ++(0,1) node[above] {$m$};
            \begin{knot}[clip width=10]
                \strand[strand m] (0,0)
                to ++(1,-1)
                to ++(-2,-2) node[below left] {$m$};
                \strand[strand a] (0,0)
                to ++(-1,-1)
                to ++(2,-2) node[below right] {$a$};
            \end{knot}
        \end{tikzineqn}
        := \
        \begin{tikzineqn}[scale=\tscale] 
            \TrivalentVertex{m}{a}{m} 
            \DrawLongBead{mid}{bottom left}{\sigma_2(a)}
        \end{tikzineqn} \quad\quad \
        \begin{tikzineqn}[scale=\eqnscale]
            \draw[strand a] (0,0) to ++(0,1) node[above] {$a$};
            \begin{knot}[clip width=10]
                \strand[strand m] (0,0)
                to ++(1,-1)
                to ++(-2,-2) node[below left] {$m$};
                \strand[strand m] (0,0)
                to ++(-1,-1)
                to ++(2,-2) node[below right] {$m$};
            \end{knot}
        \end{tikzineqn}
        := \
        \begin{tikzineqn}[scale=\tscale] 
            \TrivalentVertex{m}{m}{a} 
            \DrawLongBead{mid}{bottom right}{\sigma_3(a)}
        \end{tikzineqn} .
    \]
Note that in this case, the functions $\sigma_i$ for $i > 0$ are \textit{a priori} complex valued. From the graphical calculus computations
(see Appendix \ref{sec:appendixComputationEx} for
an example),
we get the following
equations from the forward hexagon diagrams:
\begin{align}
  &\sigma_0(x,y)\sigma_0(x,z)=\sigma_0(x,yz)
  \label{RCHexagon1} \\
  &\sigma_1(x)\sigma_0(x,y)=\chi(y,x)\sigma_1(x)^y
  \label{RCHexagon2} \\
  &\sigma_0(x,y)\sigma_1(x)=\sigma_1(x)^y\chi(x,y)
  \label{RCHexagon3} \\
  &\sigma_2(y)\chi(x,y)\sigma_2(x)=\sigma_2(xy)
  \label{RCHexagon4} \\
  &\chi(x,y)^y\sigma_1(x)^{gxy}\sigma_1(x)=\sigma_0(x,xy)  
  \label{RCHexagon5} \\
  &\sigma_2(x)^{gxy}\sigma_3(xy)=\sigma_3(y)^x\chi(x,y)^y 
  \label{RCHexagon6} \\
  &\sigma_3(xy)\sigma_2(x)^{gxy} =\sigma_3(y)^x\chi(x,y)^{gx} 
  \label{RCHexagon7} \\
  &\chi(x,y)^{-g}\sigma_3(x)^y\sigma_3(y)^x
  =2\tau\sum_{|z|=|gxy|}\chi(x,z)^{-g}
  \chi(z,y)^{-g}\sigma_2(z)^z
  \label{RCHexagon8}
\end{align}
and the following from the backward hexagon diagrams:
\begin{align}
  &\sigma_0(xy,z)=\sigma_0(x,z)\sigma_0(y,z)
  \label{RCHexagon9}\\
  &\sigma_1(xy)=\sigma_1(x)\sigma_1(y)\chi(x,y)^{-1}
  \label{RCHexagon10}\\
  &\sigma_2(y)^x\chi(x,y)^{-1}=\sigma_0(x,y)\sigma_2(y)
  \label{RCHexagon11}\\
  &\sigma_2(y)^x\chi(y,x)^{-1}=\sigma_2(y)\sigma_0(x,y)
  \label{RCHexagon12}\\
  &\sigma_3(y)\chi(x,y)^{-gx}=\sigma_1(x)\sigma_3(xy)
  \label{RCHexagon13}\\
  &\sigma_3(y)\chi(x,y)^{-y}=\sigma_1(x)\sigma_{3}(xy)
  \label{RCHexagon14}\\
  &\sigma_0(xy,x)=\sigma_2(x)^{gxy}\chi(x,y)^{-y}\sigma_2(x)
  \label{RCHexagon15}\\
  &\sigma_3(x)\sigma_3(y)\chi(x,y)^{xy}=2\tau\sum_{|z|=|gxy|}\chi(x,z)^{gz}\chi(z,y)^{gz}\sigma_1(z)
  \label{RCHexagon16}
\end{align}

We first obtain a few useful equations through algebraic simplification.
Evaluating at $y=x$ in \eqref{RCHexagon10} we get
\begin{equation}
     \sigma_1(x)^2=\chi(x,x) \label{RCReduced2}.
\end{equation}
Rearranging \eqref{RCHexagon3} we get
\begin{equation}
\sigma_0(x,y)=\chi(x,y)\frac{\sigma_1(x)^{y}}{\sigma_1(x)},
\label{RCReduced1}
\end{equation}
which we combine with evaluating \eqref{RCHexagon5} at $y=1$ to get
\begin{equation}
\sigma_1(x)^g=\sigma_1(x).
\label{RCReduced3}
\end{equation}
Lastly, evaluating \eqref{RCHexagon16} at $x=y=1$ yields
\begin{equation}
\sigma_3(1)^2=2\tau \sum_{|z|=|g|} \sigma_1(z).
\label{RCReduced6}
\end{equation}
Using these, we will prove a few lemmas which we will use
to reduce the hexagon equations down to a equivalent set of simpler
equations.

\begin{lemma}\label{lem:RCChiAAReal}
For all $a\in A$, we have $\chi(a,a)=1$.
\end{lemma}

\begin{proof}
Using equations (\ref{RCHexagon3}) and
(\ref{RCHexagon11}), we can write
\[
\sigma_0(x,y)
=\chi(x,y)\frac{\sigma_1(x)^{y}}{\sigma_1(x)}
=\chi(x,y)^{-1}\frac{\sigma_2(y)^{x}}{\sigma_2(y)}.
\]
Setting $x=a$ and $y=w$, we get
\[
\chi(a,w)^2
=\frac{\sigma_1(a)}{\sigma_1(a)^{w}}
\cdot \frac{\sigma_2(w)^{a}}{\sigma_2(w)}.
\]
Since $|a|=0$, we have
\[
1=\chi(a,w)^2
=\frac{\sigma_1(a)}{\sigma_1(a)^{w}}
\implies
\sigma_1(a)=\overline{\sigma_1(a)}.
\]
This tells us that $\sigma_1(a)\in \mathbb{R}$,
which gives us that $\chi(a,a)>0$ by
(\ref{RCReduced2}).
\end{proof}

\begin{corollary} \label{cor:RCHyperbolicPairing}
The bicharacter $\chi|_{A\times A}$
is hyperbolic, and thus for some choice of basis for $A$, is equal to the standard hyperbolic pairing $h^{n}$ on $A \cong K_4^{n}$ 
for some $n\in \mathbb{Z}_{\ge 0}$.
\end{corollary}

\begin{corollary} \label{cor:RCSelfPairingis1}
If $\mathcal{C}_{\mathbb{C}}(D_A,g,\tau,\chi)$ admits a braiding, then up to monoidal equivalence, $\chi$ is a real-valued symmetric bicharacter with $\chi(x,x)=1$ for all $x\in D_A$.
\end{corollary}
\begin{proof}
By Lemma \ref{lem:RCChiProperties} and Lemma \ref{lem:RCChiAAReal},
it suffices to check that $\chi(w,w)=1$ and use the cocycle condition.
When $g$ is trivial, this follows from Lemma \ref{lem:RCChiWWPositive}.
When $g$ is nontrivial, this is implied by \eqref{RCReduced2} and \eqref{RCReduced3} which show us that $\chi(w,w)$ is the square
of a real number.
\end{proof}

\begin{remark}\label{rmk:RCSigma1Real}
In particular, this tells us that $\sigma_1$ is always $\{\pm 1\}$-valued
by \eqref{RCReduced2}, and hence that $\sigma_0=\chi$ by \eqref{RCReduced1}.
Note also that $\chi=\chi^{-1}$ is $\{\pm 1\}$-valued, since
$\chi(x,y)^2=\chi(x^2,y)=\chi(1,y)=1$ for all $x,y\in D_A$.
\end{remark}

\begin{remark}
Note that although we know that the restriction of $\chi$ is nondegenerate on $A \times A$, it is necessarily degenerate on $D_A$, thanks to Lemma \ref{lem:CanonicalW}.
Hence the classification results for bilinear forms used previously
to show that certain forms are hyperbolic do not apply here.
\end{remark}

\begin{lemma}\label{lem:RCSigma3Squared1}
The scalar $\sigma_3(1)^2$ is real, and it can be computed by the formula
\[\sigma_3(1)^2=2^{n+1}\tau\sigma_1(w)^{|g|}\sgn(\sigma_1|_{A_0}).\]
Consequently, $\sigma_3(1)^4 = 1$. 
\end{lemma}
\begin{proof}
Recall that we have
\[
\sigma_3(1)^2=2\tau \sum_{|z|=|g|} \sigma_1(z)\,.
\]
from \eqref{RCReduced6}.
When $g$ is nontrivial, each summand is of the form
\[\sigma_1(aw)=\sigma_1(a)\sigma_1(w)\chi(a,w)=\sigma_1(a)\sigma_1(w)\,,\]
for some unique $a\in A$.  After possibly factoring out the term $\sigma_1(w)$, both cases for $g$ then follow from Proposition \ref{prop:OrbitEquivalenceCharacterization}. 
\end{proof}

\begin{corollary}
The function $\sigma_2$ is real-valued on all of $D_A$.
\end{corollary}
\begin{proof}
Comparing \eqref{RCHexagon6} and \eqref{RCHexagon13} at $y=1$
we get 
\begin{equation}
\sigma_2(x)=\sigma_1(x)^{gx}\frac{\sigma_{3}(1)^{g}}{\sigma_3(1)^{gx}}
=\sigma_1(x)\frac{\sigma_{3}(1)^{g}}{\sigma_3(1)^{gx}}.
\end{equation}
By Lemma \ref{lem:RCSigma3Squared1},
$\sigma_{3}(1)$ is purely real or imaginary, so
$\frac{\sigma_{3}(1)^{g}}{\sigma_3(1)^{gx}}\in \{\pm 1\}$.
\end{proof}

In summary, we have:
\begin{proposition} \label{prop:RCBraidingConstraintsFinal}
The braiding coefficients 
$\sigma_0$, $\sigma_1$ and $\sigma_2$ in the real-complex
category admitting a braiding are necessarily real-valued.
The hexagon equations are equivalent to the following:
\begin{align}
& \sigma_0(x,y)=\chi(x,y)
\label{RCVeryReduced1} \\
& \sigma_1(x)^2=\chi(x,x)
\label{RCVeryReduced2} \\
& \sigma_1(xy)=\sigma_1(x)\sigma_1(y)\chi(x,y)
\label{RCVeryReduced3} \\
& \sigma_3(1)^2=2\tau \sum_{|z|=|g|} \sigma_1(z)
\label{RCVeryReduced4} \\
& \sigma_3(x)=\sigma_3(1)\sigma_1(x)
\label{RCVeryReduced5} \\
& \sigma_3(x) = \sigma_3(x)^g 
\label{RCVeryReduced6} \\
& \sigma_2(x)=\sigma_1(x)\frac{\sigma_{3}(1)}{\sigma_3(1)^{x}} 
\label{RCVeryReduced7}
\end{align}
\end{proposition}

\begin{proof}
First, it remains to check that \eqref{RCVeryReduced5},
\eqref{RCVeryReduced6} and \eqref{RCVeryReduced7} follow from the hexagon
equations. The first and last equations follow from setting $y = 1$ in \eqref{RCHexagon14} and \eqref{RCHexagon7}, respectively.  We postpone the derivation of \eqref{RCVeryReduced6}. 

For the converse, we wish to derived the original hexagon equations from
the reduced ones.
We may rewrite \eqref{RCHexagon4} as
\[
\sigma_1(y)\chi(x,y)\sigma_1(x)
\frac{\sigma_3(1)^2}{\sigma_3(1)^{x}\sigma_3(1)^{y}}
\stackrel{?}{=} \sigma_1(xy) \frac{\sigma_{3}(1)}{\sigma_3(1)^{xy}},
\]
and that it holds in each of the cases $|x|=0$, $|y|=0$ and $|x|=|y|=1$
(in the last case using Lemma \ref{lem:RCSigma3Squared1}).
Similarly \eqref{RCHexagon6} and \eqref{RCHexagon7} follow from the fact
that $\sigma_3(1)^2$ is conjugate invariant. The derivation of \eqref{RCHexagon16} is exactly the same as in the split real case. 

The rest, except for \eqref{RCHexagon8}, follow from straightforward algebraic checks. We now show that \eqref{RCHexagon8} is equivalent to \eqref{RCVeryReduced6} in the presence of the other reduced hexagon equations. To begin, we can expand both sides of \eqref{RCHexagon8} using the definition of $\sigma_2$ and $\sigma_3$ and the properties of $\chi$ to arrive at the equivalent form:

\begin{align*}
\chi(x, y)\sigma_3(1)^x\sigma_3(1)^y\sigma_1(x)\sigma_1(y) &= 2\tau \sum_{|z| = |gxy|} \chi(x, z)\chi(z, y) \sigma_1(z) \frac{\sigma_3(1)^{gxy}}{\sigma_3(1)} \\
&\overset{\eqref{RCHexagon16}}{=} \sigma_3(x)\sigma_3(y)\chi(x,y)\frac{\sigma_3(1)^{gxy}}{\sigma_3(1)}
\end{align*}
Canceling terms we arrive at 
$$\sigma_3(1)^x\sigma_3(1)^y = \sigma_3(1)\sigma_3(1)^{gxy}$$
Since $\sigma_3(1)$ is a 4th root of unity, we have $(\sigma_3(1)^x\sigma_3(1)^y)/(\sigma_3(1)\sigma_3(1)^{xy}) = 1$, so that $\sigma_3(1)^{xy}$ is $g$-fixed for all $x, y$, and thus $\sigma_3(1)$ and $\sigma_3(x)$ are as well. 
\end{proof}
\subsection{Classification of Braidings in the Real/Complex Case}
Recalling Corollary \ref{cor:RCHyperbolicPairing}, we know that any real/complex Tambara-Yamagami category admitting a braiding has $D_A \cong K_4^n \rtimes (\mathbb{Z}/2\mathbb{Z})\langle w \rangle$. Moreover, in all cases we can assume $\chi(x,x) = 1$.
\begin{theorem} \label{thm:RCGTrivialBijectionClassification}
	Braidings on $\mathcal{C}_{\mathbb{C}}(K_4^n \rtimes \mathbb{Z}/2\mathbb{Z}, \id, \chi, \tau)$ are in bijection with pairs $(\sigma, \epsilon) \in \QF(\chi) \times \{\pm 1\}$.
\end{theorem}
\begin{proof}
In this case, since $g = \id$ is trivial, the constraints of Proposition \ref{prop:RCBraidingConstraintsFinal} are the same as in the split real case. The proof of this theorem is therefore the same as Theorem \ref{thm:split-class-sols} (without the requirement that $\sigma_3(1)$ is real).
\end{proof}
\begin{theorem}\label{thm:RCGNontrivialBijectionClassification}
	Braidings on $\mathcal{C}_{\mathbb{C}}(K_4^n \rtimes \mathbb{Z}/2\mathbb{Z}, \bar{\cdot}, \chi, \tau)$ are in bijection with pairs $(\sigma, \epsilon) \in \QF(\chi) \times \{\pm 1\}$ satisfying $$\sgn(\sigma|_{K_4^n})\sgn(\tau)\sigma(w) = 1.$$ 
\end{theorem}
\begin{proof}
We produce the data $(\sigma, \epsilon)$ in an identical way to the previous classification theorems. In this case, there is an extra constraint, namely that $\sigma_3$ is real, which holds if and only if $\sigma_3(1)$ is real. By Lemma \ref{lem:RCSigma3Squared1} and the definition of $\epsilon$, we have
$$\sigma_3(1)  = \epsilon \sqrt{2^{n + 1}\tau\sigma_1(w)\sgn(\sigma|_{K_4^n})},$$
which shows the constraint $\sgn(\sigma|_{K_4^n})\sgn(\tau)\sigma(w) = 1$ is necessary and sufficient for $\sigma_3$ to be real.
\end{proof}
\begin{notation}
    We denote a braiding on $\mathcal{C}(D_A, g ,\chi, \tau)$ by $\mathcal{C}_{\mathbb{C}, g}(\sigma, \epsilon)$. Note that $\tau$ is not necessarily determined by context, and the constraint $\sgn(\sigma|_{K_4^n})\sgn(\tau)\sigma(w)$ is also suppressed when $g$ is nontrivial. Moreover, we write $\sgn(\sigma) := \sgn(\sigma|_{K_4^n})$. No confusion should arise, since the sign of a quadratic form on $D_A$ is not defined. 
\end{notation}
The remainder of this section is dedicated to determining which of these braidings are equivalent, and some corollaries of this process. 

\begin{definition}
Let $f \in \Aut(D_A),~ \xi \in \Gal(\mathbb{C}/\mathbb{R})$ and $\lambda \in S^1$. We let $F(f,\xi,\lambda)$ be the candidate monoidal endofunctor of $\mathcal{C}_{\mathbb{C}}(A, g, \chi, \tau)$ whose underlying action on grouplike simples is $f$, fixes $m$ and applies $\xi$ to $\End(m)$. The tensorator coefficients are:

$$J_{a,b} = \id_{f(a)f(b)}, \quad J_{a,m} = \id_{f(a) \otimes m}, \quad J_{m,a} = \frac{\lambda}{\lambda^a}\id_m \otimes \id_{f(a)}, \quad J_{m,m} = \id_m \otimes \lambda \id_m.$$

We stress that in general, $F(f, \xi, \lambda)$ is not a monoidal functor. The consistency equations (simplified for our context from \cite[Theorem 6.10]{pss23}) are
\begin{align}
\chi\big(f(a), f(b)\big) &= \frac{\lambda \cdot \lambda^{ab}}{\lambda^a \cdot \lambda^b}\cdot \chi(a,b) \label{eqn:RCEndomorphismConsistency1}\\
\lambda^g &= \lambda. \label{eqn:RCEndomorphismConsistency2}
\end{align}
 Still, in the cases where  $F(f, \xi, \lambda)$ is monoidal, the composition rule can be seen to be 
$$F(f, \xi, \lambda) \circ  F(f', \xi', \lambda') \cong  F\big(f \circ f', \xi\circ \xi', \lambda \cdot \xi(\lambda')\big)$$
\end{definition}
\begin{remark}
The proof of \cite[Theorem 6.10]{pss23} shows that the functors $F(f, \xi, \lambda)$ satisfying the two consistency equations \eqref{eqn:RCEndomorphismConsistency1}, \eqref{eqn:RCEndomorphismConsistency2} are a complete set of representatives for $\pi_0\Aut_{\otimes}(\mathcal{C}_{\bb C}(D_A, g, \chi, \tau))$. 
\end{remark}
\begin{lemma} \label{lem:RCFunctorClassification}
We have $$\pi_0\Aut_{\otimes}\big(\mathcal{C}_{\bb C}(D_A, g, \chi, \tau)\big) \cong \Aut(D_A, \chi) \times K_4$$ whenever $\chi$ is real-valued. When $g$ is nontrivial, the functors $F(f, \xi, \pm 1)$ form a complete set of representatives. When $g$ is trivial, we instead take $F(f, \xi, 1)$ and $F(f, \xi, i)$ as representatives. 
\end{lemma}
\begin{proof}
 We first observe the function $f$ and automorphism $\xi$ are invariants of the underlying functor. We next extract the consistency equations from \cite[page 35]{pss23} for a monoidal equivalence $\mu \colon F(f,\xi, \lambda) \to F(f, \xi, \lambda')$. In the notation used in \textit{loc. cit.}, our assumptions are that $\theta, \theta',\varphi, \varphi'$ are identically 1. The consistency equations thus trivialize to: 
 \begin{align*}
 \mu_a&= \frac{\mu_m^a}{\mu_m} \\
\frac{\lambda'}{(\lambda')^a} &= \frac{\lambda}{\lambda^a} \\ 
\lambda' &= \frac{\mu_m^{ga}\mu_m}{\mu_a}\lambda
 \end{align*} 

We begin with the case when $g$ is nontrivial. In this case, the monoidal functor consistency equations  \eqref{eqn:RCEndomorphismConsistency1}, \eqref{eqn:RCEndomorphismConsistency2} imply $\lambda$ is real and $f \in \Aut(D_A, \chi)$. Substituting the first consistency equation for $\mu$ into the third (with $a = w$) shows that $F(f, \xi, 1)$ is not monoidally isomorphic to  $F(f, \xi, -1)$. 

When $g$ is trivial, we can set $a = b = w$ in \eqref{eqn:RCEndomorphismConsistency2} and use that $\chi(f(w), f(w)) = \chi(w,w) = 1$ (Corollary \ref{cor:RCSelfPairingis1}) to conclude $\lambda^4 = 1$. The second of the three consistency conditions implies that whether or not $\lambda$ is real is a monoidal invariant. It remains to show that  the two functors $F(f, \xi, \pm 1)$ are isomorphic, and likewise for $F(f, \xi, \pm i)$. This can be achieved by setting $\mu_m = i$ and then defining $\mu_a$ according to the first consistency equation. The last equation holds since $g$ is trivial. Equation \eqref{eqn:RCEndomorphismConsistency1}, together with the restrictions on $\lambda$ now implies $f \in \Aut(D_A, \chi)$.
\end{proof}

\begin{proposition} \label{prop:RCFunctorBraided}
The monoidal functor $F(f, \xi, \lambda)$ is a braided equivalence $\mathcal{C}_{\mathbb{C}, g}(\sigma, \epsilon) \to \mathcal{C}_{\mathbb{C}, g}(\sigma', \epsilon')$ if and only if $f \cdot \sigma|_{K_4^n} = \sigma'|_{K_4^n}$, and 
\begin{align}
	\sigma'(w) &= \lambda^2\sigma(w)\label{eqn:FinalRCBraidingSquare1}\\ 
	\sigma_3'(1) &= \sigma_3(1)^\xi. \label{eqn:FinalRCBraidingSquare2}
\end{align}
\end{proposition}
\begin{proof}
 The conditions for $F(f, \xi, \lambda)$ to be a braided equivalence  $\mathcal{C}_{\mathbb{C}, g}(\sigma, \epsilon)  \to \mathcal{C}_{\mathbb{C}, g}(\sigma', \epsilon')$ are: 
\begin{align}
	\chi\big(f(a), f(b)\big) &= \chi(a,b)^\xi \label{eqn:RCBraidingSquare1}\\ 
	\sigma_1'\big(f(a)\big) &= \frac{\lambda^a}{\lambda}\sigma_1(a)^\xi \label{eqn:RCBraidingSquare2}\\ 
	\sigma_2'\big(f(a)\big) &= \frac{\lambda}{\lambda^a}\sigma_2(a)^\xi \label{eqn:RCBraidingSquare3}\\ 
	\sigma_3'\big(f(a)\big) &= \sigma_3(a)^\xi. \label{eqn:RCBraidingSquare4}
\end{align}
The first of these equations always holds since $f \in \Aut(A, \chi)$. Additionally, since $f$ fixes $w$, $f$ must  take conjugating elements to conjugating elements. We may also assume $\lambda^4 = 1$, so that $\lambda/\lambda^a = \lambda^a/\lambda$. These facts allow the derivation of Equation \eqref{eqn:RCBraidingSquare3} from Equations \eqref{eqn:RCBraidingSquare2} and \eqref{eqn:RCBraidingSquare4}. Finally, using that $\sigma_{1}$ is real, we can drop the $\xi$ in \eqref{eqn:RCBraidingSquare2}, as well as prove that \eqref{eqn:RCBraidingSquare4} holds for all $a$ if and only if it holds at $1$, which is exactly  \eqref{eqn:FinalRCBraidingSquare2}. Evaluating \eqref{eqn:RCBraidingSquare2} on elements in $A$ gives $f \cdot \sigma = \sigma'$, and evaluating at $w$ gives  \eqref{eqn:FinalRCBraidingSquare1}. These conditions are indeed equivalent to \eqref{eqn:RCBraidingSquare2}, as
\[
\pushQED{\qed}
\sigma_1'\big(f(aw)\big) = \sigma_1'\big(f(a)\big)\sigma_1'(w) = \frac{\lambda}{\lambda^{aw}}\sigma_1(a)\sigma_1(w) = \frac{\lambda}{\lambda^{aw}}\sigma_1(aw).
\qedhere
\popQED
\]
\let\qed\relax
\end{proof}

As with the rest of this section, the case when $|g|=1$ is significantly easier since the structure constants are $g$ fixed. 
\begin{theorem} When $n > 0$, there are exactly three equivalence classes of braidings on
$\mathcal{C}_{\mathbb{C}}(K_4^n \rtimes \mathbb{Z}/2\mathbb{Z}, \id, \chi, \tau)$. When $n = 0$ and $\tau < 0$, there is a unique equivalence class, and when $n = 0$ and $\tau > 0$, there are precisely two. These braidings are distinguished as follows: 
\begin{itemize}
	\item The braidings $\mathcal{C}_\mathbb{C, \id}(\sigma, \epsilon)$ are all equivalent if $\sgn(\sigma) = -\sgn(\tau)$. 
	\item If $\sgn(\sigma) = \sgn(\tau)$, then there are exactly two equivalence classes of braidings, distinguished by $\epsilon$. 
\end{itemize}
\end{theorem}
\begin{proof}
First, observe that only one of the two distinguished cases can occur when $n = 0$. We begin with the first case. Suppose we are given $\mathcal{C}_\mathbb{C, \id}(\sigma, \epsilon)$ and $\mathcal{C}_\mathbb{C, \id}(\sigma', \epsilon)$ with $\sgn(\sigma) =\sgn(\sigma') = -\sgn(\tau)$. In this case $\sigma_3(1)$ and $\sigma_3'(1)$ are square roots of negative reals, and are thus purely imaginary. So, we can choose an $\xi \in \Gal(\mathbb{C}/\mathbb{R})$ such that $\sigma_3(1)^\xi = \sigma_3'(1)$. Moreover, we can also find a 4th root of unity $\lambda$ such that $\lambda^2\sigma(w) = \sigma'(w)$. Finally, since the restrictions of $\sigma$ and $\sigma'$  to $K_4^n$, have the same sign, they are orbit equivalent and thus there exists an $f \in \Aut(K_4^n, \chi|_{K_4^n})$ with $f \cdot \sigma = \sigma'$ on $K_4^n$. By Lemma \ref{lem:CanonicalW}, $f$ has a unique extension (also denoted $f$) to $\Aut(A, \chi)$. Then $F(f, h, \lambda)$ is a braided equivalence $\mathcal{C}_\mathbb{C, \id}(\sigma, \epsilon) \to \mathcal{C}_\mathbb{C, \id}(\sigma', \epsilon')$ by Proposition \ref{prop:RCFunctorBraided}. 

In the second case, the value $\sigma_3(1)$ is real and thus fixed by all braided functors, and thus $\epsilon$ is a braided invariant. It remains to show that the value of $\sigma(w)$ can be changed. We choose $\lambda$ with $\lambda^2\sigma(w) = \sigma'(w)$, and $f$ satisfying $f \cdot \sigma = \sigma'$ on $K_4^n$, extend $f$ to $A$, and deduce that $F(f, h, \lambda)$ is the desired equivalence using Proposition \ref{prop:RCFunctorBraided}.
\end{proof}
If we let  $(\sigma, \epsilon) =  (\sigma', \epsilon')$ in Proposition \ref{prop:RCFunctorBraided}, we conclude: 
\begin{corollary}
Suppose $\sgn(\sigma) = -\sgn(\tau)$. Then 
$$\pi_0\Aut_{\text{br}}\big(\mathcal{C}_\mathbb{C, \id}(\sigma, \epsilon)\big) \cong \GO_{\sgn(\sigma)}.$$
If $\sgn(\sigma) = \sgn(\tau)$, then 
$$\pi_0\Aut_{\text{br}}\big(\mathcal{C}_\mathbb{C, \id}(\sigma, \epsilon)\big) \cong \GO_{\sgn(\sigma)}\times \mathbb{Z}/2\mathbb{Z}.$$
\end{corollary}

\begin{theorem} When $n \geq 0$, there are exactly four equivalence classes of braidings on $\mathcal{C}_{\mathbb{C}}(K_4^n \rtimes \mathbb{Z}/2\mathbb{Z}, \bar \cdot, \chi, \tau)$. When $n = 0$, there are two. Two braidings $\mathcal{C}_{\mathbb{C}, \bar \cdot}(\sigma, \epsilon)$ and  $\mathcal{C}_{\mathbb{C}, \bar \cdot}(\sigma', \epsilon')$ are equivalent if and only if $\sgn(\sigma) = \sgn(\sigma')$ and $\epsilon = \epsilon'$. 
\end{theorem}
\begin{proof}
The ``only if'' direction follows from Proposition \ref{prop:RCFunctorBraided}, noting that in this case all $F(f, \xi, \lambda)$ have $\lambda^2 = 1$, and moreover that $\sigma_3(1)$ is real and so $\epsilon$ is fixed. Note that in this case the value $\sigma(w)$ is determined by the sign of $\sigma$ (restricted to $K_4^n)$ and so is automatically preserved. 

The functor required for the converse can be constructed from any $f$ such that $f \cdot \sigma = \sigma'$ as the monoidal functor $F(f, \id, 1)$, again by Proposition \ref{prop:RCFunctorBraided}.
\end{proof}
Again choosing $(\sigma, \epsilon) =  (\sigma', \epsilon')$ in Proposition \ref{prop:RCFunctorBraided}: 
\begin{corollary}
$$\pi_0\Aut_{\text{br}}\big(\mathcal{C}_{\mathbb{C}, \bar \cdot}(\sigma, \epsilon)\big) \cong \GO_{\sgn(\sigma)} \times K_4$$
\end{corollary}

\begin{lemma}
	There are exactly two families of twist morphisms for any $\mathcal{C}_{\mathbb{C}, \bar \cdot}(\sigma, \epsilon)$, corresponding to a sign $\rho \in \{\pm 1\}$. These twists are indeed ribbon structures (in the sense of \cite[Definition 8.10.1]{EGNO15}).
\end{lemma}

\section{Braidings on Split Complex Tambara-Yamagami Categories} \label{sec:SplitComplex}

In this section, we use the results of sections \ref{sec:QFAnalysis} and \ref{sec:SplitReal} to determine the number of braidings on split complex Tambara-Yamagami categories. While the classification in terms of equivalence classes of quadratic forms  was determined by Galindo (\cite{GALINDO_2022}) already, the precise number of equivalence classes was not. Moreover, most previous computations were done in the case when the rank of the underlying group is small. We show here that there are fewer equivalence classes of Tambara-Yamagami categories in these cases than in general. 

This process does not require any new computations. We begin by recalling the discussion of \cite[Section 2.5]{SchopierayNonDegenExtension}, which computes the number of equivalence classes of split complex Tambara-Yamagami categories with underlying group of rank $\leq 2$. 

Let $\ell$ be the nontrivial bicharacter on $\mathbb{Z}/2\mathbb{Z}$. There are two quadratic forms with coboundary $\ell$; these are inequivalent. Moreover, there are exactly three equivalence classes of quadratic forms on $K_4$ inducing $\ell^{2}$.
 Now let $\mathcal{C}_\mathbb{C}((\mathbb{Z}/2\mathbb{Z})^n, \chi, \tau)$ be a split complex Tambara-Yamagami category. Due to the fact that $\chi$ is symmetric, we can use the results of Wall \cite[Section 5]{wall63} to deduce that if $n$ is even, there are exactly two choices for $\chi$ and if $n$ is odd there is exactly one.
Indeed, when $n > 0$ is even, the representatives are $h^{ n/2}$ and $h^{(n - 2)/2} \oplus \ell ^{ 2}$. When $n$ is odd, the representative is $h^{ (n-1)/2} \oplus \ell$.

The following theorem both relies on, and strengthens the results of Galindo (\cite{GALINDO_2022}).

\begin{theorem}
Let $\mathcal{C}_\mathbb{C}((\mathbb{Z}/2\mathbb{Z})^n, \chi, \tau)$ be a split complex Tambara-Yamagami category ($\chi$ and $\tau$ are fixed). Then 

\begin{itemize}
	\item If $n > 0$ is even and $\chi \cong  h^{ n/2}$, there are exactly four equivalence classes of braidings on $\mathcal{C}_\mathbb{C}((\mathbb{Z}/2\mathbb{Z})^n, \chi, \tau)$. When $n = 0$, there are two. These are classified precisely by a free choice of a quadratic form $\sigma$ inducing $\chi$, together with a sign $\epsilon$. The formulas for the braidings are identical to Definition \ref{defn:ExplicitSplitRealBraidings}. These categories are symmetric if and only if they are defined over the reals, which occurs precisely when $\sgn(\sigma) = \sgn(\tau)$. Moreover, in this case 
	$$\pi_0\Aut_{\text{br}}\Big(\mathcal{C}_\mathbb{C}\big((\mathbb{Z}/2\mathbb{Z})^n, \chi, \tau, \sigma, \epsilon\big)\Big) \cong \GO_{\sgn \sigma}^{n / 2}.$$
	\item If $n \geq 4$ is even and $\chi \cong h^{(n - 2)/2} \oplus \ell^{ 2}$, there are exactly eight equivalence classes of braidings on  $\mathcal{C}_\mathbb{C}((\mathbb{Z}/2\mathbb{Z})^n, \chi, \tau)$. When $n = 2$, there are six. These are classified precisely by a free choice of a quadratic form $\zeta$ inducing $ h^{(n - 2)/2} \oplus \ell^{ 2}$, together with a sign $\epsilon$. These categories are never symmetric and are never defined over the reals. In this case, 
	$$\pi_0\Aut_{\text{br}}\big(\mathcal{C}_\mathbb{C}((\mathbb{Z}/2\mathbb{Z})^n, \chi, \tau, \zeta, \epsilon)\big) \cong \text{Stab}_{\Aut((\mathbb{Z}/2\mathbb{Z})^n, \chi)}(\zeta).$$
	\item If $n \geq 3$ is odd and $\chi \cong h^{ (n-1)/2} \oplus \ell$, there are exactly eight equivalence classes of braidings on $\mathcal{C}_\mathbb{C}((\mathbb{Z}/2\mathbb{Z})^n, \chi, \tau)$. If $n = 1$, then there are four. These are classified precisely by a free choice of a quadratic form $\sigma$ inducing $ h^{(n - 2)/2}$, a quadratic form $\nu$ inducing $\ell$, and a sign $\epsilon$. These categories are never symmetric and are never defined over the reals. In this case 
		$$\pi_0(\Aut_{\text{br}}(\mathcal{C}_\mathbb{C}((\mathbb{Z}/2\mathbb{Z})^n, \chi, \tau, \sigma,\nu ,\epsilon))) \cong \GO_{\sgn \sigma}^{(n - 1)/2}.$$
\end{itemize}
\end{theorem}
\begin{corollary}
A split complex braided Tambara-Yamagami category is symmetric if and only if it is defined over the reals. 
\end{corollary}
\begin{proof}
By \cite[Theorem 4.9]{GALINDO_2022}, we are reduced to calculating the number of orbits of quadratic forms inducing the three possible bicharacters, together with their stabilizers. We have already done this for $\chi = h^{ n}$ in Proposition \ref{prop:OrbitEquivalenceCharacterization} which gives most of the claims in this case. Indeed if $\chi = h^{ n}$ , the braiding coefficients $\sigma_1$ and $\sigma_2$ are always real. Thus, the braiding is symmetric if and only if the  function $\sigma_3(x) = \sigma_3(1)\sigma_1(x)$ is pointwise a sign. This occurs exactly when $\sigma_3(1)$ is real (so that the braiding is defined over the reals), which is again equivalent to $\sgn(\sigma) = \sgn(\tau)$. 

We tackle the case when $n$ is odd next. It is not too hard to see that extension by the identity of $\mathbb{Z}/2\mathbb{Z}$ gives an isomorphism 
$$ \Aut(K_4^{(n - 1)/2}, h^{ {(n - 1)/2}}) \cong \Aut(K_4^{(n - 1)/2} \times \mathbb{Z}/2\mathbb{Z}, h^{ {(n - 1)/2}} \oplus \ell).$$
In particular, the quadratic forms inducing $ h^{ {(n - 1)/2}} \oplus \ell$ decompose as products of quadratic forms on $K_4^{(n - 1)/2}$ and $\mathbb{Z}/2\mathbb{Z}$ inducing $h^{ {(n - 1)/2}}$ and $\ell$ respectively, and this decomposition is respected by $\Aut(K_4^{(n - 1)/2} \times \mathbb{Z}/2\mathbb{Z}, h^{ {(n - 1)/2}} \oplus \ell).$ This implies the results in the odd case, noting that any quadratic form inducing $\ell$ is complex valued and therefore not pointwise self-inverse. 

The last case is when the multiplicity of $\ell$ in $\chi$ is 2. This case follows from Proposition \ref{prop:StabilizerCombinatorics2ElectricBoogaloo} and the arguments above.

To conclude the statements about the groups of braided autoequivalences,  observe that Proposition \ref{prop:RealFunctorBraided} remains valid over the complex numbers, and all endofunctors of the split Tambara-Yamagami categories in question are still of the form $F(f)$. When the multiplicity of $\ell$ in $\chi$ is 2, the sign of $\sigma$ is not (in general) well defined and so we choose not to pursue a better description of its stabilizer.  
\end{proof}
\section{Non-existence of Braidings on Complex/Complex Tambara-Yamagami Categories }\label{sec:CrossedBraided}
In this section we explain why braidings are not possible in the complex/complex case, where the endomorphism algebra of every simple object is isomorphic to the complex numbers.
We quickly recall the relevant definitions and results from Sections 3 and 4 of \cite{pss23}:

Let $x$ be a simple object in a complex/complex Tambara-Yamagami category. Then there are two $\mathbb{R}$-linear embeddings of $\End(1)$ into $\End(x)$, given by (conjugation by) the left and right unitors respectively.  If these embeddings coincide, $x$ is said to be \textit{Galois trivial}, and \textit{Galois nontrivial} otherwise. We remark that Galois nontriviality has the graphical interpretation that passing a bead over a strand labeled by $x$ has the effect of complex conjugation. Section 3 of \cite{pss23} contains additional background on Galois nontriviality. 

With this definition in hand, \cite[Section 4]{pss23} establishes the following: 

\begin{fact}
A complex/complex Tambara-Yamagami category either contains a Galois nontrivial simple object, or is fusion over the complex numbers. In the former case, $m$ is the only Galois nontrivial simple object.  
\end{fact}

Consequently, we confine our attention to the case when $m$ is Galois nontrival and the ambient category is only fusion over the reals. \begin{lemma}\label{lem:noComplexComplexBraidings}
There are no braidings on any complex/complex Tambara-Yamagami category where $m$ is the only Galois nontrivial simple object. 
\end{lemma}
\begin{proof} 
Let $a$ be a Galois trivial simple object (such as the monoidal unit). By naturality of the braiding and Galois nontriviality of $m$, we have $$ic_{a, m} = 
\begin{tikzineqn}
	\AMBraidCrossing
    \node[smallbead] at (-.5, -.5) {$i$};
\end{tikzineqn}
=
\begin{tikzineqn}
	\AMBraidCrossing
	\node[smallbead] at (.5, .5) {$i$};
\end{tikzineqn}
=
\begin{tikzineqn}
	\AMBraidCrossing
		\node[smallbead] at (0, .75) {$i$};
\end{tikzineqn}
=
\begin{tikzineqn}
	\AMBraidCrossing
	\node[smallbead] at (-1, 0) {$\bar i$};
\end{tikzineqn}
=-ic_{a,m}
$$
which proves that the braiding is zero, a contradiction.
\end{proof}

\begin{remark}
    There is a more general structure known as a $G$-crossed braiding, which is akin to an intertwining of a braiding and group action.
    Initially the authors anticipated the possibility of such a more general structure on these complex/complex categories.
    However, it was pointed out to us by Rui Wen that a similar argument to the one above shows that these categories cannot even admit $G$-crossed braidings.
\end{remark}

\appendix

\section{Graphical Calculus Computation For 
\texorpdfstring{\eqref{RCHexagon8}}{Hexagon Equations}}
\label{sec:appendixComputationEx}

Here, we provide an example for the computations used
to derive the hexagon equations
\eqref{eqn:hexR1}-\eqref{eqn:hexR16}
and \eqref{RCHexagon1}-\eqref{RCHexagon16}.
We chose \eqref{RCHexagon8} as our example
since it is the most complicated.

Our goal is to derive an equation of coefficients
from the hexagon axiom on three $m$ strands:

\begin{equation}\label{eqn:hexagonAxiom}
\begin{tikzineqn}[scale=0.5]
    \begin{knot}[clip width=5]
        \strand[strand m] % bottom left
        (-2,0) node[below,font=\tiny] {$m$}
        to ++(0,2)
        to ++(4,2)
        to ++(0,2)
        node[above,font=\tiny] {$m$};
        \strand[strand m] % bottom mid
        (-1,0) node[below,font=\tiny] {$m$}
        to ++(0,1)
        to ++(2,1)
        to ++(-3,2.01) % tikz does not like parallel lines
        to ++(0,2)
        node[above,font=\tiny] {$m$};
        \strand[strand m] % bottom right
        (2,0) node[below,font=\tiny] {$m$}
        to ++(0,2)
        to ++(-3,2)
        to ++(2,1)
        to ++(0,1)
        node[above,font=\tiny] {$m$};
    \end{knot}
    \draw[dashed] (-3,2) to (3,2);
    \draw[dashed] (-3,4) to (3,4);
\end{tikzineqn}
\quad=\quad
\begin{tikzineqn}[scale=0.5]
    \begin{knot}[clip width=5]
        \strand[strand m] % bottom left
        (-2,0) node[below,font=\tiny] {$m$}
        to ++(0,1)
        to ++(4,4)
        to ++(0,1)
        node[above,font=\tiny] {$m$};
        \strand[strand m] % bottom mid
        (-1,0) node[below,font=\tiny] {$m$}
        to ++(0,1)
        to ++(-1,1)
        to ++(0,4)
        node[above,font=\tiny] {$m$};
        \strand[strand m] % bottom right
        (2,0) node[below,font=\tiny] {$m$}
        to ++(0,4)
        to ++(-1,1)
        to ++(0,1)
        node[above,font=\tiny] {$m$};
    \end{knot}
    \draw[dashed] (-3,3) to (3,3);
    % \draw[dashed] (-3,4) to (3,4);
    % \draw[dashed] (-3,1) to (3,1);
    % \draw[dashed] (-3,5) to (3,5);
\end{tikzineqn}
\end{equation}
where the dashed line denotes the associator
\[
    \alpha_{m,m,m}\colon (m\otimes m)\otimes m
    \to m\otimes (m\otimes m).
\]
By the Yoneda lemma, the hexagon axiom can
be equivalently encoded by precomposing
both sides of \eqref{eqn:hexagonAxiom}
on $\Hom(m\otimes( m\otimes m), c)$
for each
$c\in \mathcal{C}_{\mathbb{C}}(D_A,g,\chi,\tau)$.
Since we have
\[
    m\otimes( m\otimes m) \cong 
    m\otimes\left(2\cdot \bigoplus_{x\in D_A} x\right)
    \cong m^{\oplus 2|D_A|},
\]
by Schur's lemma the only non-trivial Hom space
that we need to consider is
\[
    \Hom(m\otimes(m\otimes m), m),
\]
which is of dimension $2|D_A|$ as a $\mathbb{C}$-bimodule
and has
\begin{equation}\label{eqn:mmmBasisVectors}
\left\{\TetraTransform{right}{m}{m}{x}{m}{m}{m}\right\}_{x\in D_A}
\end{equation}
as a generating set.
Hence, we reduce to composing each of the generators 
in \eqref{eqn:mmmBasisVectors} by both sides of
\eqref{eqn:hexagonAxiom}, i.e.
\begin{equation}\label{eqn:mmmHexagonPrecomp}
\begin{tikzineqn}[scale=0.5]
\draw[grid] (-2,-6) grid (2,4);
\begin{knot}[clip width=5]
\strand[strand m] (1,1) -- ++(-1,-1)
                        -- ++(0,-1)
                        -- ++(-1,-1)
                        -- ++(3,-2)
                        -- ++(0,-2);
\strand[strand m] (1,1) -- ++(1,-1)
                        -- ++(0,-2)
                        -- ++(-4,-2)
                        -- ++(0,-2);
\strand[strand m] (0,3.5) -- ++(0,-1.5)
                          -- ++(-2,-2)
                          -- ++(0,-2)
                          -- ++(3,-2.01)
                          -- ++(-1,-1)
                          -- ++(0,-1);                      
\strand[strand x] (0,2) -- (1,1);
\end{knot}
\node[node x] at (1,2) {$x$};
\node[node m] at (0,4) {$m$};
\node[node m] at (-2,-7) {$m$};
\node[node m] at (0,-7) {$m$};
\node[node m] at (2,-7) {$m$};
% \draw[dashed] (-3,-1) -- (3,-1);
% \draw[dashed] (-3,-5) -- (3,-5);
\draw[dashed] (-3,-1.5) -- (3,-1.5);
\draw[dashed] (-3,-4.5) -- (3,-4.5);
\end{tikzineqn}
\quad=\quad
\begin{tikzineqn}[scale=0.5]
\draw[grid] (-2,-6) grid (2,4);
\begin{knot}[clip width=5]
\strand[strand m] (1,1) -- (2,0)
                      -- (0,-2)
                      -- (0,-4)
                      -- (-2,-6);
\strand[strand m] (1,1) -- (0,0)
                      -- (2,-2)
                      -- (2,-6);
\strand[strand m] (0,3.5) -- (0,2)
                          -- (-2,0)
                          -- (-2,-4)
                          -- (0,-6);
\strand[strand x] (0,2) -- (1,1);
\end{knot}
\draw[dashed] (-3,-3) -- (3,-3);
\node[node x] at (1,2) {$x$};
\node[node m] at (0,4) {$m$};
\node[node m] at (-2,-7) {$m$};
\node[node m] at (0,-7) {$m$};
\node[node m] at (2,-7) {$m$};
\end{tikzineqn}
\end{equation}
which we will reduce using the braiding and associator
coefficients.
For convenience, we denote
\[
\gamma(x,y)_s:=\frac{\overline{s}^{gx}\tau}{\chi(x,y)^{gy}}
\]
for any $x,y\in D_A$ and any $s\in \{1,i\}$.
Recall from \cite[Section 6]{pss23}
that precomposing with the associator $\alpha_{m,m,m}$
yields the following sum
\[
\begin{tikzineqn}[scale=0.5]
    \draw[grid] (-2,-2) grid (2,4);
    \begin{knot}[clip width=5]
        \strand[strand m] (0,3.5) -- ++(0,-1.5)
        -- ++(-2,-2)
        -- ++(0,-2);
        \strand[strand m] (1,1) -- ++(-1,-1)
        -- ++(0,-2);
        \strand[strand m] (1,1) -- ++(1,-1)
        -- ++(0,-2);
        \strand[strand x] (0,2) -- (1,1);
    \end{knot}
    \node[node x] at (1,2) {$x$};
\node[node m] at (0,4) {$m$};
\node[node m] at (-2,-3) {$m$};
\node[node m] at (0,-3) {$m$};
\node[node m] at (2,-3) {$m$};
\draw[dashed] (-3,-1) -- (3,-1);
\end{tikzineqn}
\quad=\quad
\sum_{\substack{y\in D_A \\ s\in \{1,i\}}}
\begin{tikzineqn}[scale=0.5]
    \draw[grid] (-2,-2) grid (2,4);
    \begin{knot}[clip width=5]
        \strand[strand m] (0,3.5) -- ++(0,-1.5)
        -- ++(-2,-2)
        -- ++(0,-2);
        \strand[strand m] (-1,1) -- ++(1,-1)
        -- ++(0,-2);
        \strand[strand m] (0,2) -- ++(2,-2)
        -- ++(0,-2);
        \strand[strand y] (0,2) -- (-1,1);
    \end{knot}
\node[node y] at (-1,2) {$y$};
\node[node m] at (0,4) {$m$};
\node[node m] at (-2,-3) {$m$};
\node[node m] at (0,-3) {$m$};
\node[node m] at (2,-3) {$m$};
\node[emptybead] at (0,-1) {};
\node[beadLabel, right] at (0.1,-1) {$s$};
\node[emptybead] at (2,-1) {};
\node[beadLabel, right] at (2.1,-1) {$\gamma(x,y)_s\;\;.$};
\end{tikzineqn}
\]
This can also be derived directly by composing
the generating morphism from \eqref{eqn:mmmBasisVectors}
with the description of $\alpha_{m,m,m}$ 
in Theorem \ref{thm:RealComplexFromPSS}.
We summarize the moves that we
perform in the graphical calculus computations
as follows:
\begin{enumerate}
\item We may pass ``beads" (which represent
morphisms in $\End(m)\cong \mathbb{C}$)
through braidings and associators by naturality;
\item We may ``pull" trivalent vertices through
braidings by naturality;
\item We may pull out sums by bilinearity of composition;
\item We may pass beads through trivalent vertices
as long as they pick up conjugating factors
following \eqref{eqn:passThroughVertexConjugate}
and \eqref{eqn:passThroughVertexConjugateMM}.
\end{enumerate}

\vspace{1cm}
\begin{center}
(continued on next page)
\end{center}
\newpage

\noindent The reduction for the left hand side of
\eqref{eqn:mmmHexagonPrecomp} is as below,
where indexing variables $w,z$ vary over $D_A$
and $t,u$ varies over $\{1,i\}$:
\newcommand{\mmmStringDiagramScale}{0.6}
\newcommand{\mmmStringDiagramScaleL}{0.6}
\newcommand{\LHSStepOne}{
\begin{tikzineqn}[scale=\mmmStringDiagramScale]
\draw[grid] (-2,-6) grid (2,4);
\begin{knot}[clip width=5]
\strand[strand m] (1,1) -- ++(1,-1)
                        -- ++(0,-2)
                        -- ++(-4,-2)
                        -- ++(0,-2);
\strand[strand m] (0,3.5) -- ++(0,-1.5)
                          -- ++(-2,-2)
                          -- ++(0,-2)
                          -- ++(3,-2.01)
                          -- ++(-1,-1)
                          -- ++(0,-1);
\strand[strand m] (1,1) -- ++(-1,-1)
                        -- ++(0,-1)
                        -- ++(-1,-1)
                        -- ++(3,-2)
                        -- ++(0,-2);
\strand[strand x] (0,2) -- (1,1);
\end{knot}
\node[node x] at (1,2) {$x$};
\node[node m] at (0,4) {$m$};
\node[node m, below] at (-2,-6) {$m$};
\node[node m, below] at (0,-6) {$m$};
\node[node m, below] at (2,-6) {$m$};
\draw[dashed] (-3,-1.5) -- (3,-1.5);
\draw[dashed] (-3,-4.5) -- (3,-4.5);
\end{tikzineqn}
}

\newcommand{\LHSStepTwo}{
\displaystyle\sum_{z,t}\hspace{-0.2cm}
\begin{tikzineqn}[scale=\mmmStringDiagramScale]
\draw[grid] (-2,-6) grid (2,4);
\begin{knot}[clip width=5]
\strand[strand m] (0,2) -- ++(2,-2)
                        -- ++(0,-2)
                        -- ++(-4,-2)
                        -- ++(0,-2);
\strand[strand m] (-1,1) -- ++(1,-1)
                        -- ++(0,-1)
                        -- ++(-1,-1)
                        -- ++(3,-2)
                        -- ++(0,-2);
\strand[strand m] (0,3.5) -- ++(0,-1.5)
                          -- ++(-2,-2)
                          -- ++(0,-2)
                          -- ++(3,-2.01)
                          -- ++(-1,-1)
                          -- ++(0,-1);
\strand[strand z] (0,2) -- (-1,1);
\end{knot}
\node[node z] at (-1,2) {$z$};
\node[node m] at (0,4) {$m$};
\node[node m, below] at (-2,-6) {$m$};
\node[node m, below] at (0,-6) {$m$};
\node[node m, below] at (2,-6) {$m$};
\draw[dashed] (-3,-4.5) -- (3,-4.5);
\node[emptybead] at (-0.5,0.5) {};
\node[beadLabel, above right] at (-0.5,0.5) {$t$};
\node[emptybead] at (1,1) {};
\node[beadLabel, above right] at (1,1) {$\gamma(x,z)_t$};
\end{tikzineqn}
}

\newcommand{\LHSStepThree}{
\displaystyle\sum_{z,t}\hspace{-1cm}
\begin{tikzineqn}[scale=\mmmStringDiagramScale]
\draw[grid] (-2,-6) grid (2,4);
\begin{knot}[clip width=5]
\strand[strand m] (0,2) -- ++(2,-2)
                        -- ++(0,-2)
                        -- ++(-4,-2)
                        -- ++(0,-2);
\strand[strand m] (-1,1) -- ++(1,-1)
                        -- ++(0,-1)
                        -- ++(-1,-1)
                        -- ++(3,-2)
                        -- ++(0,-2);
\strand[strand m] (0,3.5) -- ++(0,-1.5)
                          -- ++(-2,-2)
                          -- ++(0,-2)
                          -- ++(3,-2.01)
                          -- ++(-1,-1)
                          -- ++(0,-1);
\strand[strand z] (0,2) -- (-1,1);
\end{knot}
\node[node z] at (-1,2) {$z$};
\node[node m] at (0,4) {$m$};
\node[node m, below] at (-2,-6) {$m$};
\node[node m, below] at (0,-6) {$m$};
\node[node m, below] at (2,-6) {$m$};
\draw[dashed] (-3,-4.5) -- (3,-4.5);
\node[emptybead] at (2,-5.5) {};
\node[beadLabel, right] at (2.1,-5.5) {$t$};
\node[emptybead] at (-2,-5.5) {};
\node[beadLabel, left] at (-2.1,-5.5)
    {$\gamma(x,z)_t$};
\end{tikzineqn}
}

\newcommand{\LHSStepFour}{
\displaystyle\sum_{z,t}\hspace{-1cm}
\begin{tikzineqn}[scale=\mmmStringDiagramScale]
\draw[grid] (-2,-6) grid (2,4);
\begin{knot}[clip width=5]
\strand[strand m] (0,2) -- ++(0,1.5);
\strand[strand m] (0,2) -- ++(2,-2)
                        -- ++(-4,-3)
                        -- ++(0,-3);
\strand[strand z] (0,2) to [bend right=0] ++(-1,-1)
                        to [bend right=0] ++(0,-1)
                        to [bend right=35] ++(1,-1.5)
                        to [bend left=35] ++(1.5,-1);
\strand[strand m] (1.5,-2.5) to ++(0.5,-0.5)
                             to ++(0,-3);
\strand[strand m] (1.5,-2.5) to ++(-0.5,-0.5)
                             to ++(0,-1)
                             to ++(-1,-1)
                             to ++(0,-1);
\end{knot}
\node[node z] at (-1,2) {$z$};
\node[node m] at (0,4) {$m$};
\node[node m, below] at (-2,-6) {$m$};
\node[node m, below] at (0,-6) {$m$};
\node[node m, below] at (2,-6) {$m$};
\draw[dashed] (-3,-4.5) -- (3,-4.5);
\node[emptybead] at (2,-5.5) {};
\node[beadLabel, right] at (2.1,-5.5) {$t$};
\node[emptybead] at (-2,-5.5) {};
\node[beadLabel, left] at (-2.1,-5.5)
    {$\gamma(x,z)_t$};
\end{tikzineqn}
}

\newcommand{\LHSStepFive}{
\displaystyle\sum_{z,t} \hspace{-0.5cm}
\begin{tikzineqn}[scale=\mmmStringDiagramScale]
\draw[grid] (-2,-3) grid (2,4);
\begin{knot}[clip width=5]
\strand[strand m] (0,3.5) to ++(0,-1.5)
                          to ++(-2,-2)
                          to ++(0,-3);
\strand[strand m] (1,1) to ++(-1,-1)
                        to ++(0,-3);
\strand[strand m] (1,1) to ++(1,-1)
                        to ++(0,-3);
\strand[strand z] (1,1) to ++(-1,1);
\end{knot}
\node[node z] at (1,2) {$z$};
\node[node m] at (0,4) {$m$};
\node[node m, below] at (-2,-3) {$m$};
\node[node m, below] at (0,-3) {$m$};
\node[node m, below] at (2,-3) {$m$};
\draw[dashed] (-3,-1) -- (3,-1);
\node[emptybead] at (2,-2) {};
\node[beadLabel, right] at (2.1,-2) {$t$};
\node[emptybead] at (-2,-2) {};
\node[beadLabel, below left] at (-2,-2)
    {$\gamma(x,z)_t$};
\node[emptybead] at (-1,1) {};
\node[beadLabel, above left] at (-1,1)
    {$\sigma_2(z)$};
\end{tikzineqn}
}

\newcommand{\LHSStepSix}{
\displaystyle\sum_{z,t} \hspace{-0.5cm}
\begin{tikzineqn}[scale=\mmmStringDiagramScale]
\draw[grid] (-2,-3) grid (2,4);
\begin{knot}[clip width=5]
\strand[strand m] (0,3.5) to ++(0,-1.5)
                          to ++(-2,-2)
                          to ++(0,-3);
\strand[strand m] (1,1) to ++(-1,-1)
                        to ++(0,-3);
\strand[strand m] (1,1) to ++(1,-1)
                        to ++(0,-3);
\strand[strand z] (1,1) to ++(-1,1);
\end{knot}
\node[node z] at (1,2) {$z$};
\node[node m] at (0,4) {$m$};
\node[node m, below] at (-2,-3) {$m$};
\node[node m, below] at (0,-3) {$m$};
\node[node m, below] at (2,-3) {$m$};
\draw[dashed] (-3,-0.5) -- (3,-0.5);
\node[emptybead] at (2,-2) {};
\node[beadLabel, right] at (2.1,-2) {$t$};
\node[emptybead] at (-2,-2.4) {};
\node[beadLabel, left] at (-2.1,-2.4)
    {$\gamma(x,z)_t$};
\node[emptybead] at (-2,-1.6) {};
\node[beadLabel, right] at (-1.9,-1.6)
    {$\sigma_2(z)$};
\end{tikzineqn}
}

\newcommand{\LHSStepSeven}{
\displaystyle\sum_{w,z,t,u} \hspace{-0.8cm}
\begin{tikzineqn}[scale=\mmmStringDiagramScale]
\draw[grid] (-2,-3) grid (2,4);
\begin{knot}[clip width=5]
\strand[strand m] (0,3.5) to ++(0,-1.5)
                          to ++(-2,-2)
                          to ++(0,-3);
\strand[strand m] (-1,1) to ++(1,-1)
                        to ++(0,-3);
\strand[strand m] (0,2) to ++(2,-2)
                        to ++(0,-3);
\strand[strand w] (0,2) to ++(-1,-1);
\end{knot}
\node[node w] at (-1,2) {$w$};
\node[node m] at (0,4) {$m$};
\node[node m, below] at (-2,-3) {$m$};
\node[node m, below] at (0,-3) {$m$};
\node[node m, below] at (2,-3) {$m$};
\node[emptybead] at (2,-2) {};
\node[beadLabel, right] at (2.1,-2) {$t$};
\node[emptybead] at (-2,-2) {};
\node[beadLabel, below left] at (-2,-2)
    {$\gamma(x,z)_t$};
\node[emptybead] at (-2,-1) {};
\node[beadLabel, right] at (-1.9,-1)
    {$\sigma_2(z)$};
\node[emptybead] at (-0.5,0.5) {};
\node[beadLabel, above right] at (-0.5,0.5) {$u$};
\node[emptybead] at (1,1) {};
\node[beadLabel, above right] at (1,1) {$\gamma(z,w)_u$};
\end{tikzineqn}
}

\newcommand{\LHSStepEight}{
\displaystyle\sum_{w,z,t,u}
\begin{tikzineqn}[scale=\mmmStringDiagramScale,
    xscale=1.1]
\draw[grid] (-2,-3) grid (2,4);
\begin{knot}[clip width=5]
\strand[strand m] (0,3.5) to ++(0,-1.5)
                          to ++(-2,-2)
                          to ++(0,-3);
\strand[strand m] (-1,1) to ++(1,-1)
                        to ++(0,-3);
\strand[strand m] (0,2) to ++(2,-2)
                        to ++(0,-3);
\strand[strand w] (0,2) to ++(-1,-1);
\end{knot}
\node[node w] at (-1,2) {$w$};
\node[node m] at (0,4) {$m$};
\node[node m, below] at (-2.5,-3) {$m$};
\node[node m, below] at (0,-3) {$m$};
\node[node m, below] at (2.5,-3) {$m$};
\node[emptybead] at (2,-2) {};
\node[beadLabel, right] at (2.1,-2) {$t$};
\node[emptybead] at (0,-2) {};
\node[beadLabel, below left] at (0.1,-2)
    {$\gamma(x,z)_t^{gw}$};
\node[emptybead] at (0,-1) {};
\node[beadLabel, right] at (0.1,-1)
    {$\sigma_2(z)^{gw}$};
\node[emptybead] at (-0.5,0.5) {};
\node[beadLabel, above right] at (-0.5,0.5) {$u$};
\node[emptybead] at (1,1) {};
\node[beadLabel, above right] at (1,1) {$\gamma(z,w)_u$};
\end{tikzineqn}
}

\newcommand{\LHSStepNine}{
\displaystyle\sum_{z,t,u}
\begin{tikzineqn}[scale=\mmmStringDiagramScale,
    xscale=1.1]
\draw[grid] (-2,-3) grid (2,4);
\begin{knot}[clip width=5]
\strand[strand m] (0,3.5) to ++(0,-1.5)
                          to ++(-2,-2)
                          to ++(0,-3);
\strand[strand m] (-1,1) to ++(1,-1)
                        to ++(0,-3);
\strand[strand m] (0,2) to ++(2,-2)
                        to ++(0,-3);
\strand[strand y] (0,2) to ++(-1,-1);
\end{knot}
\node[node y] at (-1,2) {$y$};
\node[node m] at (0,4) {$m$};
\node[node m, below] at (-2.5,-3) {$m$};
\node[node m, below] at (0,-3) {$m$};
\node[node m, below] at (2.5,-3) {$m$};
\node[emptybead] at (2,-2) {};
\node[beadLabel, right] at (2.1,-2) {$t$};
\node[emptybead] at (0,-2) {};
\node[beadLabel, below left] at (0.1,-2)
    {$\gamma(x,z)_t^{gy}$};
\node[emptybead] at (0,-1) {};
\node[beadLabel, right] at (0.1,-1)
    {$\sigma_2(z)^{gy}$};
\node[emptybead] at (-0.5,0.5) {};
\node[beadLabel, above right] at (-0.5,0.5) {$u$};
\node[emptybead] at (1,1) {};
\node[beadLabel, above right] at (1,1) {$\gamma(z,y)_u$};
\end{tikzineqn}
}

% ---------------------------------------------------------

\newcommand{\RHSStepOne}{
\begin{tikzineqn}[scale=\mmmStringDiagramScale]
\draw[grid] (-2,-6) grid (2,4);
\begin{knot}[clip width=5]
\strand[strand m] (1,1) -- (2,0)
                      -- (0,-2)
                      -- (0,-4)
                      -- (-2,-6);
\strand[strand m] (1,1) -- (0,0)
                      -- (2,-2)
                      -- (2,-6);
\strand[strand m] (0,3.5) -- (0,2)
                          -- (-2,0)
                          -- (-2,-4)
                          -- (0,-6);
\strand[strand x] (0,2) -- (1,1);
\end{knot}
\draw[dashed] (-3,-3) -- (3,-3);
\node[node x] at (1,2) {$x$};
\node[node m] at (0,4) {$m$};
\node[node m, below] at (-2,-6) {$m$};
\node[node m, below] at (0,-6) {$m$};
\node[node m, below] at (2,-6) {$m$};
\end{tikzineqn}
}

\newcommand{\RHSStepTwo}{
\begin{tikzineqn}[scale=\mmmStringDiagramScale]
\draw[grid] (-2,-6) grid (2,4);
\begin{knot}[clip width=5]
\strand[strand m] (1,1) -- (2,0)
                      -- (2,-6);
\strand[strand m] (1,1) -- (0,0)
                      -- (0,-4)
                      -- (-2,-6);
\strand[strand m] (0,3.5) -- (0,2)
                          -- (-2,0)
                          -- (-2,-4)
                          -- (0,-6);
\strand[strand x] (0,2) -- (1,1);
\end{knot}
\node[emptybead] at (2,-1.5) {};
\node[beadLabel] at (3,-2) {$\sigma_3(x)$};
\draw[dashed] (-3,-3) -- (3,-3);
\node[node x] at (1,2) {$x$};
\node[node m] at (0,4) {$m$};
\node[node m, below] at (-2,-6) {$m$};
\node[node m, below] at (0,-6) {$m$};
\node[node m, below] at (2,-6) {$m$};
\end{tikzineqn}}

\newcommand{\RHSStepThree}{
\begin{tikzineqn}[scale=\mmmStringDiagramScale]
\draw[grid] (-2,-6) grid (2,4);
\begin{knot}[clip width=5]
\strand[strand m] (1,1) -- (2,0)
                      -- (2,-6);
\strand[strand m] (1,1) -- (0,0)
                      -- (0,-4)
                      -- (-2,-6);
\strand[strand m] (0,3.5) -- (0,2)
                          -- (-2,0)
                          -- (-2,-4)
                          -- (0,-6);
\strand[strand x] (0,2) -- (1,1);
\end{knot}
\node[emptybead] at (2,-4.5) {};
\node[beadLabel] at (3,-5) {$\sigma_3(x)$};
% \node[smallbead] at (2,-1.5) {$\sigma_3(x)$};
\draw[dashed] (-3,-3) -- (3,-3);
\node[node x] at (1,2) {$x$};
\node[node m] at (0,4) {$m$};
\node[node m, below] at (-2,-6) {$m$};
\node[node m, below] at (0,-6) {$m$};
\node[node m, below] at (2,-6) {$m$};
\end{tikzineqn}}

\newcommand{\RHSStepFour}{
\begin{tikzineqn}[scale=\mmmStringDiagramScale]
\draw[grid] (-2,-4) grid (2,4);
\begin{knot}[clip width=5]
\strand[strand m] (0,3.5) -- ++(0,-1.5)
                          -- ++(-2,-2)
                          -- ++(0,-2)
                          -- ++(2,-2);
\strand[strand m] (1,1) -- ++(-1,-1)
                        -- ++(0,-2)
                        -- ++(-2,-2);
\strand[strand m] (1,1) -- ++(1,-1)
                        -- ++(0,-4);
\strand[strand x] (0,2) -- (1,1);
\end{knot}
\draw[dashed] (-3,-1) -- (3,-1);
\node[node x] at (1,2) {$x$};
\node[node m] at (0,4) {$m$};
\node[node m, below] at (-2,-6) {$m$};
\node[node m, below] at (0,-6) {$m$};
\node[node m, below] at (2,-6) {$m$};
\node[emptybead] at (2,-2.5) {};
\node[beadLabel] at (3,-3) {$\sigma_3(x)$};
\end{tikzineqn}}

\newcommand{\RHSStepFive}{
\displaystyle\sum_{y,s}
\begin{tikzineqn}[scale=\mmmStringDiagramScale]
\draw[grid] (-2,-4) grid (2,4);
\begin{knot}[clip width=5]
\strand[strand m] (0,2) -- ++(2,-2)
                        -- ++(0,-4);
\strand[strand m] (-1,1) -- ++(1,-1)
                        -- ++(0,-2)
                        -- ++(-2,-2);
\strand[strand m] (0,3.5) -- ++(0,-1.5)
                          -- ++(-2,-2)
                          -- ++(0,-2)
                          -- ++(2,-2);
\strand[strand y] (0,2) -- ++(-1,-1);
\end{knot}
\node[node y] at (-1,2) {$y$};
\node[node m] at (0,4) {$m$};
\node[node m, below] at (-2,-4) {$m$};
\node[node m, below] at (0,-4) {$m$};
\node[node m, below] at (2,-4) {$m$};
\node[emptybead] at (2,-2.5) {};
\node[beadLabel] at (3,-3) {$\sigma_3(x)$};
\node[emptybead] at (1,1) {};
\node[beadLabel, above right] at (1,1) {$\gamma(x,y)_s$};
\node[emptybead] at (0,0) {};
\node[beadLabel, below left] at (0,0) {$s$};
\end{tikzineqn}}

\newcommand{\RHSStepSix}{
\displaystyle\sum_{y,s}
\begin{tikzineqn}[scale=\mmmStringDiagramScale]
\draw[grid] (-2,-4) grid (2,4);
\begin{knot}[clip width=5]
\strand[strand m] (0,2) -- ++(2,-2)
                        -- ++(0,-4);
\strand[strand m] (-1,1) -- ++(1,-1)
                        -- ++(-2,-2)
                        -- ++(0,-2);
\strand[strand m] (0,3.5) -- ++(0,-1.5)
                          -- ++(-2,-2)
                          -- ++(2,-2)
                          -- ++(0,-2);
\strand[strand y] (0,2) -- ++(-1,-1);
\end{knot}
\node[node y] at (-1,2) {$y$};
\node[node m] at (0,4) {$m$};
\node[node m, below] at (-2,-4) {$m$};
\node[node m, below] at (0,-4) {$m$};
\node[node m, below] at (2,-4) {$m$};
\node[emptybead] at (2,-2.5) {};
\node[beadLabel, below right]
    at (2,-2.5) {$\sigma_3(x)$};
\node[emptybead] at (-2,-3) {};
\node[beadLabel, below right] at (-2,-3) {$s$};
\node[emptybead] at (1,1) {};
\node[beadLabel, above right] at (1,1) {$\gamma(x,y)_s$};
\end{tikzineqn}}

\newcommand{\RHSStepSeven}{
\displaystyle\sum_{y,s}
\begin{tikzineqn}[scale=\mmmStringDiagramScale]
\draw[grid] (-2,-4) grid (2,4);
\begin{knot}[clip width=5]
\strand[strand m] (0,3.5) -- ++(0,-1.5)
                          -- ++(-2,-2)
                          -- ++(0,-2)
                          -- ++(2,-2);
\strand[strand m] (-1,1) -- ++(1,-1)
                        -- ++(0,-2)
                        -- ++(-2,-2);
\strand[strand m] (0,2) -- ++(2,-2)
                        -- ++(0,-4);
\strand[strand y] (0,2) -- ++(-1,-1);
\end{knot}
\node[node y] at (-1,2) {$y$};
\node[node m] at (0,4) {$m$};
\node[node m, below] at (-2,-6) {$m$};
\node[node m, below] at (0,-6) {$m$};
\node[node m, below] at (2,-6) {$m$};
\node[emptybead] at (2,-2.5) {};
\node[beadLabel, below right]
    at (2,-2.5) {$\sigma_3(x)$};
\node[emptybead] at (-1.5,-3.5) {};
\node[beadLabel] at (-2,-3.5) {$s$};
\node[emptybead] at (1,1) {};
\node[beadLabel, above right] at (1,1) {$\gamma(x,y)_s$};
\end{tikzineqn}}

\newcommand{\RHSStepEight}{
\displaystyle\sum_{y,s}
\begin{tikzineqn}[scale=\mmmStringDiagramScale]
\draw[grid] (-2,-4) grid (2,4);
\begin{knot}[clip width=5]
\strand[strand m] (0,3.5) -- ++(0,-1.5)
                          -- ++(-2,-2)
                          -- ++(0,-2)
                          -- ++(0,-2);
\strand[strand m] (-1,1) -- ++(1,-1)
                        -- ++(0,-2)
                        -- ++(0,-2);
\strand[strand m] (0,2) -- ++(2,-2)
                        -- ++(0,-4);
\strand[strand y] (0,2) -- ++(-1,-1);
\end{knot}
\node[node y] at (-1,2) {$y$};
\node[node m] at (0,4) {$m$};
\node[node m, below] at (-2,-6) {$m$};
\node[node m, below] at (0,-6) {$m$};
\node[node m, below] at (2,-6) {$m$};
\node[emptybead] at (2,-2.5) {};
\node[beadLabel, below right]
    at (2,-2.5) {$\sigma_3(x)$};
\node[emptybead] at (-2,-3) {};
\node[beadLabel, below right]
    at (-2,-3) {$s$};
\node[emptybead] at (0,-3) {};
\node[beadLabel, below right]
    at (0,-3) {$\sigma_3(y)$};
\node[emptybead] at (1,1) {};
\node[beadLabel, above right] at (1,1) {$\gamma(x,y)_s$};
\end{tikzineqn}}

\newcommand{\RHSStepNine}{
\displaystyle\sum_{y,s}
\begin{tikzineqn}[scale=\mmmStringDiagramScale]
\draw[grid] (-2,-4) grid (2,4);
\begin{knot}[clip width=5]
\strand[strand m] (0,3.5) -- ++(0,-1.5)
                          -- ++(-2,-2)
                          -- ++(0,-2)
                          -- ++(0,-2);
\strand[strand m] (-1,1) -- ++(1,-1)
                        -- ++(0,-2)
                        -- ++(0,-2);
\strand[strand m] (0,2) -- ++(2,-2)
                        -- ++(0,-4);
\strand[strand y] (0,2) -- ++(-1,-1);
\end{knot}
\node[node y] at (-1,2) {$y$};
\node[node m] at (0,4) {$m$};
\node[node m, below] at (-2,-4) {$m$};
\node[node m, below] at (0,-4) {$m$};
\node[node m, below] at (2,-4) {$m$};
\node[emptybead] at (2,-2.5) {};
\node[beadLabel, below right]
    at (2,-2.5) {$\sigma_3(x)$};
\node[emptybead] at (-2,-3) {};
\node[beadLabel, below right]
    at (-2,-3) {$s$};
\node[emptybead] at (0,-3) {};
\node[beadLabel, below right]
    at (0,-3) {$\sigma_3(y)$};
\node[emptybead] at (1,1) {};
\node[beadLabel, above right] at (1,1) {$\gamma(x,y)_s$};
\end{tikzineqn}}

\newcommand{\RHSStepTenPrime}{
\displaystyle\sum_{y,s}
\begin{tikzineqn}[scale=\mmmStringDiagramScale]
\draw[grid] (-2,-4) grid (2,4);
\begin{knot}[clip width=5]
\strand[strand m] (0,3.5) -- ++(0,-1.5)
                          -- ++(-2,-2)
                          -- ++(0,-2)
                          -- ++(0,-2);
\strand[strand m] (-1,1) -- ++(1,-1)
                        -- ++(0,-2)
                        -- ++(0,-2);
\strand[strand m] (0,2) -- ++(2,-2)
                        -- ++(0,-4);
\strand[strand y] (0,2) -- ++(-1,-1);
\end{knot}
\node[node y] at (-1,2) {$y$};
\node[node m] at (0,4) {$m$};
\node[node m, below] at (-2,-6) {$m$};
\node[node m, below] at (0,-6) {$m$};
\node[node m, below] at (2,-6) {$m$};
\node[emptybead] at (2,-2.5) {};
\node[beadLabel, below right]
    at (2,-2.5) {$\sigma_3(x)$};
\node[emptybead] at (-0,-1) {};
\node[beadLabel, below left]
    at (-0,-1) {$s^{gy}$};
\node[emptybead] at (0,-3) {};
\node[beadLabel, below right]
    at (0,-3) {$\sigma_3(y)$};
\node[emptybead] at (1,1) {};
\node[beadLabel, above right] at (1,1) {$\gamma(x,y)_s$};
\end{tikzineqn}}

\newcommand{\RHSStepTen}{
\displaystyle\sum_{y,s}
\begin{tikzineqn}[scale=\mmmStringDiagramScaleL]
\draw[grid] (-2,-3) grid (2,4);
\begin{knot}[clip width=5]
\strand[strand m] (0,3.5) -- ++(0,-1.5)
                          -- ++(-2,-2)
                          -- ++(0,-3);
\strand[strand m] (-1,1) -- ++(1,-1)
                        -- ++(0,-3);
\strand[strand m] (0,2) -- ++(2,-2)
                        -- ++(0,-3);
\strand[strand y] (0,2) -- ++(-1,-1);
\end{knot}
\node[node y] at (-1,2) {$y$};
\node[node m] at (0,4) {$m$};
\node[node m, below] at (-2,-3) {$m$};
\node[node m, below] at (0,-3) {$m$};
\node[node m, below] at (2,-3) {$m$};
\node[emptybead] at (2,-1) {};
\node[beadLabel, below right]
    at (2,-1) {$\sigma_3(x)$};
\node[emptybead] at (0,-1) {};
\node[beadLabel, above right]
    at (0,-1) {$s^{gy}$};
\node[emptybead] at (0,-2) {};
\node[beadLabel, below left]
    at (0,-2) {$\sigma_3(y)$};
\node[emptybead] at (1,1) {};
\node[beadLabel, above right] at (1,1) {$\gamma(x,y)_s$};
\end{tikzineqn}}

\newcommand{\RHSStepEleven}{
\displaystyle\sum_{s}
\begin{tikzineqn}[scale=\mmmStringDiagramScaleL]
\draw[grid] (-2,-3) grid (2,4);
\begin{knot}[clip width=5]
\strand[strand m] (0,3.5) -- ++(0,-1.5)
                          -- ++(-2,-2)
                          -- ++(0,-3);
\strand[strand m] (-1,1) -- ++(1,-1)
                        -- ++(0,-3);
\strand[strand m] (0,2) -- ++(2,-2)
                        -- ++(0,-3);
\strand[strand y] (0,2) -- ++(-1,-1);
\end{knot}
\node[node y] at (-1,2) {$y$};
\node[node m] at (0,4) {$m$};
\node[node m, below] at (-2,-3) {$m$};
\node[node m, below] at (0,-3) {$m$};
\node[node m, below] at (2,-3) {$m$};
\node[emptybead] at (2,-1) {};
\node[beadLabel, below right]
    at (2,-1) {$\sigma_3(x)$};
\node[emptybead] at (0,-1) {};
\node[beadLabel, above right]
    at (0,-1) {$s^{gy}$};
\node[emptybead] at (0,-2) {};
\node[beadLabel, below left]
    at (0,-2) {$\sigma_3(y)$};
\node[emptybead] at (1,1) {};
\node[beadLabel, above right] at (1,1) {$\gamma(x,y)_s$};
\end{tikzineqn}}

\[
\setlength\arraycolsep{2pt}
\renewcommand\arraystretch{8}
\begin{array}{llllll}
& \LHSStepOne
& =
& \LHSStepTwo
& =
& \LHSStepThree \\
=
& \LHSStepFour
& =
& \LHSStepFive
& =
& \LHSStepSix \\
=
& \LHSStepSeven
& =
& \LHSStepEight
&&
\end{array}
\]

\newpage
\noindent The reduction for the right hand side of
\eqref{eqn:mmmHexagonPrecomp} is as follows,
where the indexing variable $y$ varies over $D_A$
and $s$ varies over $\{1,i\}$:
\vspace{0.5cm}
\[
\setlength\arraycolsep{2pt}
\renewcommand\arraystretch{8}
\begin{array}{llllll}
& \RHSStepOne
& =
& \RHSStepTwo
& =
& \RHSStepThree \\
=
& \RHSStepFive
& =
& \RHSStepSix
& =
& \RHSStepNine \\
&& =
& \RHSStepTen
&
\end{array}
\]

\newpage
\noindent
Combining the left and right hand side computations,
\eqref{eqn:mmmHexagonPrecomp}
is reduced to
\[
\LHSStepEight \quad=\quad \RHSStepTen \, .
\]
By linear independence of the generators
\eqref{eqn:mmmBasisVectors} (which 
can be checked by precomposing with splittings),
for any $y\in D_A$,
we may identify the $y$-summand of the left sum
with the $y$-summand of the right sum,
which gives us
\begin{equation}\label{eqn:mmmComputationResult}
    \LHSStepNine \quad=\quad \RHSStepEleven.
\end{equation}
Regarding the two rightmost strands as a
complex bimodule
$\mathbb{C}\otimes_{\mathbb{R}}\mathbb{C}$,
we can rewrite \eqref{eqn:mmmComputationResult}
as an inline formula
\begin{equation}\label{eqn:mmmComputationResultInline}
\sum_{\substack{z\in D_A\\ t,u\in \{1,i\}}}
\gamma(x,z)_t^{gy} \sigma_2(z)^{gy} u
\otimes
t \gamma(z,y)_u
= \sum_{s\in \{1,i\}}
\sigma_3(y) s^{gy}
\otimes \sigma_3(x) \gamma(x,y)_s .
\end{equation}
To extract an equation from
\eqref{eqn:mmmComputationResultInline},
we recall some algebraic facts about 
$\mathbb{C}\otimes_{\mathbb{R}} \mathbb{C}$.
There is a $\mathbb{C}$-algebra isomorphism
\begin{equation}\label{eqn:CtensorRCiso}
\phi\colon \mathbb{C}\otimes_{\mathbb{R}}\mathbb{C}
\xrightarrow[]{\cong } \mathbb{C}\oplus \mathbb{C},
\end{equation}
sending a simple tensor $\lambda\otimes \mu$ to
the pair $(\lambda\mu,\lambda\overline{\mu})$,
where the $\mathbb{C}$-action on
$\mathbb{C}\otimes_{\mathbb{R}}\mathbb{C}$
is given by acting on the left tensor factor.
The elements in
$\mathbb{C}\otimes_{\mathbb{R}} \mathbb{C}$
that correspond to the
$\mathbb{C}$-basis vectors
$(1,0)$ and $(0,1)$ in
$\mathbb{C}\oplus \mathbb{C}$
(forgetting the algebra structure) are
\[
    P := \tfrac{1}{2} (1\otimes 1 - i\otimes i)
    \text{ and }
    Q := \tfrac{1}{2} (1\otimes 1 + i\otimes i)
\]
respectively.
These basis elements satisfy the relations
\[
    P\cdot (\lambda\otimes 1) = P \cdot (1\otimes \lambda)
    \quad\text{and}\quad
    Q\cdot (\lambda\otimes 1) = 
    Q \cdot (1\otimes \overline{\lambda}),
\]
for any $\lambda\in \mathbb{C}$,
which is generalized slightly for our convenience
as follows:

\begin{lemma}[{\cite[Lemma 6.10]{pss23}}]
\label{lem:PCommutesWithScalars}
Let $u$ and $v$ represent words in the set $D_A\cup\{g\}$.  Consider the element
      \[P_{u,v}:=\tfrac12\left(1\otimes 1+i^u\otimes\overline{i}^v\right)
= \tfrac12 (1\otimes 1- (-1)^{|uv|} i\otimes i)
      \in\bb C\otimes_{\bb R}\bb C
\,.\]
      This element is an idempotent and it satisfies that $P_{u,v}\cdot(\lambda\otimes1)=P_{u,v}\cdot(1\otimes\lambda^{uv})$. 
\end{lemma}

Using the isomorphism \eqref{eqn:CtensorRCiso},
it is straightforward to check that:
\begin{lemma}\label{lem:PblankTensor1Injective}
For any $\lambda,\lambda'\in \mathbb{C}$,
we have
\[
    P_{u,v}(\lambda\otimes 1) = P_{u,v}(\lambda'\otimes 1)
    \quad\text{if and only if}\quad
    \lambda=\lambda'.
\]
\end{lemma}

We now begin our algebraic analysis of
\eqref{eqn:mmmComputationResultInline}.
For the right hand side,
we start by expanding $\gamma(x,y)_s$ to get
\[
\tau\sum_{s\in\{1,i\}} \sigma_3(y)s^{gy}\otimes
\sigma_3(x)\bar{s}^{gx}\chi(x,y)^{-gy}.
\]
The pair of conjugates $s$ and $\overline{s}$ 
on either side of the tensor product
allows us to identify a copy of
$P_{gy,gx}$, which gives us
\[
2\tau P_{gy,gx} \cdot
\big(\sigma_3(y)\otimes \sigma_3(x)\chi(x,y)^{-gy}\big),
\]
to which we apply
Lemma \ref{lem:PCommutesWithScalars} to get
\[
2\tau P_{1,xy} \cdot
\big(\sigma_3(y)
\sigma_3(x)^{xy}\chi(x,y)^{-gx}\big) 
\otimes 1.
\]

For the left hand side, we again start by
expanding the $\gamma$'s to get
\[
\tau^2 \sum_{\substack{z\in D_A\\ t,u\in \{1,i\}}}
\big(\overline{t}^{gx}\chi(x,z)^{-gz}\big)^{gy}
\sigma_2(z)^{gy}u
\otimes 
t \overline{u}^{gz}\chi(z,y)^{-gy}.
\]
The pair $u$ and $\overline{u}$ 
on either side of the tensor product
gives us a copy of $P_{1,gz}$,
so the equation becomes
\[
2 \tau^2 \sum_{\substack{z\in D_A\\ t\in \{1,i\}}} P_{1,gz} \cdot 
\big(\overline{t}^{xy} \chi(x,z)^{-yz}\big)
\sigma_2(z)^{gy}
\otimes
\big(t\chi(z,y)^{-gy}\big),
\]
to which we apply Lemma \ref{lem:PCommutesWithScalars}
to get
\[
2\tau^2
\sum_{\substack{z\in D_A\\ t\in \{1,i\}}}P_{1,gz} \cdot 
\big(\overline{t}^{xy} \chi(x,z)^{-yz}\big)
\big(t^{gz}\chi(z,y)^{-yz}\big) 
\sigma_2(z)^{gy}
\otimes 1,
\]
where we can collect the $t$ terms to obtain
\[
2\tau^2
\sum_{z\in D_A}
P_{1,gz} \cdot 
\left(\sum_{t\in \{1,i\}}
\overline{t}^{xy}t^{gz}\right)
\chi(x,z)^{-yz}
\chi(z,y)^{-yz}
\sigma_2(z)^{gy}
\otimes 1.
\]
Here, we may reduce the sum indexed by $t$ to
\[
\sum_{t\in \{1,i\}} \overline{t}^{xy}t^{gz}
= 1 + (-i)^{xy}\cdot i^{gz}
= 1 - (-1)^{|gxyz|}\cdot i^2
= \begin{cases}
   2, \text{ if } |gxyz|=0, \\
   0, \text{ if } |gxyz|=1.
\end{cases}
\]
Hence, the larger sum reduces to
\[
4\tau^2
\sum_{\substack{z\in D_A \\ |gxyz|=0}}
P_{1,gz}\cdot 
\chi(x,z)^{-yz}
\chi(z,y)^{-yz}
\sigma_2(z)^{gy}
\otimes 1,
\]
and we can use the condition on $z$ to
equate $|gz|=|xy|$ to obtain
\[
4\tau^2 P_{1,xy}\cdot 
\sum_{\substack{z\in D_A \\ |gxyz|=0}}
\chi(x,z)^{-yz}
\chi(z,y)^{-yz}
\sigma_2(z)^{gy}
\otimes 1.
\]
Thus, the condition \eqref{eqn:mmmHexagonPrecomp}
is equivalent to
\[
2\tau P_{1,xy} \cdot
\big(\sigma_3(y)
\sigma_3(x)^{xy}\chi(x,y)^{-gx}\big) 
\otimes 1 \\
=\ 4\tau^2 P_{1,xy}\cdot 
\sum_{\substack{z\in D_A \\ |gxyz|=0}}
\chi(x,z)^{-yz}
\chi(z,y)^{-yz}
\sigma_2(z)^{gy}
\otimes 1,
\]
and by Lemma \ref{lem:PblankTensor1Injective},
this reduces to
\[
\sigma_3(y) \sigma_3(x)^{xy}\chi(x,y)^{-gx} 
= 2\tau\cdot\sum_{\substack{z\in D_A \\ |gxyz|=0}}
\chi(x,z)^{-yz}
\chi(z,y)^{-yz}
\sigma_2(z)^{gy},
\]
which becomes \eqref{RCHexagon8} after applying
$(-)^{x}$.

~\\~\\
\noindent\textbf{Conflict of Interest:} The authors have no conflicts of interest to declare that are relevant to this article.

\noindent\textbf{Statement of Data Availability:} Data sharing not applicable to this article as no datasets were generated or analysed during the current study.


\newpage
\begin{thebibliography}{DMNO13}
	
	\bibitem[CHS18]{10.21468/SciPostPhys.5.1.006}
	Clay Córdova, Po-Shen Hsin, and Nathan Seiberg.
	\newblock {Time-reversal symmetry, anomalies, and dualities in (2+1)$d$}.
	\newblock {\em SciPost Phys.}, 5:006, 2018.
	
	\bibitem[DGNO10]{MR2609644}
	Vladimir Drinfeld, Shlomo Gelaki, Dmitri Nikshych, and Victor Ostrik.
	\newblock On braided fusion categories. {I}.
	\newblock {\em Selecta Math. (N.S.)}, 16(1):1--119, 2010.
	
	\bibitem[Dic58]{MR104735}
	Leonard~Eugene Dickson.
	\newblock {\em Linear groups: {W}ith an exposition of the {G}alois field
		theory}.
	\newblock Dover Publications, Inc., New York, 1958.
	\newblock With an introduction by W. Magnus.
	
	\bibitem[Die71]{MR310083}
	Jean~A. Dieudonn\'e.
	\newblock {\em La g\'eom\'etrie des groupes classiques}, volume Band 5 of {\em
		Ergebnisse der Mathematik und ihrer Grenzgebiete [Results in Mathematics and
		Related Areas]}.
	\newblock Springer-Verlag, Berlin-New York, 1971.
	\newblock Troisi\`eme \'edition.
	
	\bibitem[DMNO13]{MR3039775}
	Alexei Davydov, Michael M\"uger, Dmitri Nikshych, and Victor Ostrik.
	\newblock The {W}itt group of non-degenerate braided fusion categories.
	\newblock {\em J. Reine Angew. Math.}, 677:135--177, 2013.
	
	\bibitem[DNO13]{MR3022755}
	Alexei Davydov, Dmitri Nikshych, and Victor Ostrik.
	\newblock On the structure of the {W}itt group of braided fusion categories.
	\newblock {\em Selecta Math. (N.S.)}, 19(1):237--269, 2013.
	
	\bibitem[EG12]{MR2946231}
	Pavel Etingof and Shlomo Gelaki.
	\newblock Descent and forms of tensor categories.
	\newblock {\em Int. Math. Res. Not. IMRN}, (13):3040--3063, 2012.
	
	\bibitem[EGNO15]{EGNO15}
	Pavel Etingof, Shlomo Gelaki, Dmitri Nikshych, and Victor Ostrik.
	\newblock {\em Tensor categories}, volume 205 of {\em Mathematical Surveys and
		Monographs}.
	\newblock American Mathematical Society, Providence, RI, 2015.
	\newblock \mathscinet{MR3242743} \doi{10.1090/surv/205}.
	
	\bibitem[EM22]{EDIEMICHELL2022108364}
	Cain Edie-Michell.
	\newblock Auto-equivalences of the modular tensor categories of type a, b, c
	and g.
	\newblock {\em Advances in Mathematics}, 402:108364, 2022.
	
	\bibitem[ENO10]{MR2677836}
	Pavel Etingof, Dmitri Nikshych, and Victor Ostrik.
	\newblock Fusion categories and homotopy theory.
	\newblock {\em Quantum Topol.}, 1(3):209--273, 2010.
	\newblock With an appendix by Ehud Meir, \mathscinet{2677836}
	\doi{10.4171/QT/6} \arxiv{0909.3140}.
	
	\bibitem[FH21]{MR4268163}
	Daniel~S. Freed and Michael~J. Hopkins.
	\newblock Reflection positivity and invertible topological phases.
	\newblock {\em Geom. Topol.}, 25(3):1165--1330, 2021.
	
	\bibitem[Gal22]{GALINDO_2022}
	César Galindo.
	\newblock Trivializing group actions on braided crossed tensor categories and
	graded braided tensor categories.
	\newblock {\em Journal of the Mathematical Society of Japan}, 74(3), July 2022.
	
	\bibitem[GLM24]{galindo2024modular}
	César Galindo, Simon Lentner, and Sven Möller.
	\newblock Modular $\mathbb{Z}_2$-crossed tambara-yamagami-like categories for
	even groups, 2024.
	
	\bibitem[IT21]{Izumi_2021}
	Masaki Izumi and Henry Tucker.
	\newblock Algebraic realization of noncommutative near-group fusion categories.
	\newblock {\em Algebra \& Number Theory}, 15(5):1077–1093, June 2021.
	
	\bibitem[Jac09]{jacobson2009basic}
	N.~Jacobson.
	\newblock {\em Basic Algebra I: Second Edition}.
	\newblock Basic Algebra. Dover Publications, 2009.
	
	\bibitem[JS93]{MR1250465}
	Andr\'e{} Joyal and Ross Street.
	\newblock Braided tensor categories.
	\newblock {\em Adv. Math.}, 102(1):20--78, 1993.
	
	\bibitem[KW14]{kong2014braidedfusioncategoriesgravitational}
	Liang Kong and Xiao-Gang Wen.
	\newblock Braided fusion categories, gravitational anomalies, and the
	mathematical framework for topological orders in any dimensions, 2014.
	
	\bibitem[MH73]{MR506372}
	John Milnor and Dale Husemoller.
	\newblock {\em Symmetric bilinear forms}, volume Band 73 of {\em Ergebnisse der
		Mathematik und ihrer Grenzgebiete [Results in Mathematics and Related
		Areas]}.
	\newblock Springer-Verlag, New York-Heidelberg, 1973.
	
	\bibitem[Nik07]{Nikshych2007NongrouptheoreticalSH}
	Dmitri Nikshych.
	\newblock Non-group-theoretical semisimple hopf algebras from group actions on
	fusion categories.
	\newblock {\em Selecta Mathematica}, 14:145--161, 2007.
	
	\bibitem[PSS23]{pss23}
	Julia Plavnik, Sean Sanford, and Dalton Sconce.
	\newblock Tambara-yamagami categories over the reals: The non-split case, 2023.
	
	\bibitem[Rob17]{MR3675715}
	Bryan~W. Roberts.
	\newblock Three myths about time reversal in quantum theory.
	\newblock {\em Philos. Sci.}, 84(2):315--334, 2017.
	
	\bibitem[San25]{MR4806973}
	Sean Sanford.
	\newblock Fusion categories over non-algebraically closed fields.
	\newblock {\em J. Algebra}, 663:316--351, 2025.
	
	\bibitem[Sch85]{MR770063}
	Winfried Scharlau.
	\newblock {\em Quadratic and {H}ermitian forms}, volume 270 of {\em Grundlehren
		der mathematischen Wissenschaften [Fundamental Principles of Mathematical
		Sciences]}.
	\newblock Springer-Verlag, Berlin, 1985.
	
	\bibitem[Sch00]{MR1803370}
	Winfried Scharlau.
	\newblock On the history of the algebraic theory of quadratic forms.
	\newblock In {\em Quadratic forms and their applications ({D}ublin, 1999)},
	volume 272 of {\em Contemp. Math.}, pages 229--259. Amer. Math. Soc.,
	Providence, RI, 2000.
	
	\bibitem[Sch23]{SchopierayNonDegenExtension}
	Andrew Schopieray.
	\newblock Nondegenerate extensions of near-group braided fusion categories.
	\newblock {\em Revista de la Unión Matemática Argentina}, pages 413--438, 07
	2023.
	
	\bibitem[Sie00]{sie00}
	Jacob~A. Siehler.
	\newblock Braided near-group categories, 2000.
	
	\bibitem[Tam00]{Tambara2000}
	D.~Tambara.
	\newblock Representations of tensor categories with fusion rules of
	self-duality for abelian groups.
	\newblock {\em Israel Journal of Mathematics}, 118(1):29--60, 2000.
	
	\bibitem[TY98]{ty98}
	Daisuke Tambara and Shigeru Yamagami.
	\newblock Tensor categories with fusion rules of self-duality for finite
	abelian groups.
	\newblock {\em Journal of Algebra}, 209(2):692--707, 1998.
	
	\bibitem[Wal63]{wall63}
	C.~T.~C. Wall.
	\newblock Quadratic forms on finite groups, and related topics.
	\newblock {\em Topology}, 2:281--298, 1963.
	
\end{thebibliography}
\end{document}